\documentclass[a4paper]{article}
\usepackage[utf8]{inputenc}
\usepackage[a4paper,left=3cm,right=3cm,top=2cm,bottom=3cm,bindingoffset=5mm]{geometry}
\usepackage[english]{babel}
\usepackage[T1]{fontenc}
\usepackage{graphicx, subfig}
\graphicspath{{img/}}
\usepackage{fancyhdr}
\usepackage{lmodern}
\usepackage{color}
\usepackage{tocbibind}

\usepackage{amsfonts}
\usepackage{amsmath}
\usepackage{amsthm}
\usepackage{hyperref}
\title{Circularity and Symmetries of $p$ and $p^{2}$-polygons}
\author{Rolf Haag \\
\href{rhaag.98@gmail.com}{rhaag.98@gmail.com}}
\date{2025/26}
\usepackage{float}
\begin{document}
\maketitle
\begin{abstract}
The symmetry of polygons can be characterized by the number of symmetry axes they have. For $n$-polygons with $p$ or $p^2$ vertices $p\geq3$ there exist few symmetry categories, depending from the number of symmetry-axes the have. Further in the case of $n=p^2$ there exist $p^2$-polygons with no axis, but the property of $p$-circularity. This work investigates the corresponding equivalence classes and their enumeration. The total number of equivalence classes and the number of the regular ones are known. Formulas get established for the number of equivalence classes for the other categories of symmetry. We show complete sets of representatifs in some cases. Together, these results provide a comprehensive description of the symmetry structure of polygons with a prime number or prime square number of vertices and lay the groundwork for extending the classification and enumeration to polygons with a more composed numbers of vertices.
\end{abstract}

\tableofcontents

\part{Commons}

\section{Introduction}
\label{sec:introduction}

Polygons are among the most fundamental objects of Euclidean geometry and have been studied since antiquity. Beyond their geometric significance, they have also appeared in cultural and symbolic contexts, for example, the pentagram was regarded by the Pythagoreans as a figure of harmony and healing.

From a mathematical perspective, classical investigations have focused on properties such as angles, diagonals and area. Euler, using the totient function, derived and proved a formula for the number of distinct regular polygons that can be formed by placing $n$ equidistant vertices on a circle. In the twentieth century, Golomb and Welch \cite{Golomb1960} extended the problem to polygons in general, where the vertices need not to form regular shapes. They considered polygons equivalent if they differed only by a rotation in the plane, and similar if they could also be related by reflection through an axis. Using a combinatorial and group-theoretic approach, based on the action of the dihedral group $D_n$ and Burnside’s Lemma, they obtained formulas for the number of distinct equivalence-classes of $n$-polygons.\cite{Golomb1960} Their work thus provided a complete enumeration of fundamentally distinct $n$-polygons. However, their treatment did not address the degree of symmetry of a polygon: all symmetric and asymmetric cases were treated uniformly, without finer classification.

The study of symmetry degrees in polygons is of intrinsic interest for several reasons: Symmetry is a central organizing principle in mathematics, with implications ranging from group theory to crystallography. In polygons, the degree of symmetry can be defined in terms of the size of the subgroup of the dihedral group $D_n$,  that leaves the polygon invariant. At one extreme, regular polygons exhibit maximal symmetry, while at the other extreme, generic polygons are completely asymmetric, invariant only under the identity transformation. Between these extremes lie intermediate cases of partial symmetry. A systematic classification of polygons according to these symmetry degrees enriches our understanding of the combinatorial structure of polygons and provides insight into the balance between order and disorder in geometric configurations.

Earlier works, \cite{Haag2019},\cite{Haag2019a} provided investigations of two cases and discovered an unexpected connection to perfect numbers. These case studies suggested deeper structural relationships but did not yield a general framework for classifying polygons by symmetry. The aim of the present paper is to develop such a framework for polygons with $p$ and with $p^{2}$ vertices, where $p \geq 3$ is a prime number. The restriction to prime numbers is natural, since the symmetry structure is considerably simpler in the prime case, while composite values of $n$  introduce additional complications. This allows for a clear and complete classification, which can then serve as a foundation for the more general case of arbitrary $n$.

\section{Definition of a $n$-polygon}
\label{sec:definition_of_a_n_-polygon}

Let be $n\geq 5$ a natural number. $n$ vertices are regularly distributed on a circle. We consider the Hamiltonian cycles through the $n$ vertices. In this paper such Hamiltonian cycles are called $n$-polygons. The usual polygons are the special case where all edges have minimal length.

$ S ^ 1 \subset \mathbb{R} ^ 2 = \mathbb{C} $ is the unit circle in the Euclidean plane. The finite subset
\begin{center}
$V_n:=\{v_k:=e^{2 \pi i k/n} \mid k = 0,1,\ldots,n-1 \} \subset S^1$
\end{center}
represents the vertices of a $n$-polygon.
To describe the $n$ polygon $P(\sigma)$ we use the $n$ cycles $\sigma = (\sigma_1, \sigma_2, \cdots, \sigma_n)$ consisting of the $n$ numbers $\lbrace 0,1, \cdots , n-1 \rbrace$ in any order. The associated $n$-polygon $P(\sigma)$ is given by the path $\overline{v_{\sigma_1}v_{\sigma_2} \cdots v_{\sigma_n}v_{\sigma_1}}$ or more precisely by combining the links $\overline{v_{\sigma_i}v_{\sigma_{i+1}}}$, $i= 1,2, \cdots , n$, where $\sigma_{n+1}= \sigma_1$.
Each of the $n$ edges is assigned its length $e_i$.  $e_i = 1$ means that the $i-$th edge runs counterclockwise from the vertex $V_i$ to the following vertex $V_{i + 1}$. $e_i = 2$ means that the $i-$th edge runs counterclockwise from the vertex $V_i$  to the vertex $V_{i + 2}$ and so on. $e_i = n$ is not possible, since this would mean the connection of the vertex $V_i$ to itself. Therefore, only the numbers between $1$ and $n-1$ are allowed to describe the length of the edges. The length of an edge $e_i$ is referred to briefly as a side of the n-polygon. Therefore, an $n-$polygon can also get described by the $n$-cycle of its sides: $(e_1, e_2, \ldots, e_i, \ldots, e_n)$.\\\\
\newpage
\textbf{Example of a $9$-polygon with one symmetry axis}
\begin{figure}[H]
\begin{center}
\includegraphics[width=0.4\textwidth]{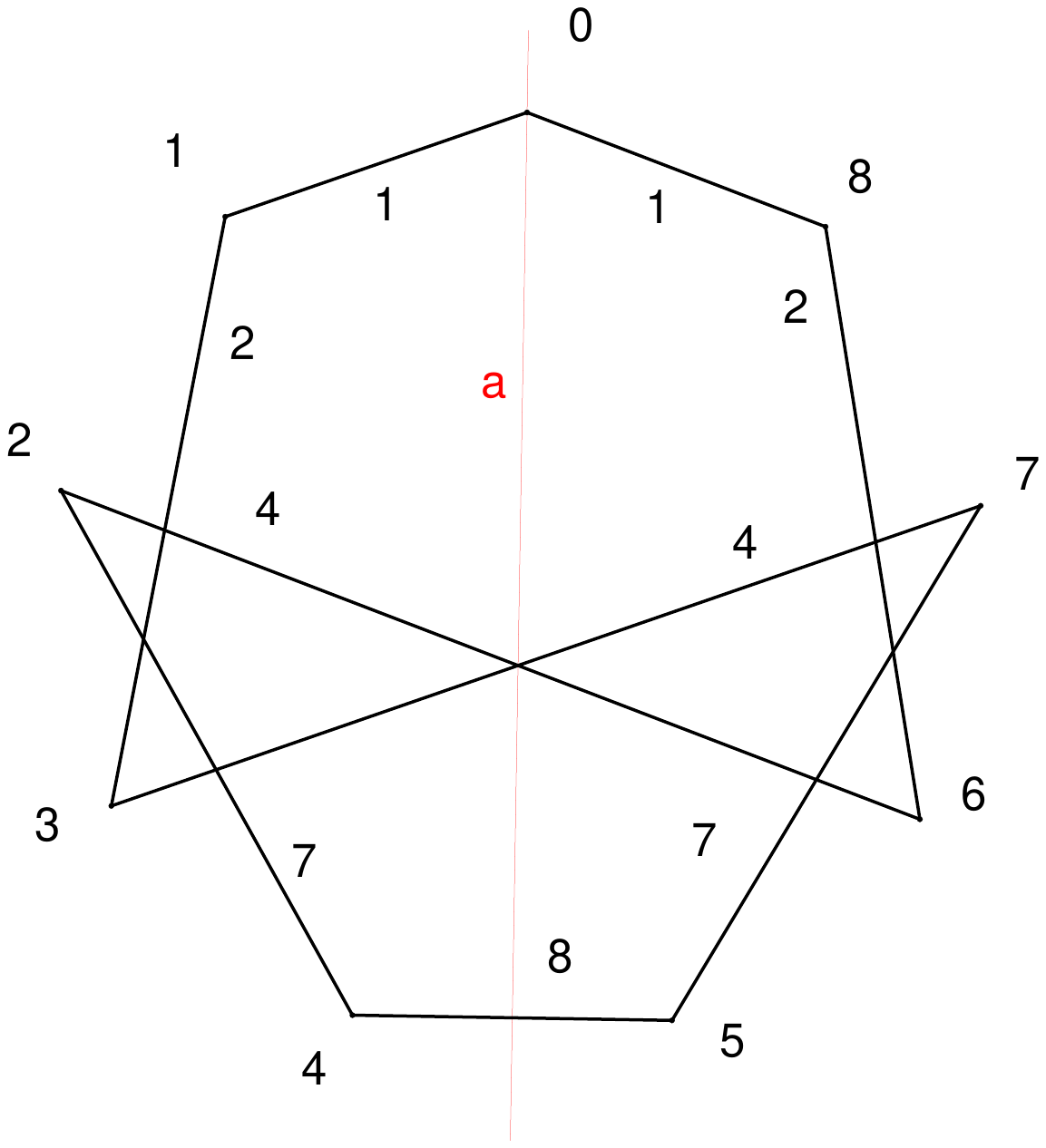}
\caption{Cycle of vertices (0 1 3 7 5 4 2 6 8) and cycle of edges (1 2 4 7 8 7 4 2 1)}
\end{center}
\end{figure}

\section{Definition of a $p$-circular $p^{2}$-polygon}
\label{sec: Definition_of_a_p-circular_p^2-polygon}
Für eine genaue Klassifizierung der $p^{2}$-polygone benötigen wir die nachstehende Definition:\\\\
A $p$-circular $p^{2}$-polygon is a $p^{2}$-polygon with the two properties:
\begin{enumerate}
\item{it has no symmetry-axes},
\item{After each rotation by an angle of $\dfrac{ 360}{p}$ it is coincident to itself.}
\end{enumerate}
\newpage
\begin{center}
\textbf{Two examples of a $p$-circular $p^{2}$-polygon}
\end{center}
\begin{figure} [H]
\begin{center}
\includegraphics[width=0.4\textwidth]{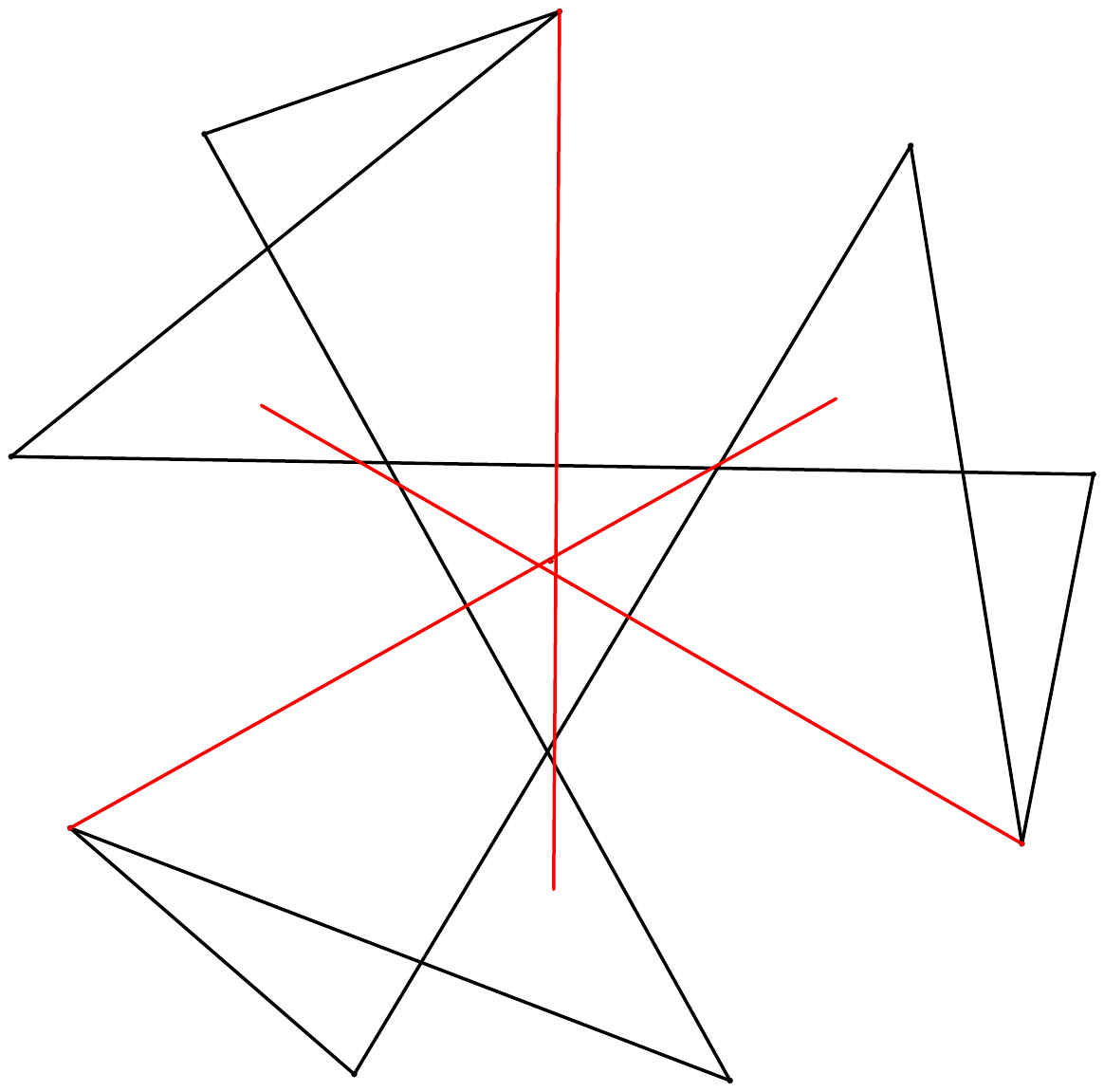}
\caption{Example of a 3-circular polygon with 9 vertices}
\end{center}
\end{figure}

\begin{figure}[H]
\begin{center}
\includegraphics[width=0.4\textwidth]{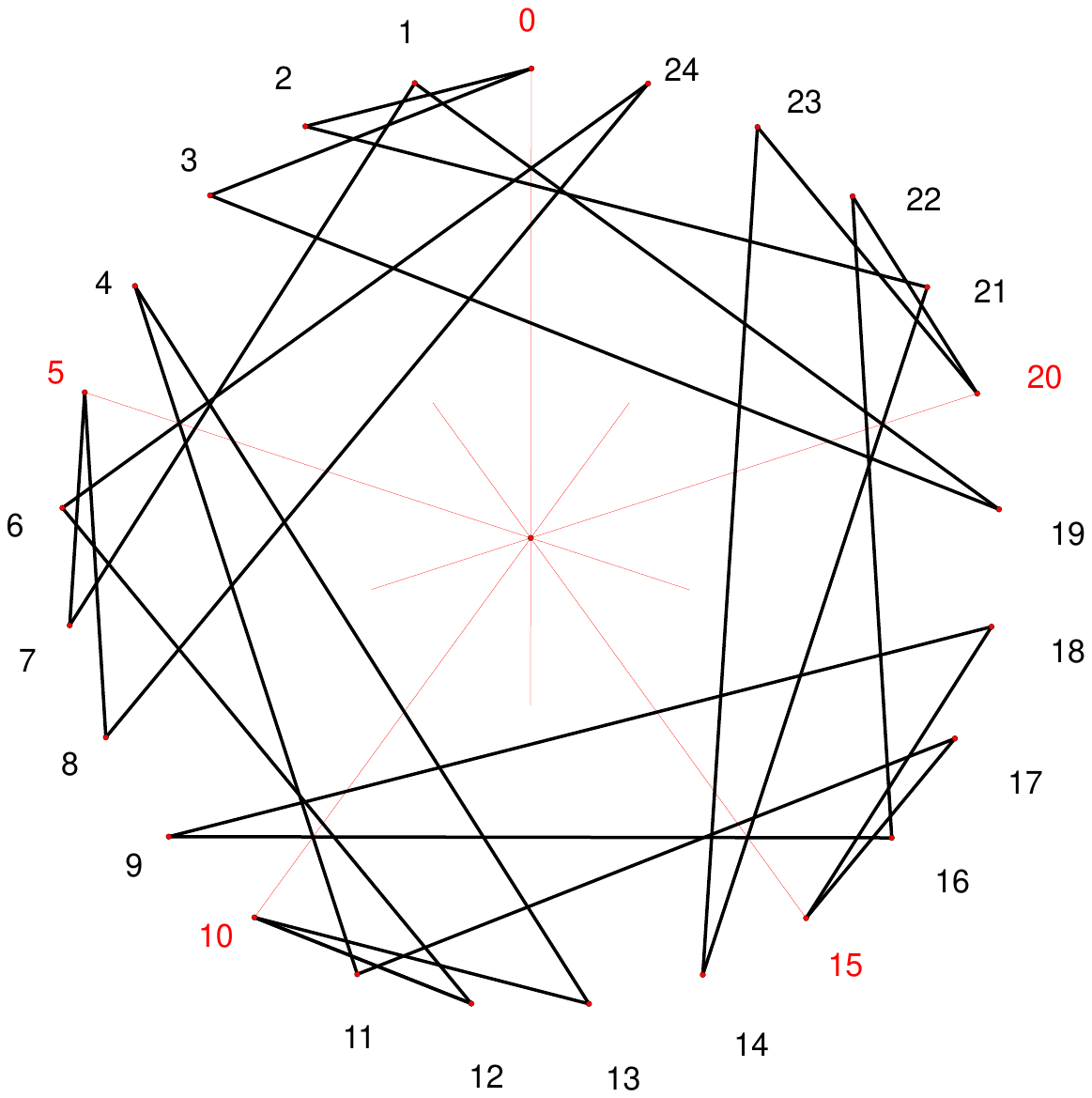}
\caption{Example of a 5-circular 25-polygon}
\end{center}
\end{figure}

\section{Definition of the used equivalence relation}
\label{sec: definition_of_the_used_equivalence_relation}
We denote by $C(n)$ the set of all $n$-polygons and define the following equivalence relation on $C(n)$:

Two $n$-polygons $P_1(n)$ and $P_2(n)$ are said to be equivalent, denoted $P_1(n)\equiv P_2(n)$, if they are obtainable from one another by a rotation.
This means, that the two considered $n$-polygons have the same shape and belong to the same equivalence-class.
Two $n$-polygons $P_1(n)$ and $P_2(n)$  are not equivalent, denoted $P_1(n)\not\equiv P_2(n)$ if they are not obtainable by a rotation. \\\\
\part{$p$-polygons}
\section{Question to deal with in this part of the article}
\label{sec:questions_to_deal_with_in_this_part_of_the_article}

Let $p \geq 5$ be a prime number, then $p$ has the only divisors $p$ and $1$. Therefore there exist only three different kinds of symmetry-categories for $p$-polygons:
\begin{enumerate}
\item[$\bullet$]equivalence classes with $p$ symmetry-axes, (socalled regular or star-polygons)
\item[$\bullet$]equivalence-classes with $1$ symmetry-axis,
\item[$\bullet$]equivalence-classes of the with no symmetry-axis, i.e. the completely asymmetric polygons.
\end{enumerate}

We define and will calculate the following numbers: 
\begin{enumerate}
\item[$\bullet$] $\vert X(p)\vert$ of all equivalence-classes,
\item[$\bullet$] $ \vert X_p(p) \vert$ of the equivalence-classes with $p$ symmetry-axes,
\item[$\bullet$] $ \vert X_1(p) \vert$ of  the equivalence-classes with $1$ symmetry-axis,
\item[$\bullet$] $ \vert X_a(p) \vert$ of  the equivalence-classes of polygons with no symmetry-axis, the completely asymmetrical $p$-polygons.
\end{enumerate}

In addition, we give a table of results and illustrate them by some sets of representatifs for the prime number $p=5$ and $p=7$.

\section{Results}
\label{sec:results}
\subsection{Preparations}
\label{subsec:preparations}
Besides the numbers, which were defined above and which are the goals of the calculations of this part of the article, we need  the following number during the proofs:
$\vert X_{1+}(p) \vert$ is the number of equivalence-classes of $p$-polygons with at least $1$ symmetry-axis.
\subsection{Number of equivalence-classes of all p-polygons}
\label{subsec:number_of_equivalence-classes_of_all_p_polygons}
Let be $p\geq5$ a prim number.
$ \vert X(p) \vert = \dfrac{(p-1)!+(p-1)^{2}}{2p}$
\subsection{The regular p-polygons}
\label{subsec:the_regular_p_polygons}
Let be $p\geq5$ a prim number. 
The number $\vert X_p(p) \vert=\dfrac{\varphi}{2}=\dfrac{p-1}{2}$, while $\varphi(p)$ denotes the Euler divisor function of $p$.\cite{Coxeter}

\subsection{p-polygons with at least 1 symmetry-axis and with exact 1 axis}
\label{subsec:p-polygons_with_at_least_1_symmetry-axis_and_with_exact_1_axis}
Let be $p\geq5$ a prim number. 
A $p$-polygon with at least $1$ symmetry-axis has either $1$ or $p$ symmetry-axes.
 $\vert X_{1+}(p) \vert$ is the number of the equvalence-classes of such $p$-polygons. $\vert X_1(p) \vert$ is the number of the equvalence-classes of $p$-polygons with exact $1$ axe of symmetry.

$\vert X_{1+}(p) \vert=\dfrac{p-1}{2} \cdot2^{\dfrac{p-3}{2}}\cdot \left(\dfrac{p-3}{2}\right)!$

 \begin{center}
 \boxed{\vert X_1(p) \vert=\dfrac{p-1}{2} \cdot [2^{\dfrac{p-3}{2}} \cdot \left(  \dfrac{p-3}{2} \right)!-1]}
 \end{center}

\subsection{The asymmetrical p-polygons}
\label{subsec:the_asymmetrical_p-polygons}

$\vert X_0(p)\vert =\dfrac{\left(p-1\right)^{2}+\left(p-1\right)!-p\cdot \left(p-1\right)\cdot 2^{\dfrac{p-3}{2}}\cdot\left(\dfrac{p-3}{2}\right)!}{2p}$
\section{Proofs}
\label{sec:proofs}
\subsection{Number of equivalence-classes of all p-polygons}
\label{subsec:number_of_equivalence-classes_of_all_p_polygons}
\begin{proof}
This number $\vert X(n) \vert$ is proved by S.W.Golomb and L.R.Welch \cite{Golomb1960} for odd values of $n$:

 $\vert X(n) \vert=\dfrac{1}{2n^2}\left(   \sum \limits_{d \mid n}  \varphi^2 \left( \dfrac{n}{d} \right)\cdot d! \cdot \left( \frac{n}{d} \right)^d \right)$,\\\\

where $d$ are divisors of $n$ and $\varphi$ denotes the Euler divisor function.

Let be $p \geq  5$ a prime number, then the divisors of $p$ are only $1$ and $p$ and the formula of Golomb and Welch changes into:
\begin{center}
$\vert X(p) \vert =\dfrac{(p-1)!+(p-1)^{2}}{2p}$
\end{center}
\end{proof}

\subsection{Number of equivalence-classes of the $p$-polygons with at least one axe}
\label{subsec:number_of_equivalence-classes_of_the_p-polygons_with_at_least_ one_axe}

\begin{proof}
A representatif of each equivalence-class of $p$-polygons with at least one axe of symmetry will get constructed using the following construction figure:
\begin{center}
\begin{figure}[H]
\includegraphics[width=\textwidth]{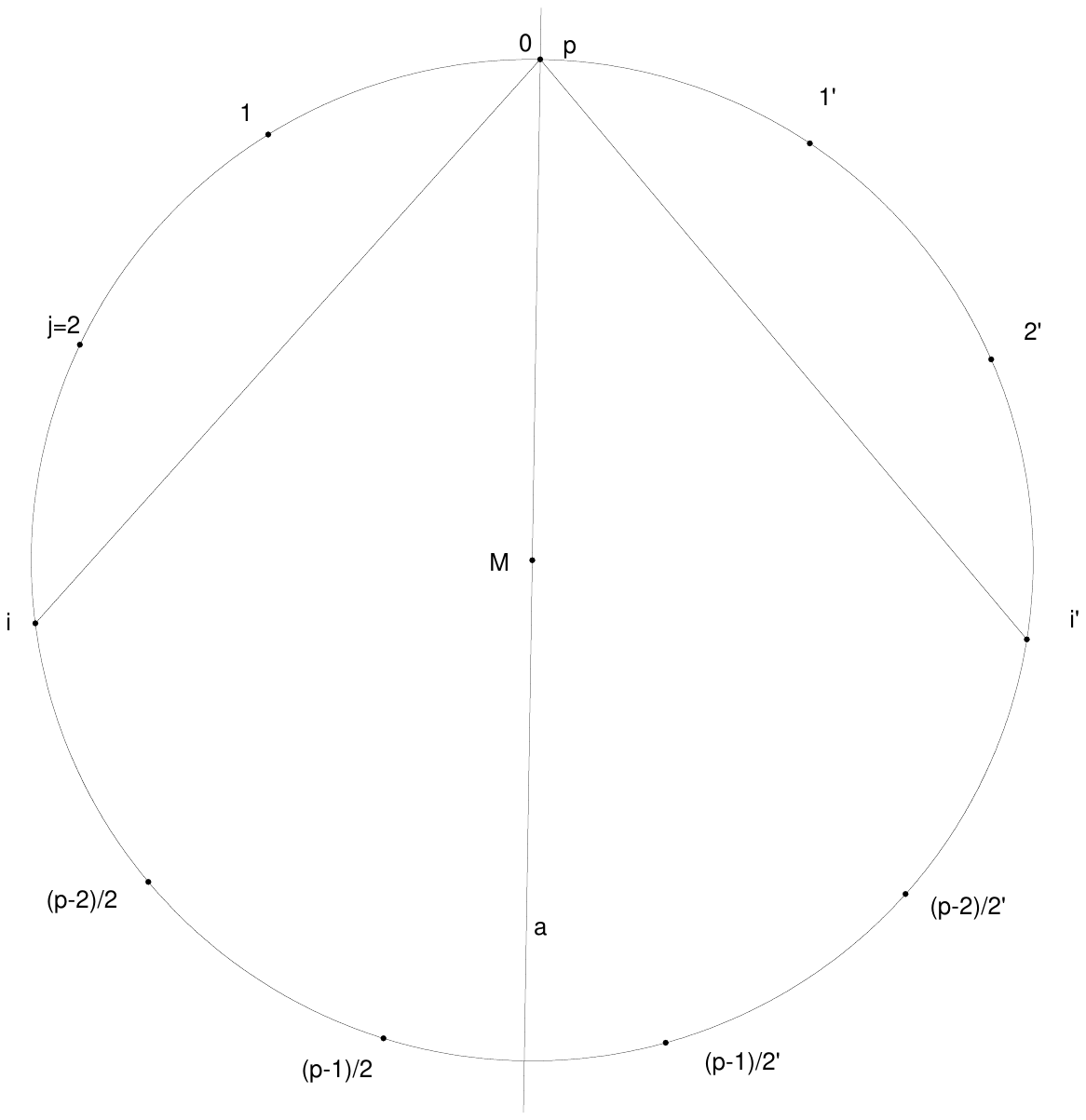}
\caption{Construction-figure 1}
\end{figure}
\end{center}
The representatifs start and end in the vertex labeled with $0$ and $p$, i.e. $v_0=v_p$. One axe is vertical and labeled with $a$. The $\dfrac{p-1}{2}$ left-handed vertices are $v_1,v_2,\dots v_i, \dots v_{\dfrac{p-2}{2}}, v_{\dfrac{p-1}{2}}$, while the right-handed vertices are $v_{1'}=v_{p-1},v_{2'}= v_{p-2}, \dots v_{i'}=v_{p-i} \dots v_{\left(\dfrac{p-2}{2}\right)'}=v_{\dfrac{p+2}{2}},  v_{\left(\dfrac{p-1}{2}\right)'}=v_{\dfrac{p+1}{2}}$.
We choose a first pair of vertices, which are situated symmetrically to the axe $a$ , e.g $v_i$ and $v_{i'}$. We connect those two vertices with the vertex $v_0=v_p$.
So we get an  open chain $C(3)$ of three vertices, symmetric by the axe $a$: $C(3)=\overline{v_iv_0v_{i'}}$. 
\begin{center}
\begin{figure}[H]
\includegraphics[width=\textwidth]{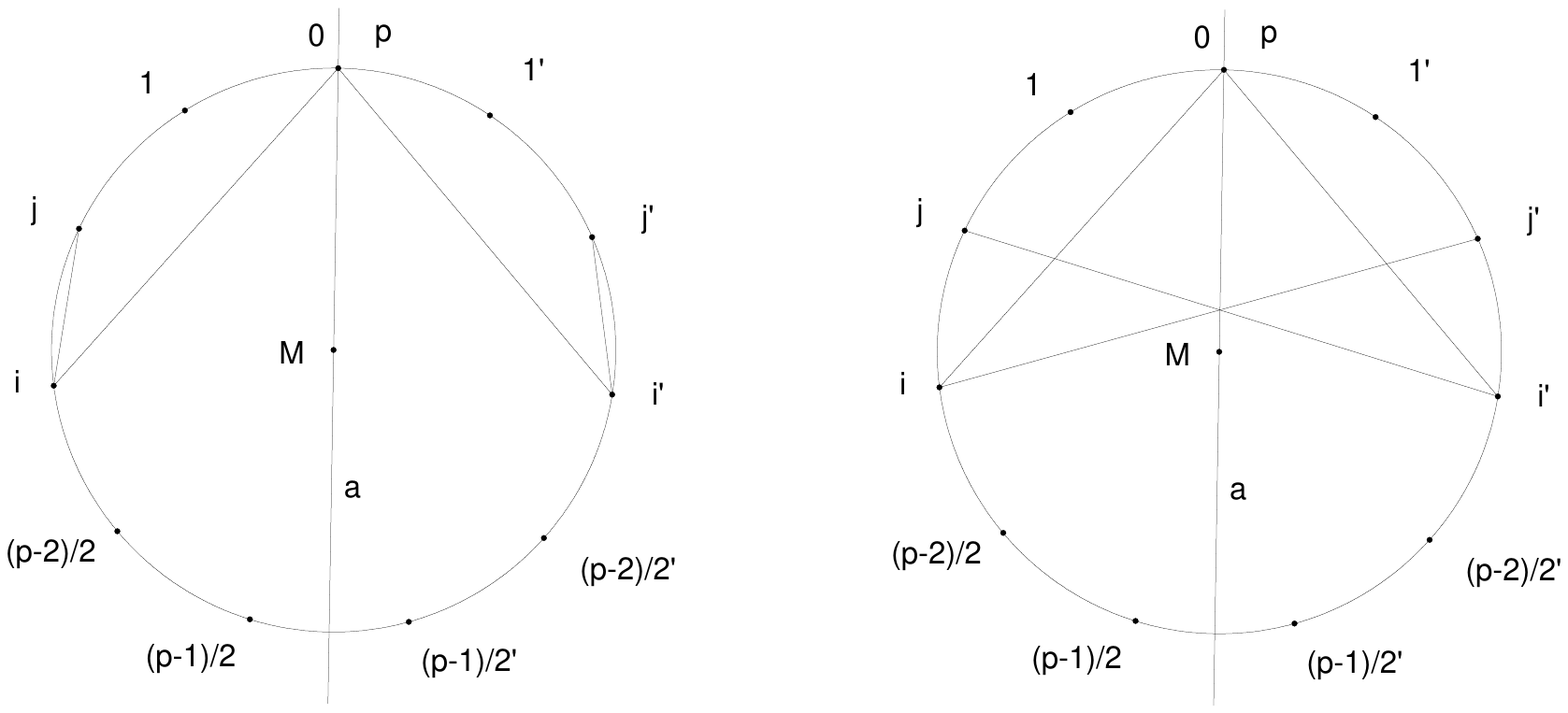}
\caption{construction figure 2}
\end{figure}
\end{center}
If there exist vertices not included in $C(3)$ we choose a pair of vertices $v_j$ and $v_{j'}$ with $j \neq i$. We have two possibilities to connect $v_j$ and $v_{j'}$ with the ends $v_i$ and $v_{i'}$:We get two open chains of 5 vertices:$C_1(5)=\overline{v_jv_iv_0v_{i'}v_{j'}}$ and $C_2(5))=\overline{v_{j'}v_iv_0v_{i'}v_j}$.
If $p=5$ there exist no vertices which are not yet included into the construction and we finish the construction by connecting $v_{j'}$ and $v_j$ to get closed chains of $5$ vertices, representing the equivalence-classes of the $5$-polygons with at least one symmetry axe.\\\\
If $p>5$ we repeat the construction until all $p$ vertices are included and finally close the chains to get the $p$-polygons with at least one axe of symmetry by connecting the two vertices, which were involved into the chains $C(p)$ in the last step.\\\\\

 \textbf{Counting the possibilities:}
\begin{enumerate}
 \item Choosing  the first pair of vertices: $\dfrac{p-1}{2}$ possibilities
 \item Choosing the second pair of vertices: $\dfrac{p-3}{2}$ possibilities
 \item Connecting the second pair with the first pair: $2$ possibilities
 \item The number of chains of 5 vertices: $\vert C(5)\vert =\dfrac{p-1}{2}\cdot\dfrac{p-3}{2}\cdot 2$
 \item Choosing the third pair of vertices: $\dfrac{p-5}{2}$ possibilities
 \item Connecting the third pair with the second pair: $2$ possibilities
 \item ...
\item Choosing the last pair of vertices: $1=p-(p-1)$ possibility
\item Connecting the last pair of vertices with the second to last: $2$ possibilities
\item Closing the chain: 1 possibility
\end{enumerate}
Together: $\vert P_{1+}(p) \vert = 2^{\dfrac{p-3}{2}}\cdot\dfrac{p-3}{2}!\cdot \left(\dfrac{p-1}{2}\right)$

Due to Leonhard Euler:$\vert P_p(p)\vert=\dfrac{\varphi}{2}=\dfrac{p-1}{2}$.\\\\
We conclude:
\begin{center}
\fbox{$\vert P_1(p)\vert=\dfrac{p-1}{2} \cdot \left[2^{\dfrac{p-3}{2}} \cdot \left(\dfrac{p-3}{2}\right)!-1\right]$}
\end{center}
\end{proof}

\subsection{Number of the equivalence-classes of asymmetrical $p$-polygons}
\label{subsec:number_of_the_equivalence-classes_of_asymmetrical_p-polygons}
\begin{proof}
$\vert P_0(p) \vert = \vert P(p) \vert - \vert P_{1+}(p) \vert $\\\\
Therefore:\\\\
$\vert P_0(p) \vert = \dfrac{(p-1)!+(p-1)^{2}}{2p}-2^{\dfrac{p-3}{2}}\cdot\left(\dfrac{p-3}{2}\right)!\cdot \left(\dfrac{p-1}{2}\right)$

It follows:
\begin{center}
\fbox{$\vert X_0(p)\vert =\dfrac{\left(p-1\right)^{2}+\left(p-1\right)!-p\cdot \left(p-1\right)\cdot 2^{\dfrac{p-3}{2}}\cdot\left(\dfrac{p-3}{2}\right)!}{2p}$}
\end{center}
\end{proof}
\section{Table of the results}
\label{sec:table_of_the_results}
\begin{center}
\begin{table}[H]
\centering
\begin{tabular}{|| c | c | c | c | c ||}
 p & $\vert X(p) \vert$ & $\vert X_p(p)\vert$ & $\vert X_1(p)\vert$ & $\vert X_0(p)\vert$\\ \hline
5 & 4 & 2 & 2 & 0\\ \hline
7 & 54 & 3 & 21 & 30\\ \hline
11&164950 &5 & 1915 & 163030 \\ \hline
13 & 18423144 & 6 & 23034 & 1840010\\ \hline
\end{tabular}
\caption{Result for $5\leq p \leq 13$}
\end{table}
\end{center}

\section{Representatifs of the equivalence-classes}
\label{sec:representatifs_of_the_equivalence-classes}
\subsection{Representation of the $5$-polygons}
\label{subsec:representation_of_the_5-polygons}
\begin{center}
\begin{figure}[H]
\centering
\begin{tabular}{| c | c | c | c |}
\hline
\includegraphics[width=0.2\textwidth]{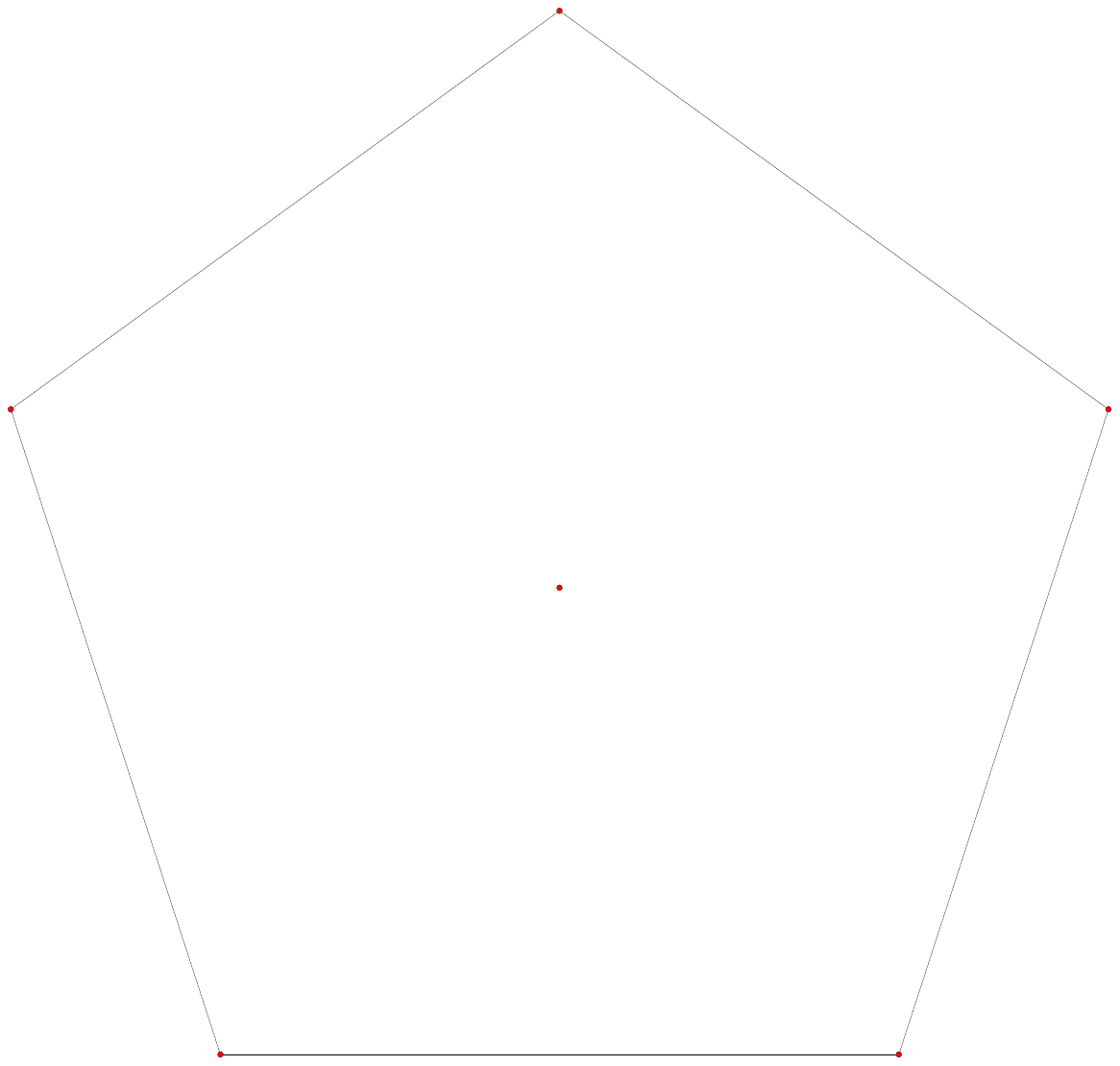} & \includegraphics[width=0.2\textwidth]{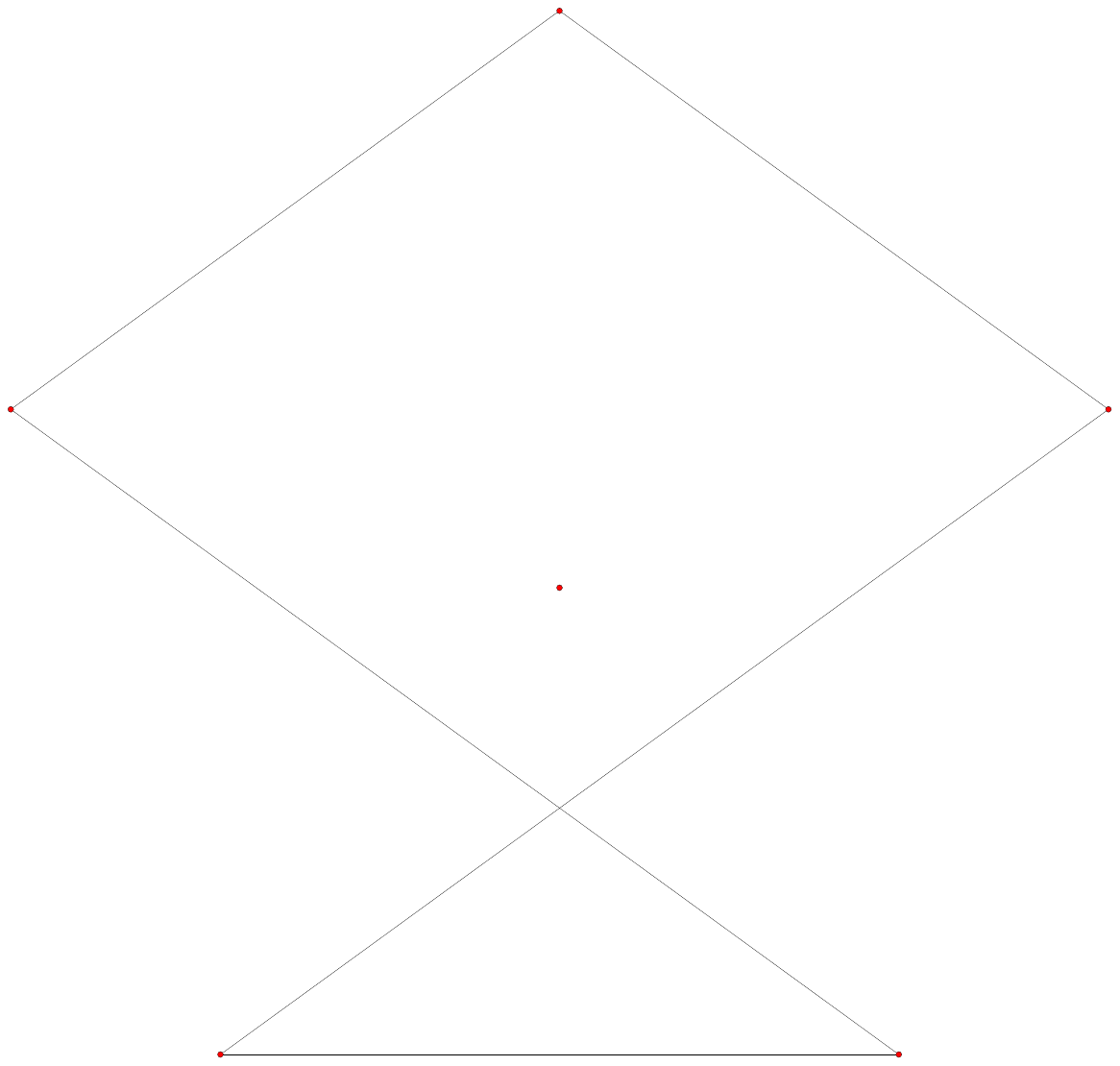} & \includegraphics[width=0.2\textwidth]{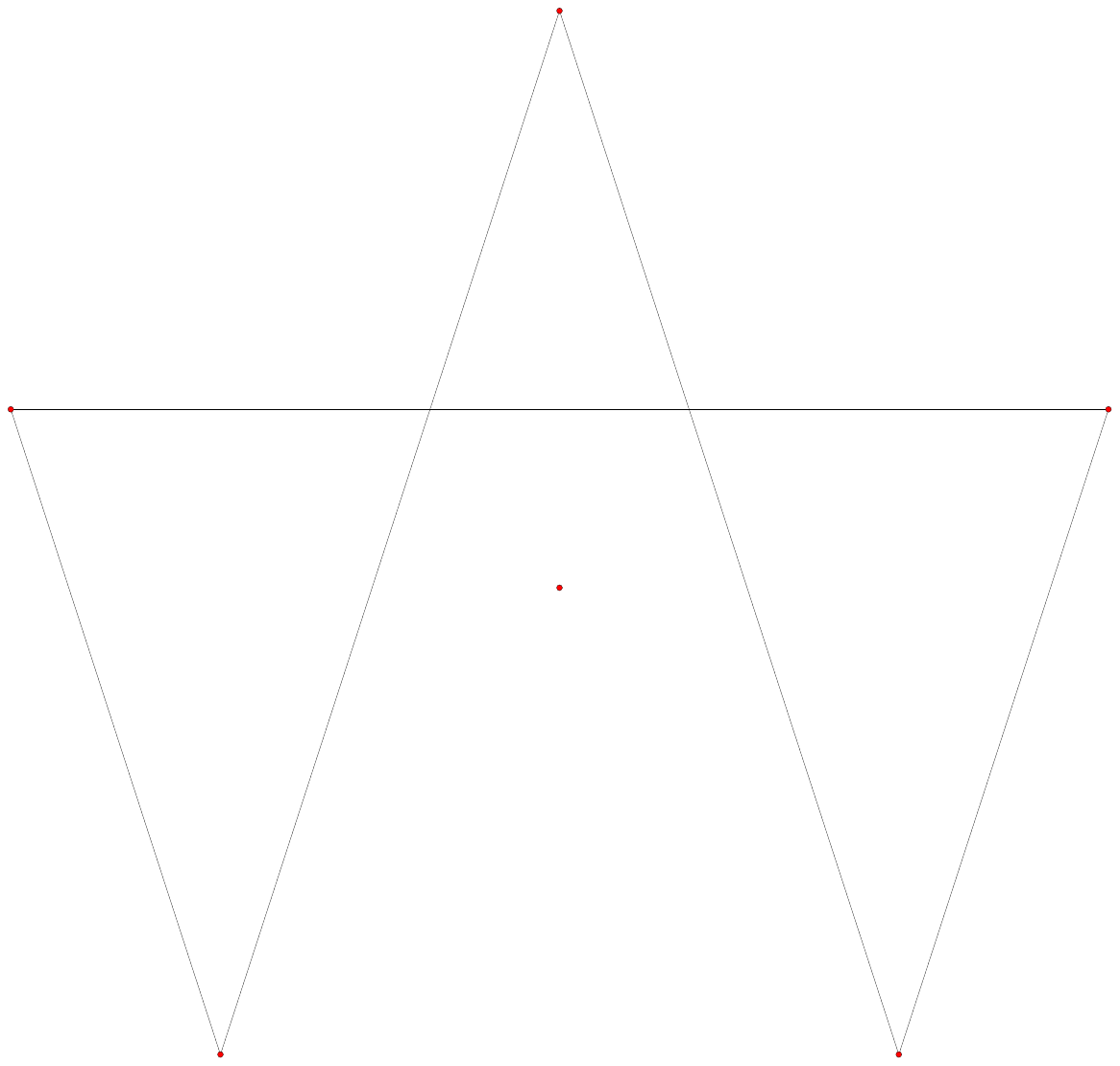} & \includegraphics[width=0.2\textwidth]{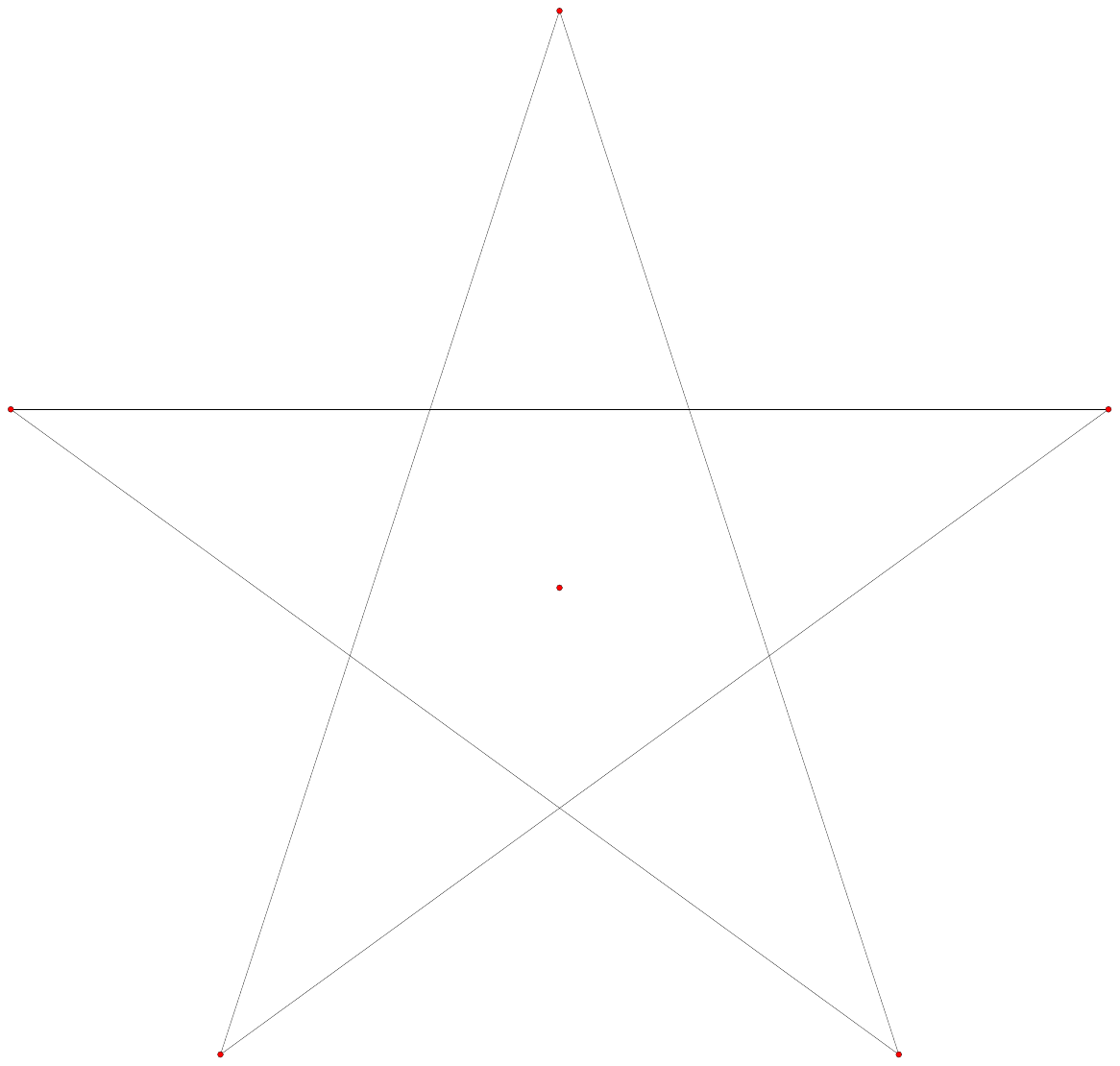}\\ \hline
\end{tabular}
\caption{representatifs of the 4 equivalence-classes of the $5$-polygons}
\end{figure}
\end{center}

\newpage

\subsection{Representation of the asymmetrical $7$-polygons}
\label{subsec:representation_of_the_asymmetrical_7-polygons}
\begin{center}
\begin{figure}[H]
\centering
\begin{tabular}{| c | c | c |}
\hline
\includegraphics[width=0.3\textwidth]{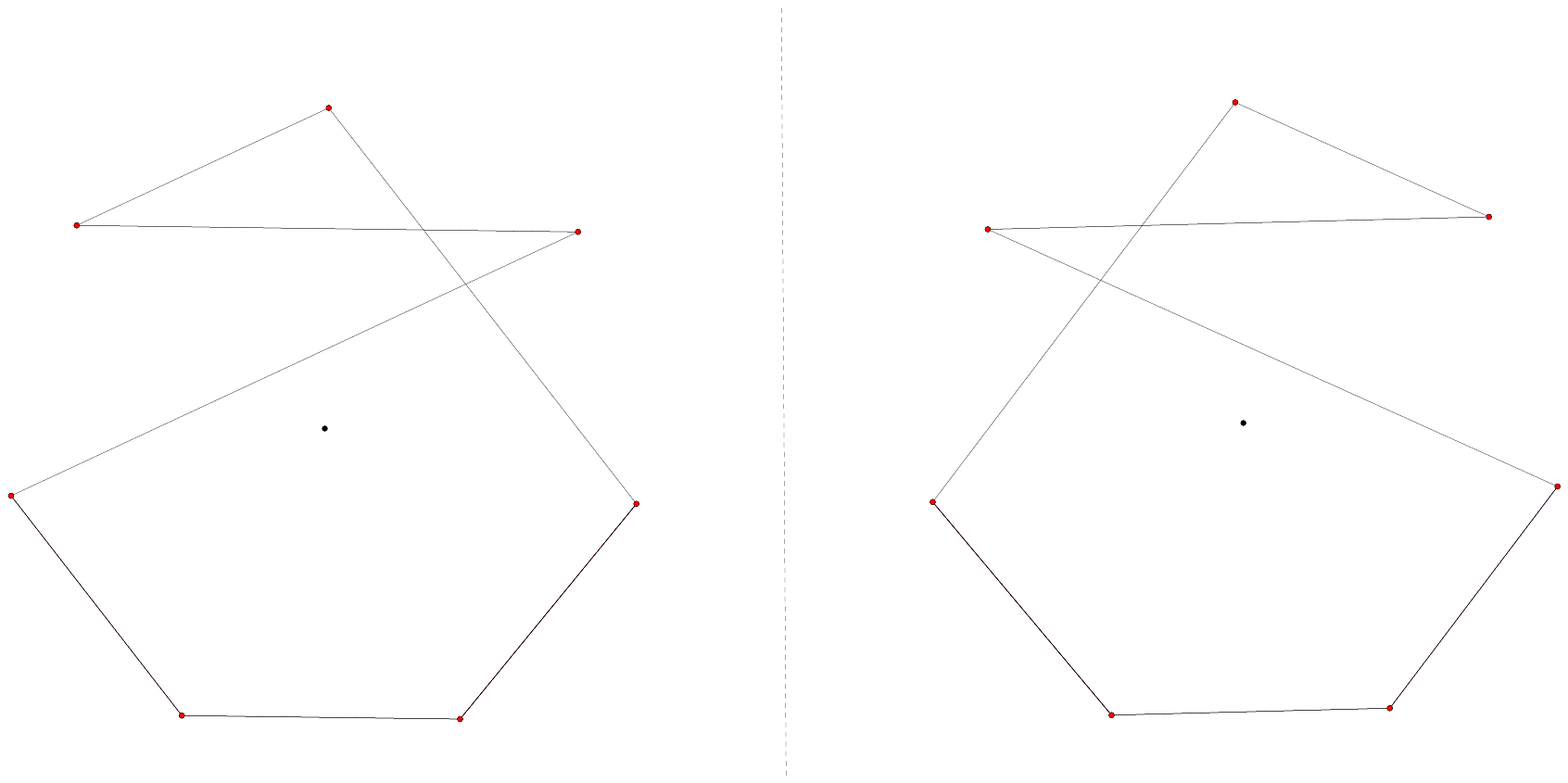} & \includegraphics[width=0.3\textwidth]{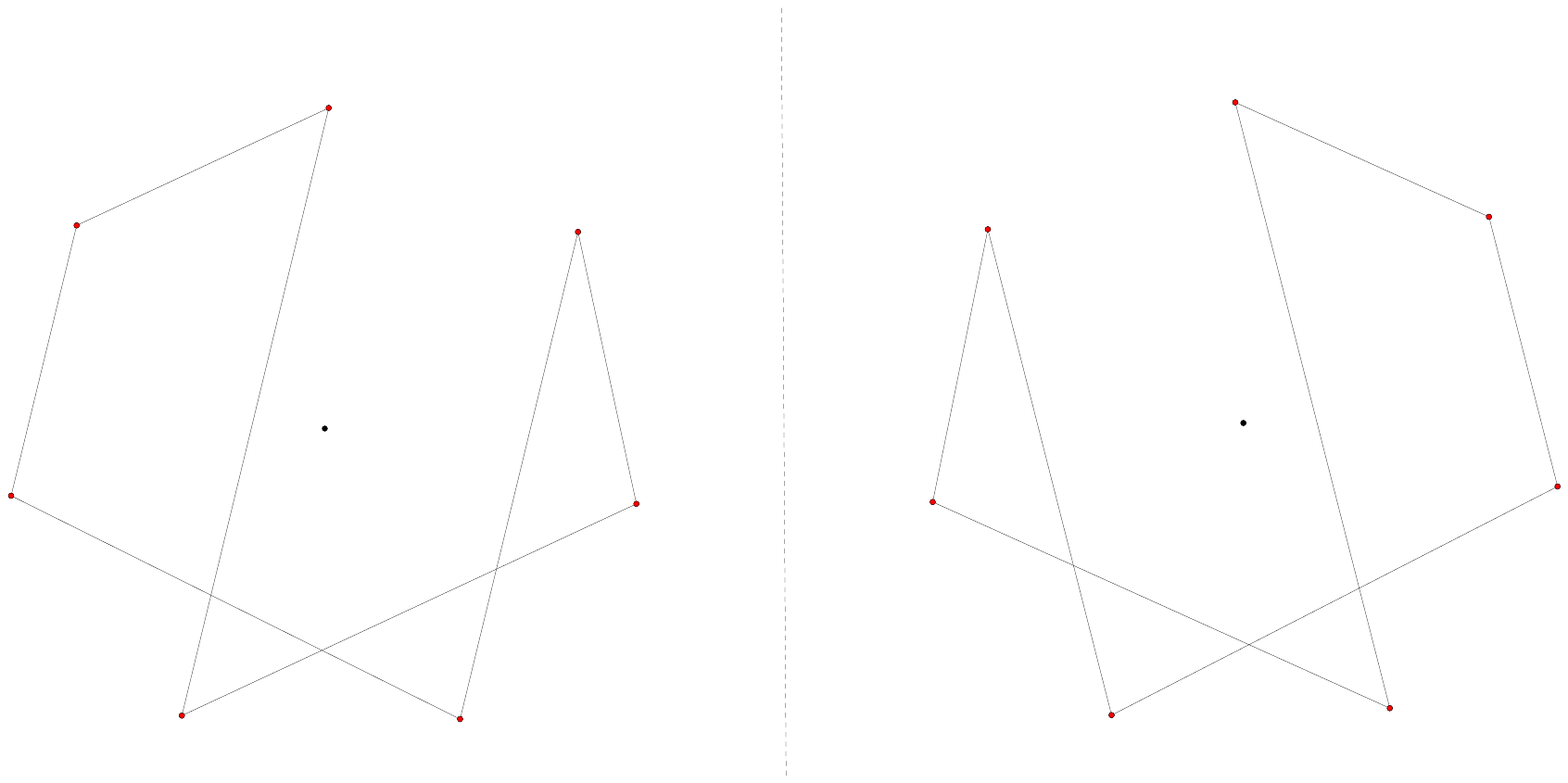} & \includegraphics[width=0.3\textwidth]{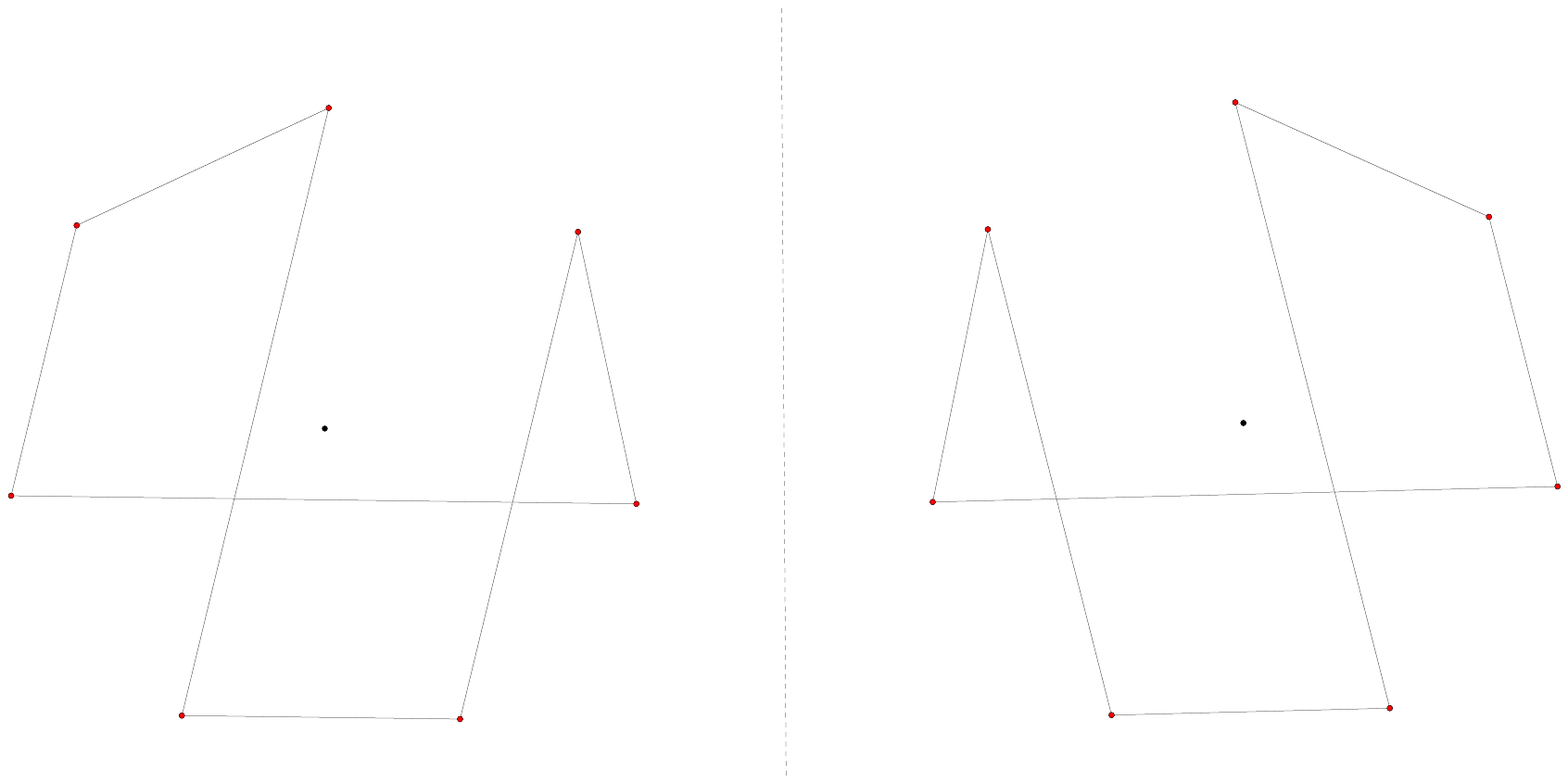}\\ \hline \includegraphics[width=0.3\textwidth]{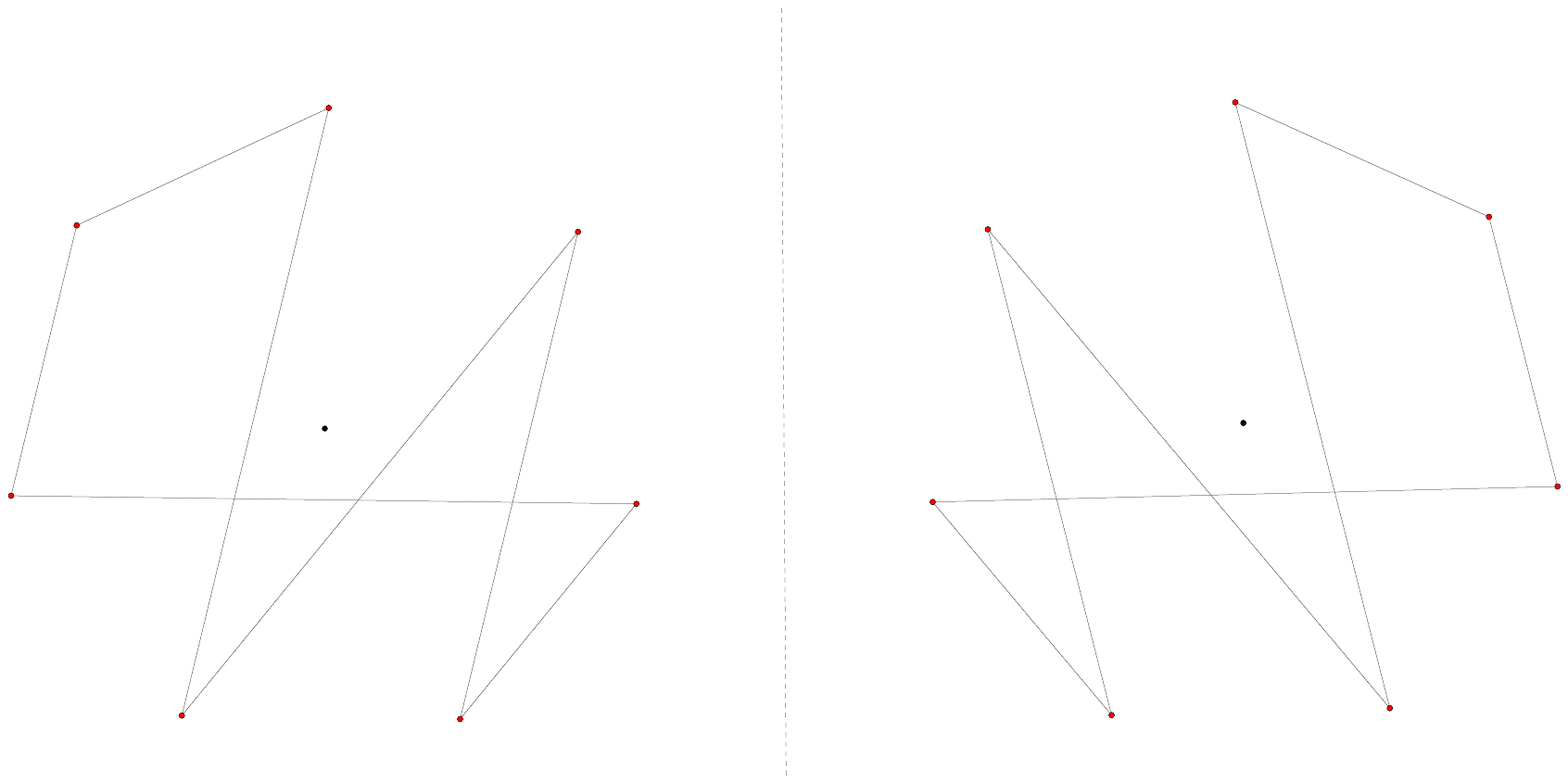} &
\includegraphics[width=0.3\textwidth]{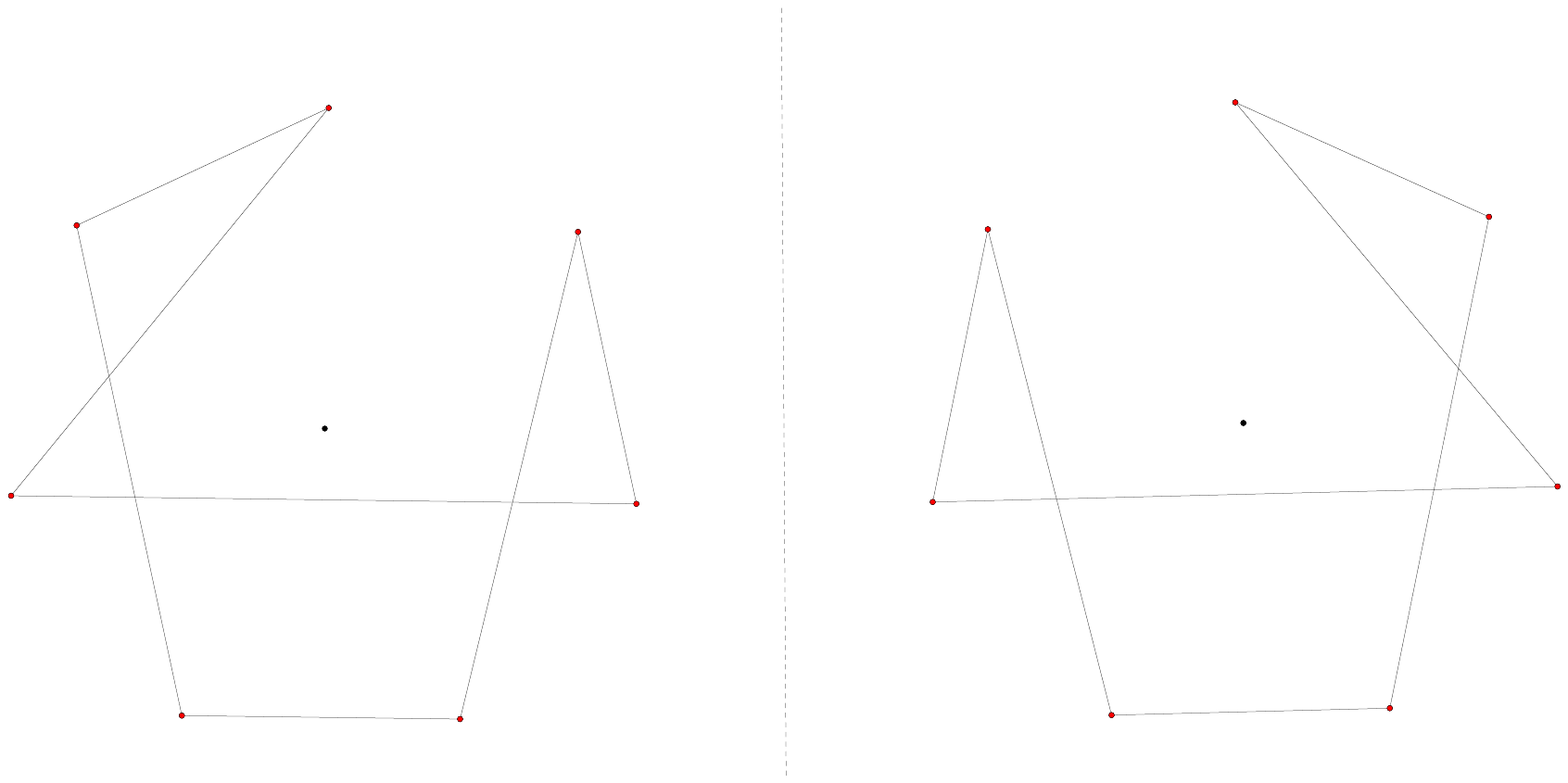} & \includegraphics[width=0.3\textwidth]{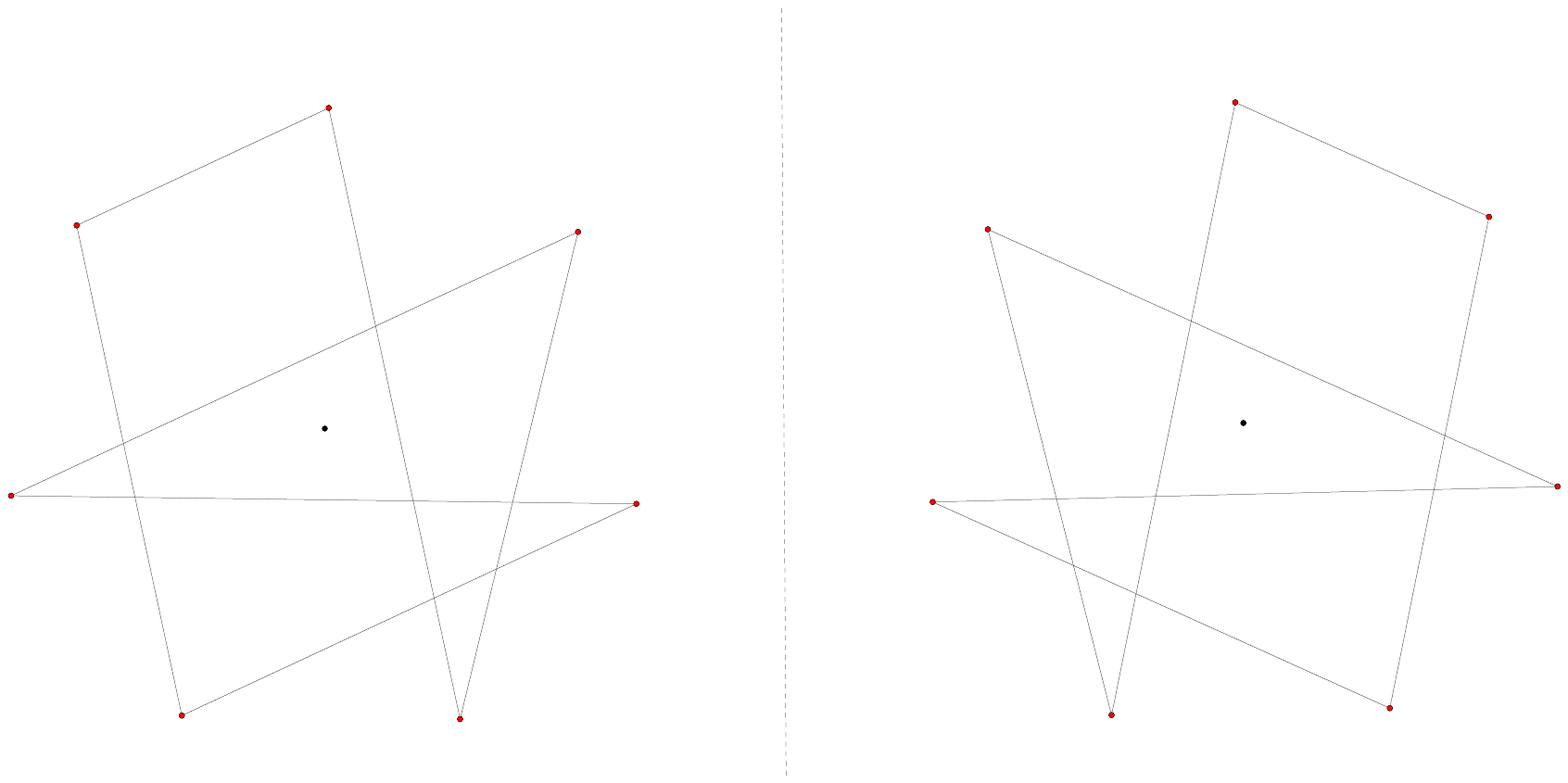}\\ \hline \includegraphics[width=0.3\textwidth]{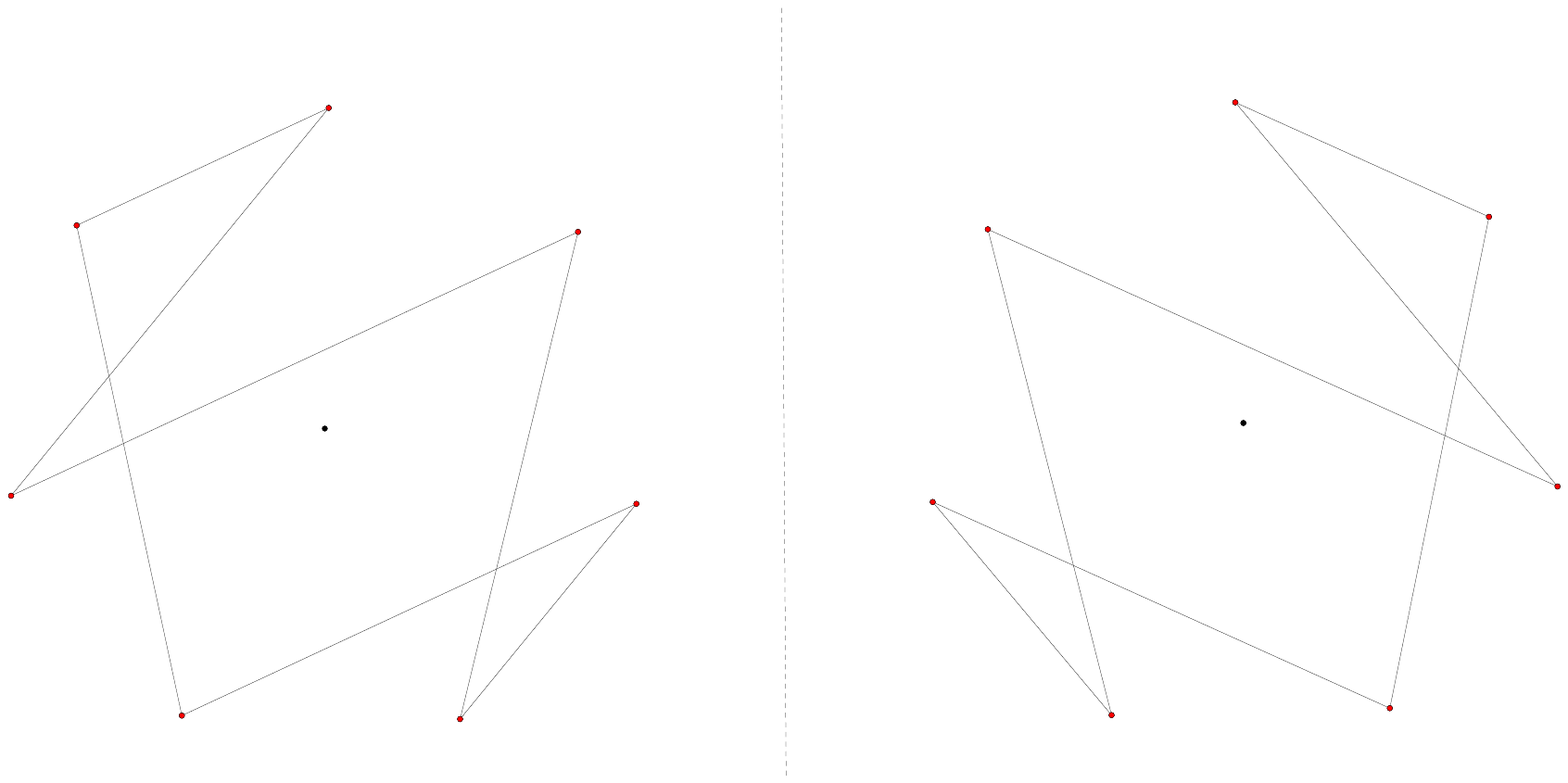} & \includegraphics[width=0.3\textwidth]{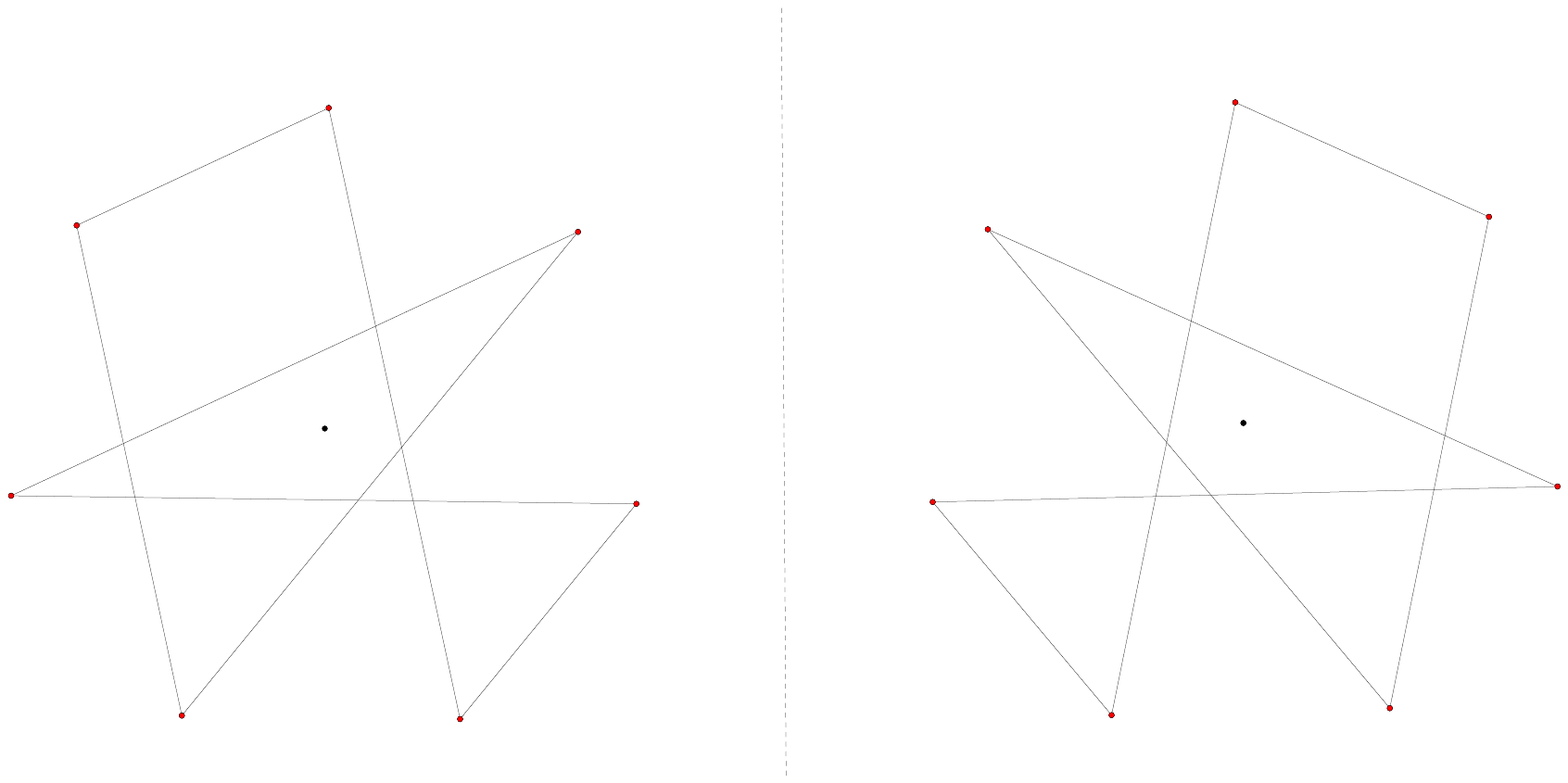} &
\includegraphics[width=0.3\textwidth]{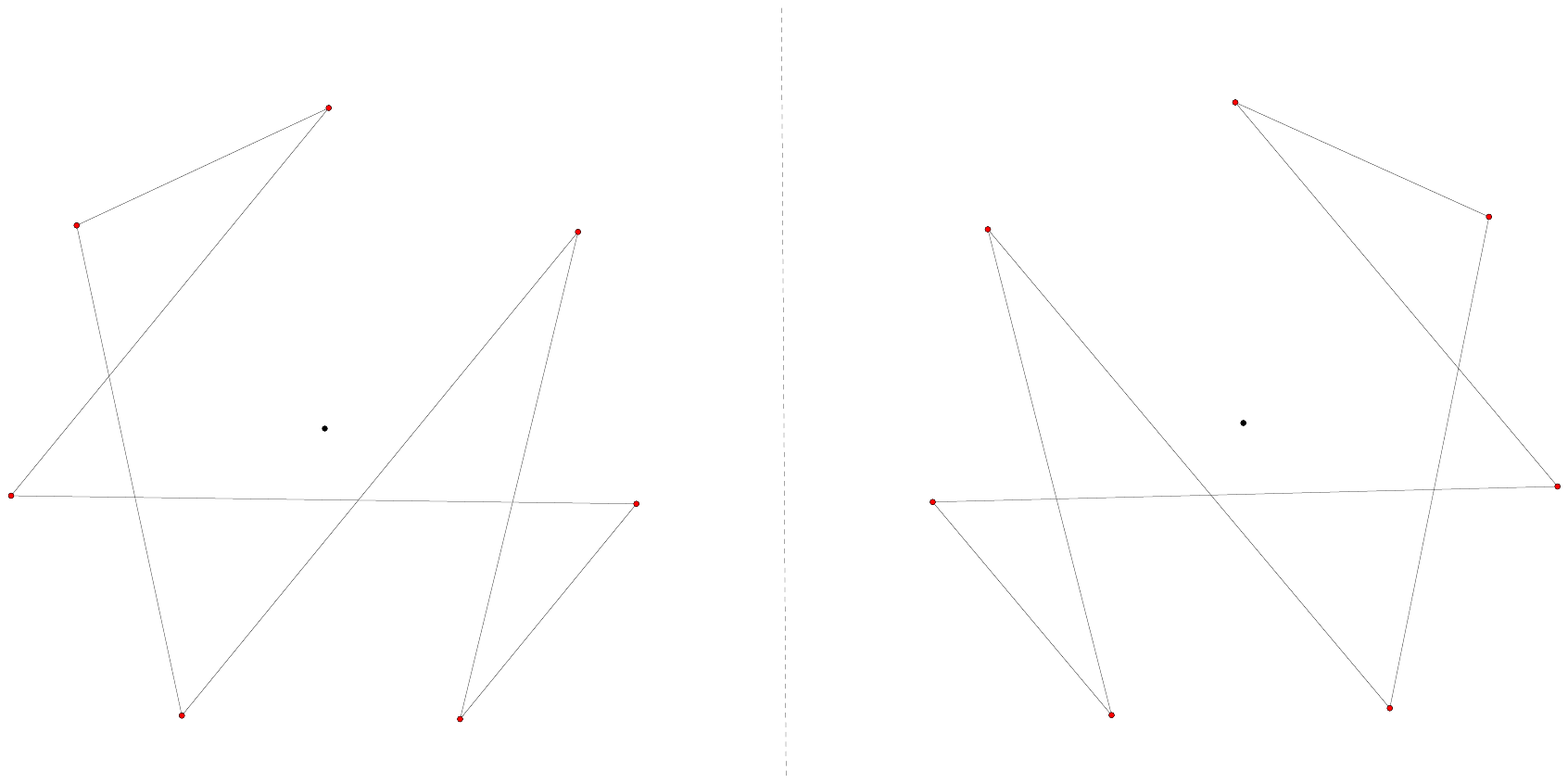}\\ \hline \includegraphics[width=0.3\textwidth]{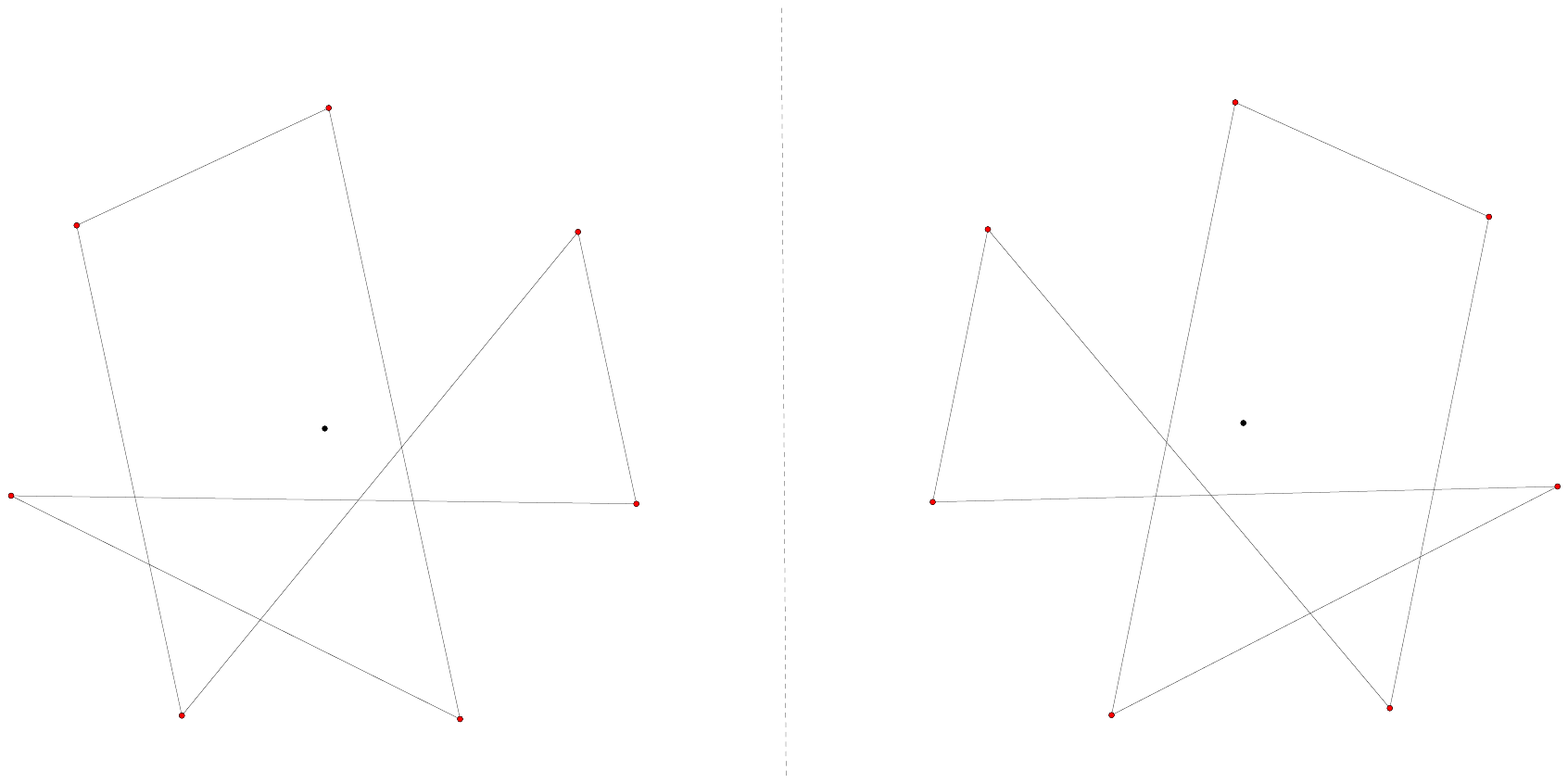} & \includegraphics[width=0.3\textwidth]{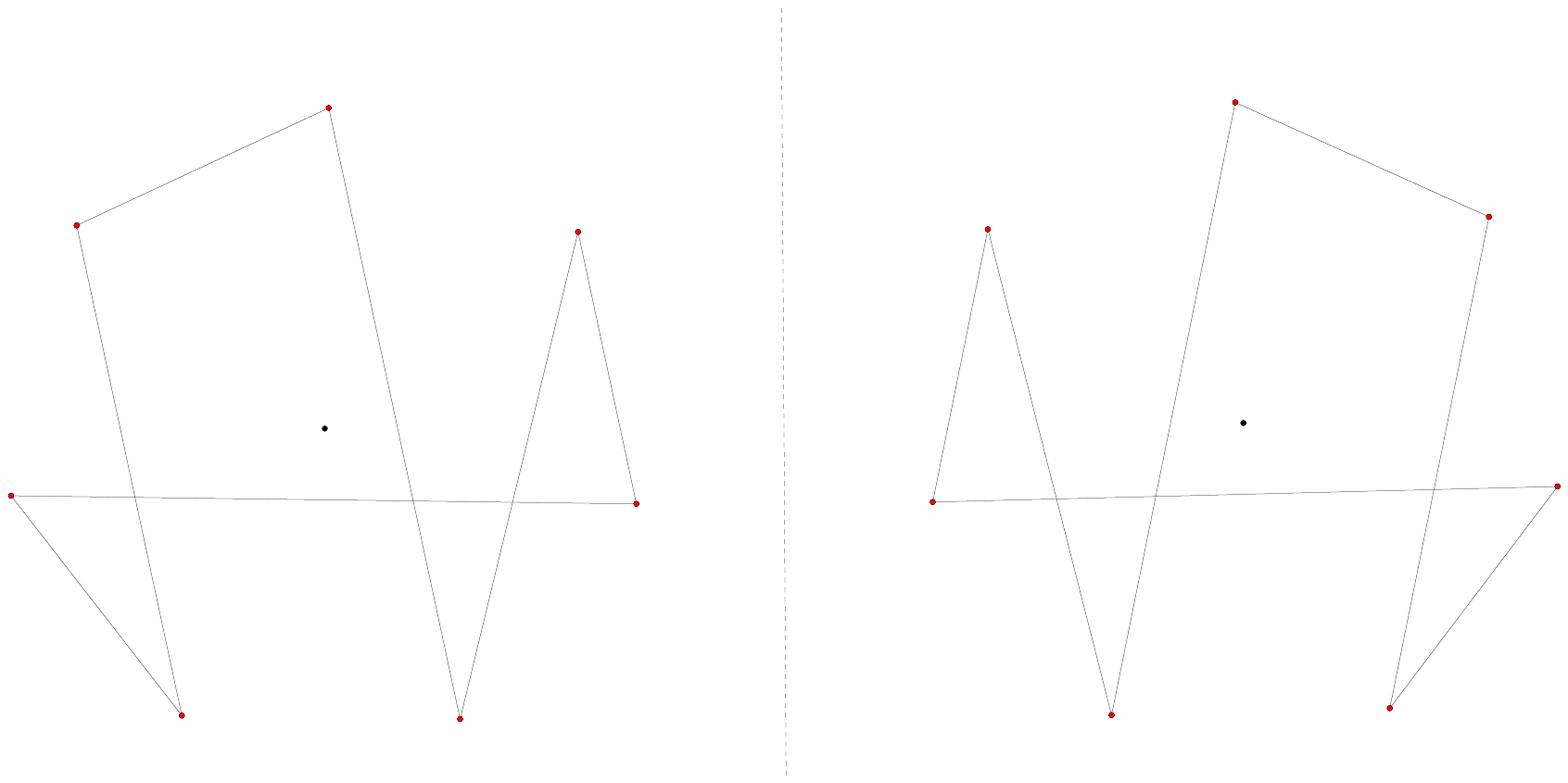} & \includegraphics[width=0.3\textwidth]{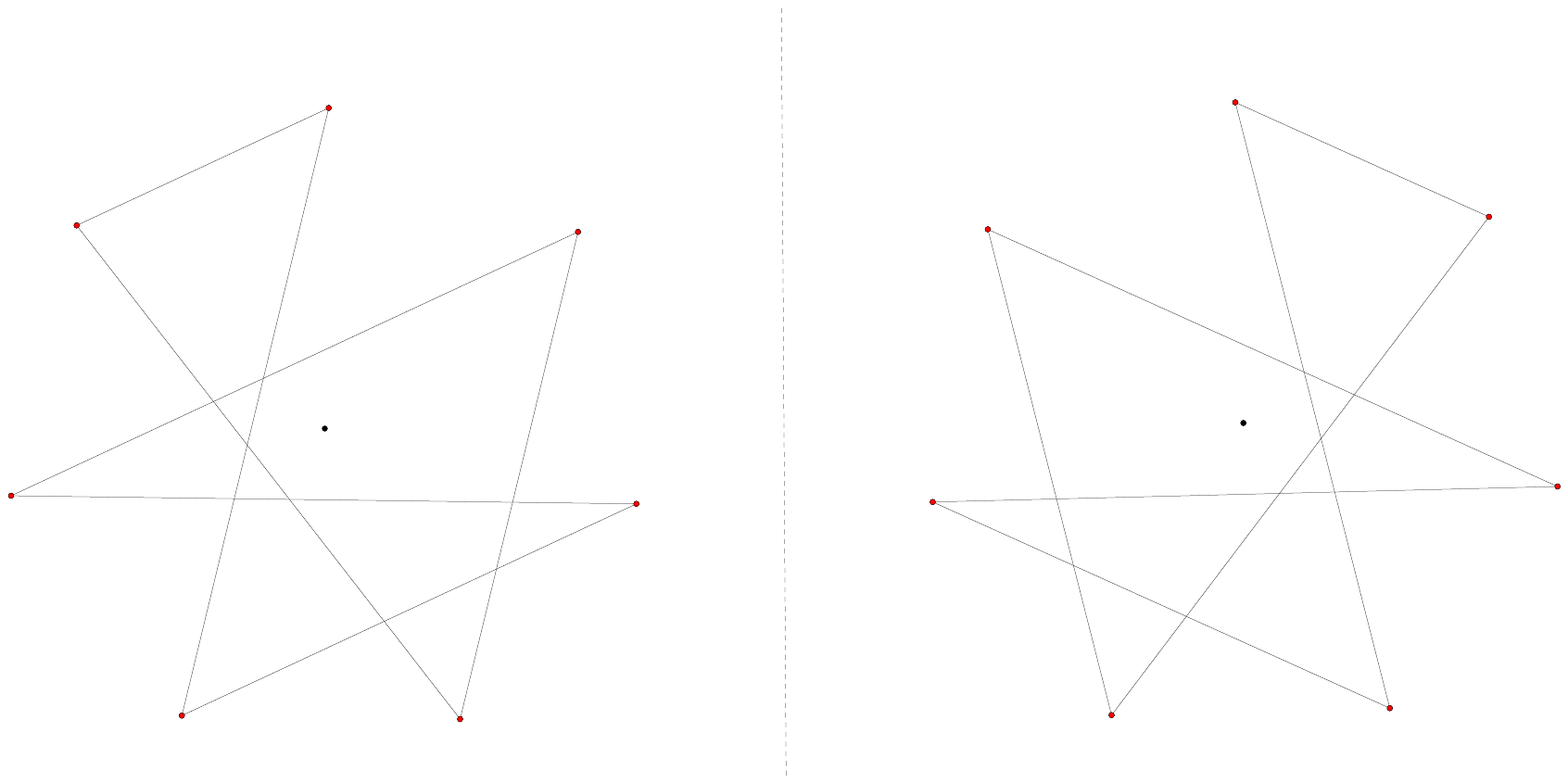}\\ \hline
\includegraphics[width=0.3\textwidth]{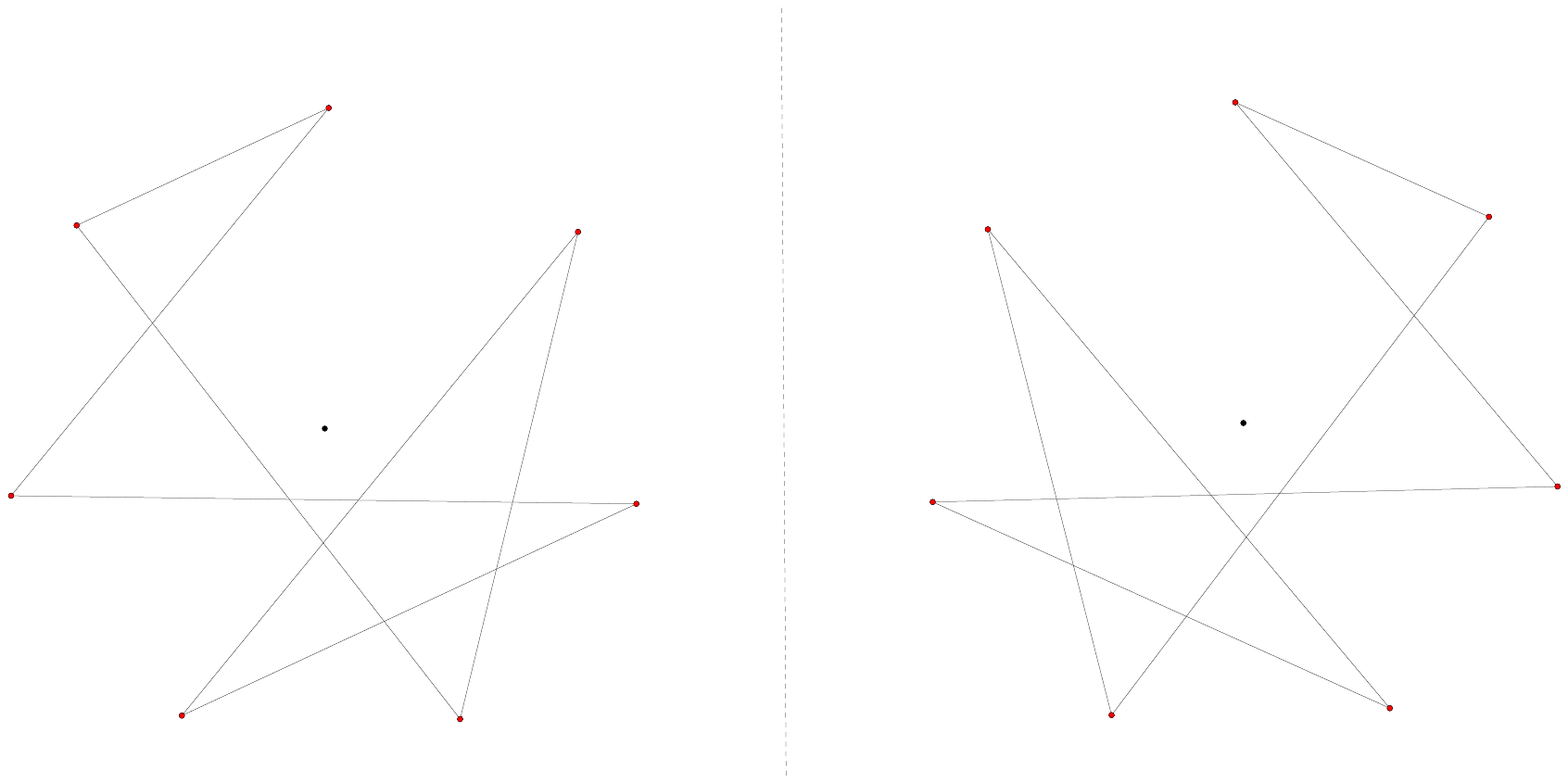} & \includegraphics[width=0.3\textwidth]{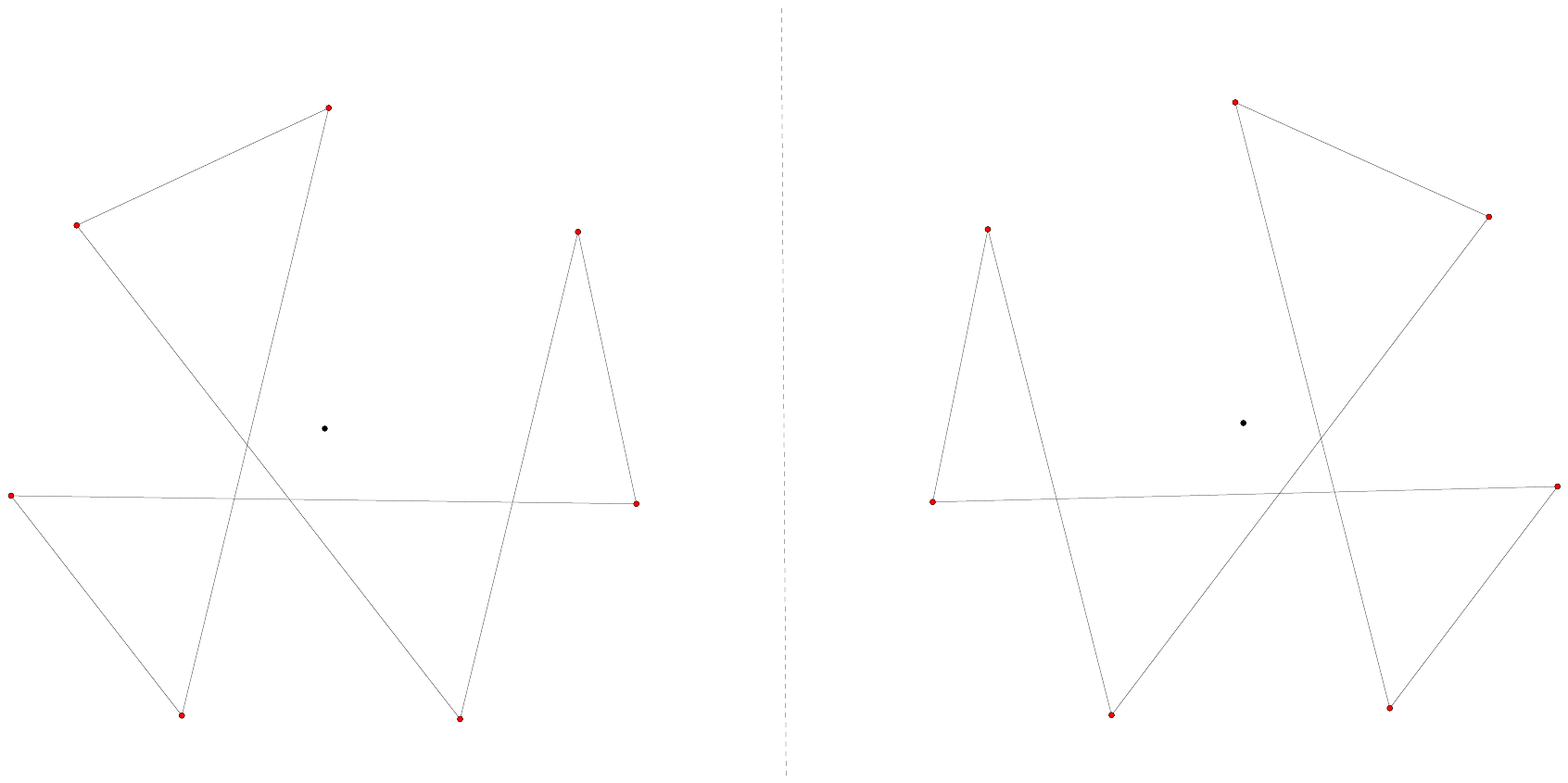} & \includegraphics[width=0.3\textwidth]{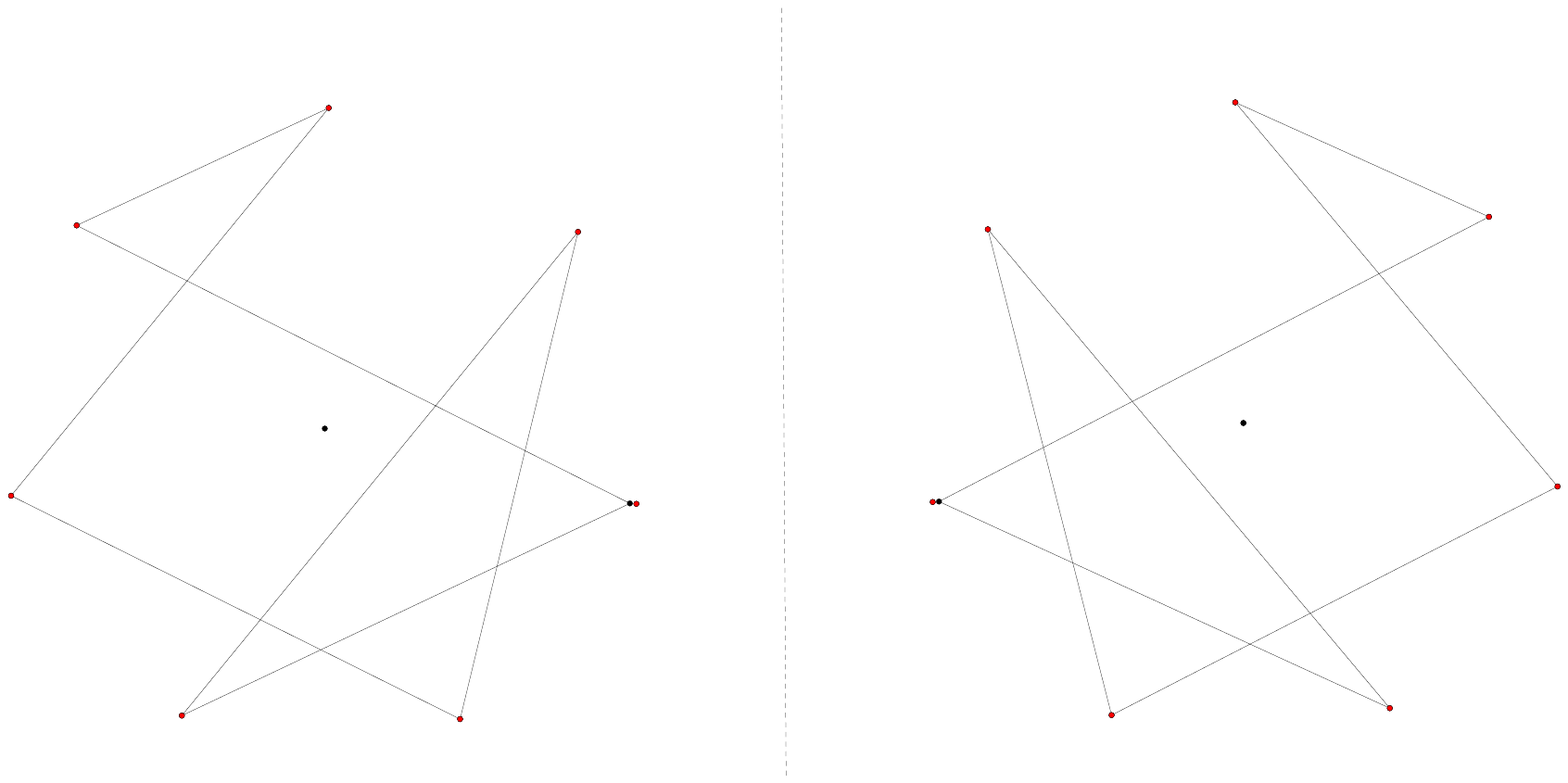}\\ \hline
\end{tabular}
\caption{representatifs of the 30 equivalence-classes of asymmetrical $7$-polygons}
\end{figure}
\end{center}
Between the two mirrored representatifs is always drawn a small vertical line. They belong to two different equivalence-classes.

\newpage
\subsection{Representation of the 24 symmetrical $7$-polygons}
\label{subsec:representation_of_the_24_symmetrical_7-polygons}
\begin{center}
\begin{figure}[H]
\centering
\begin{tabular}{| c | c | c | c |}
\hline
\includegraphics[width=0.2\textwidth]{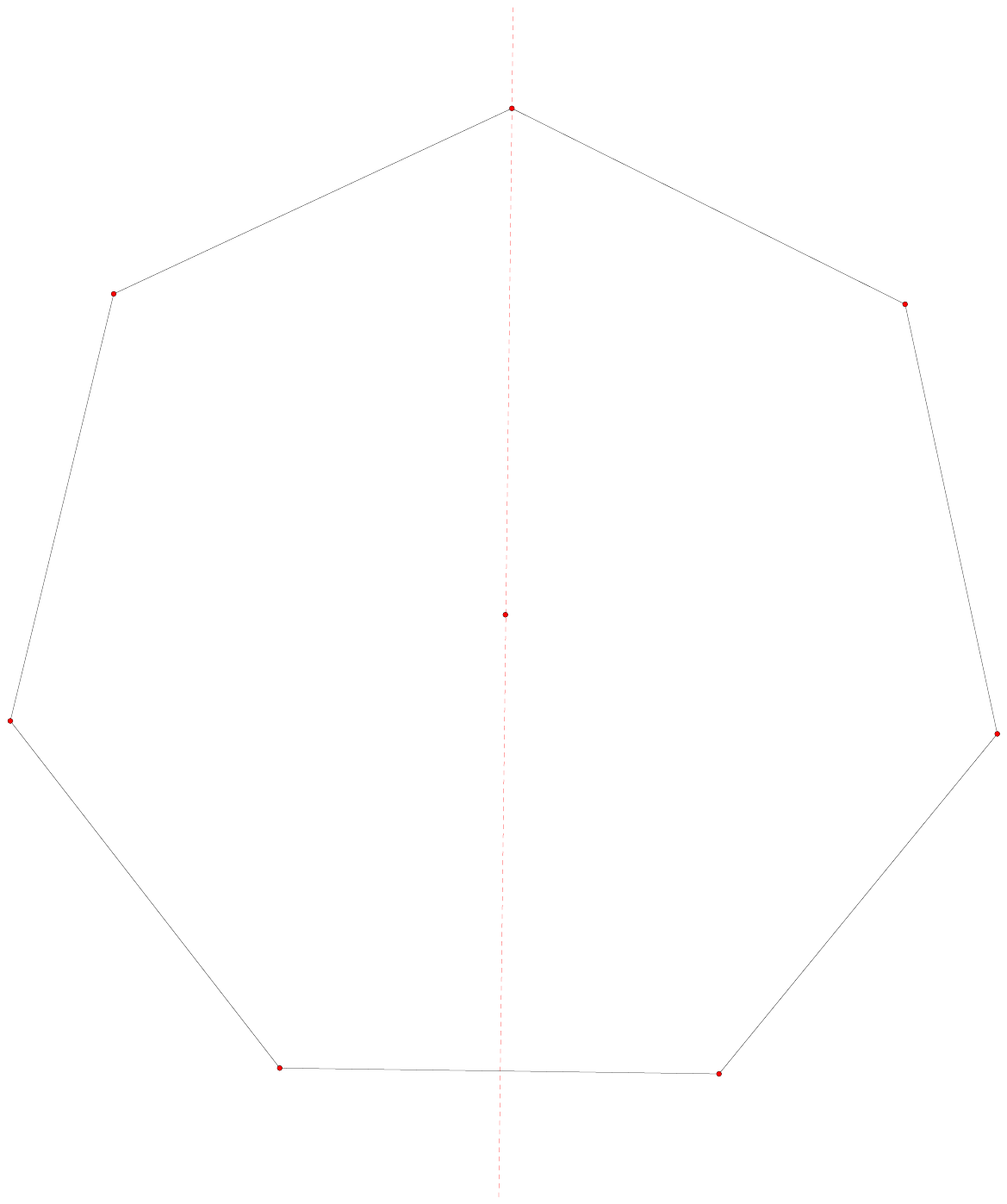} & \includegraphics[width=0.2\textwidth]{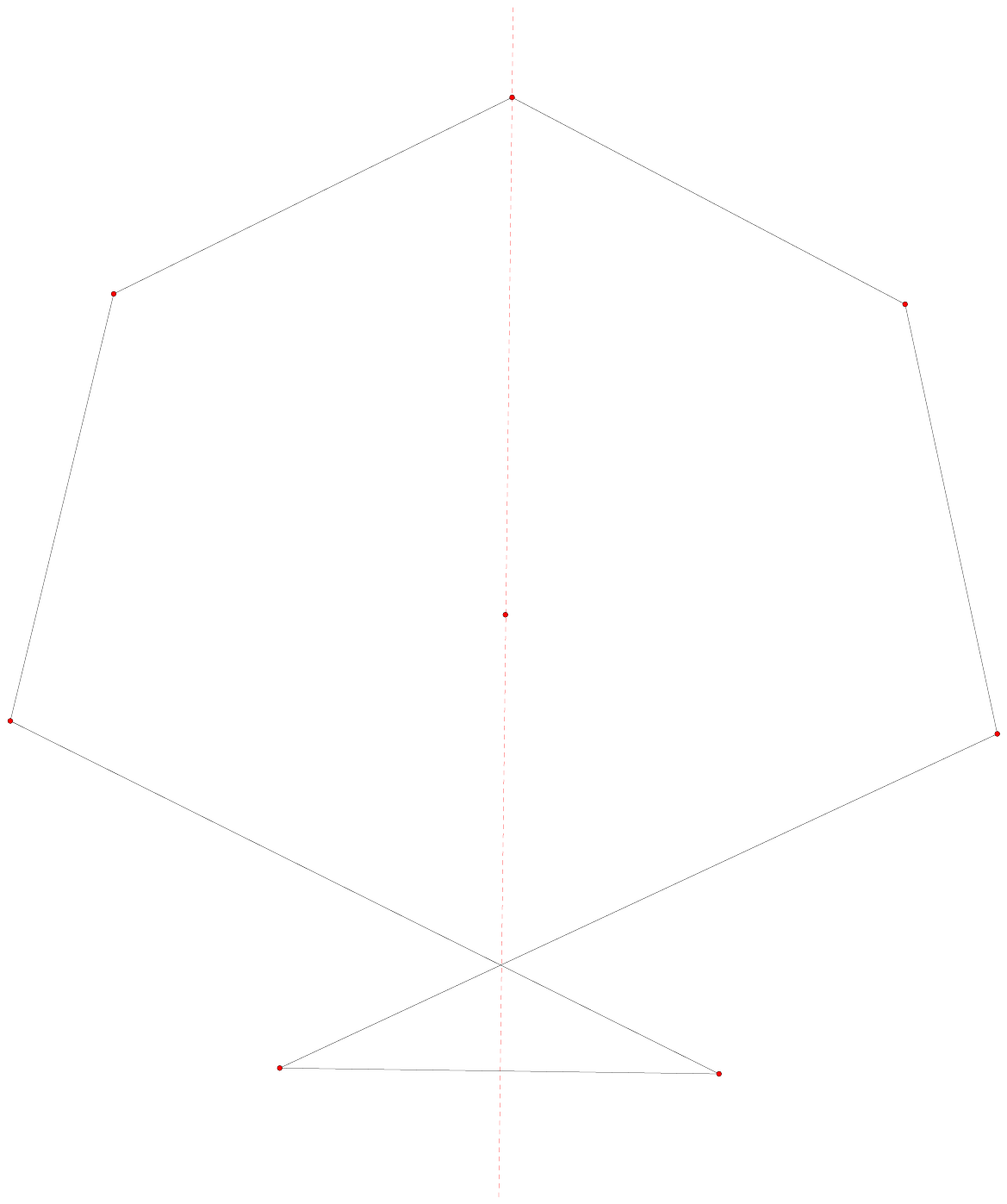} & \includegraphics[width=0.2\textwidth]{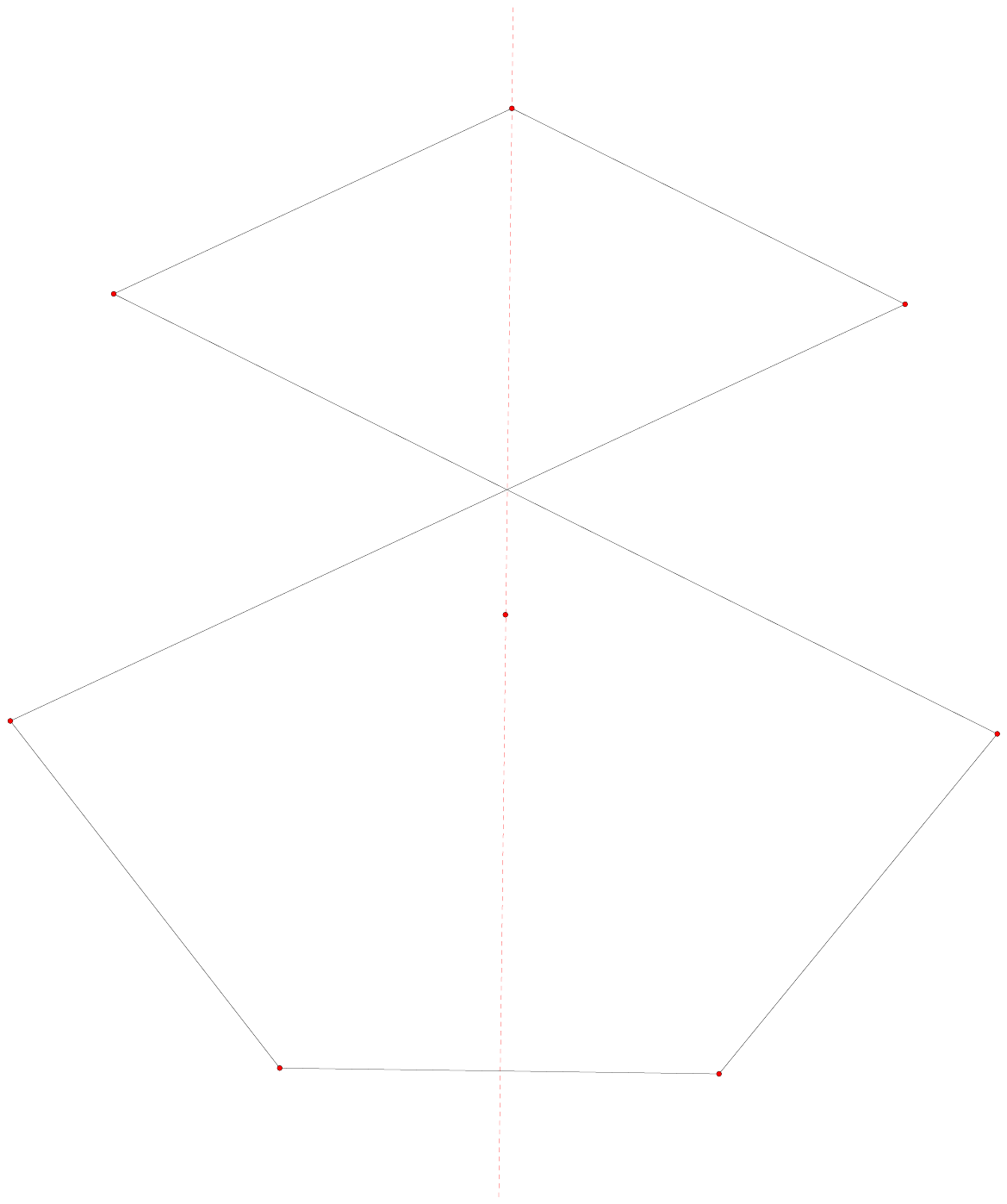} & \includegraphics[width=0.2\textwidth]{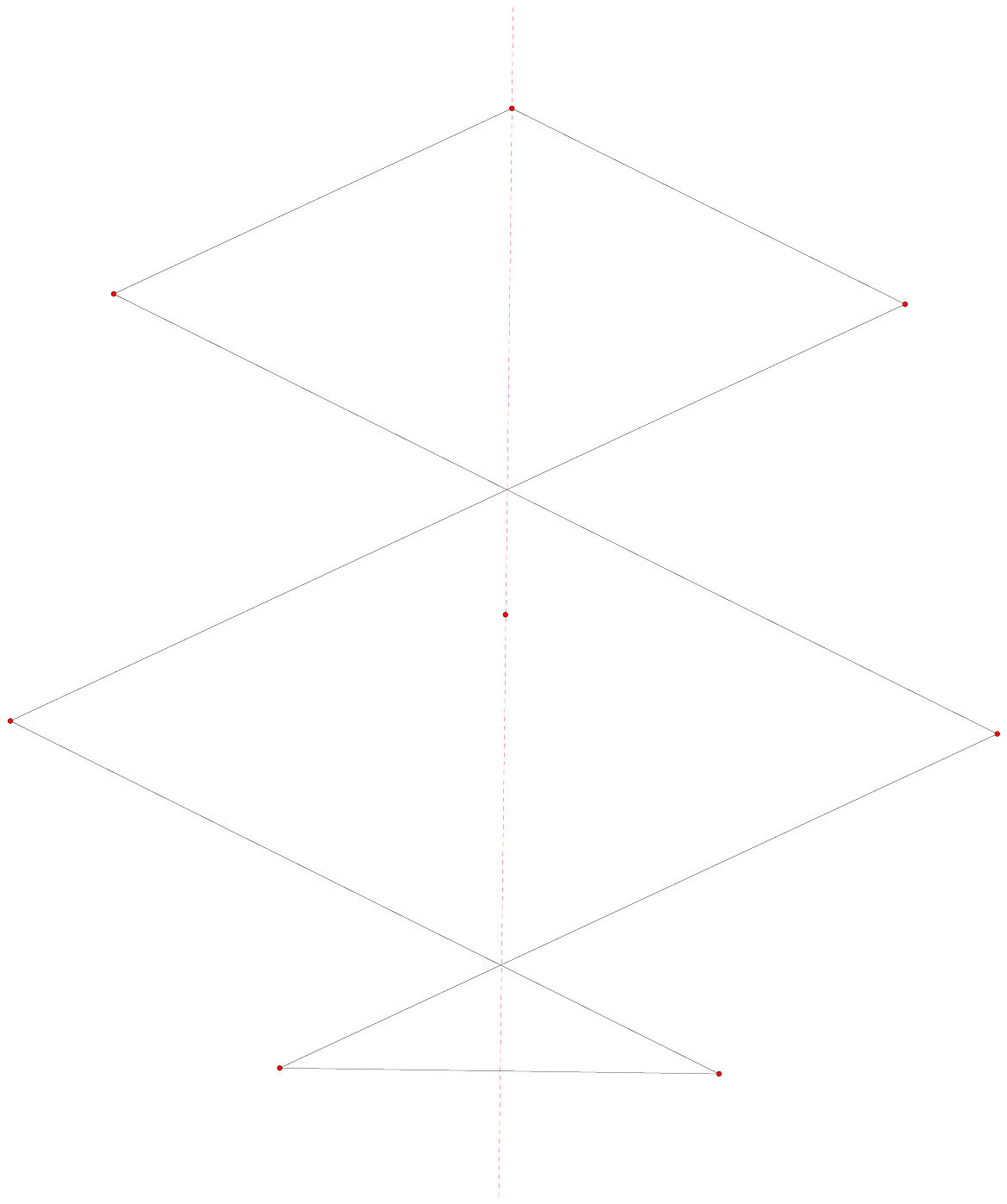}\\ \hline
\includegraphics[width=0.2\textwidth]{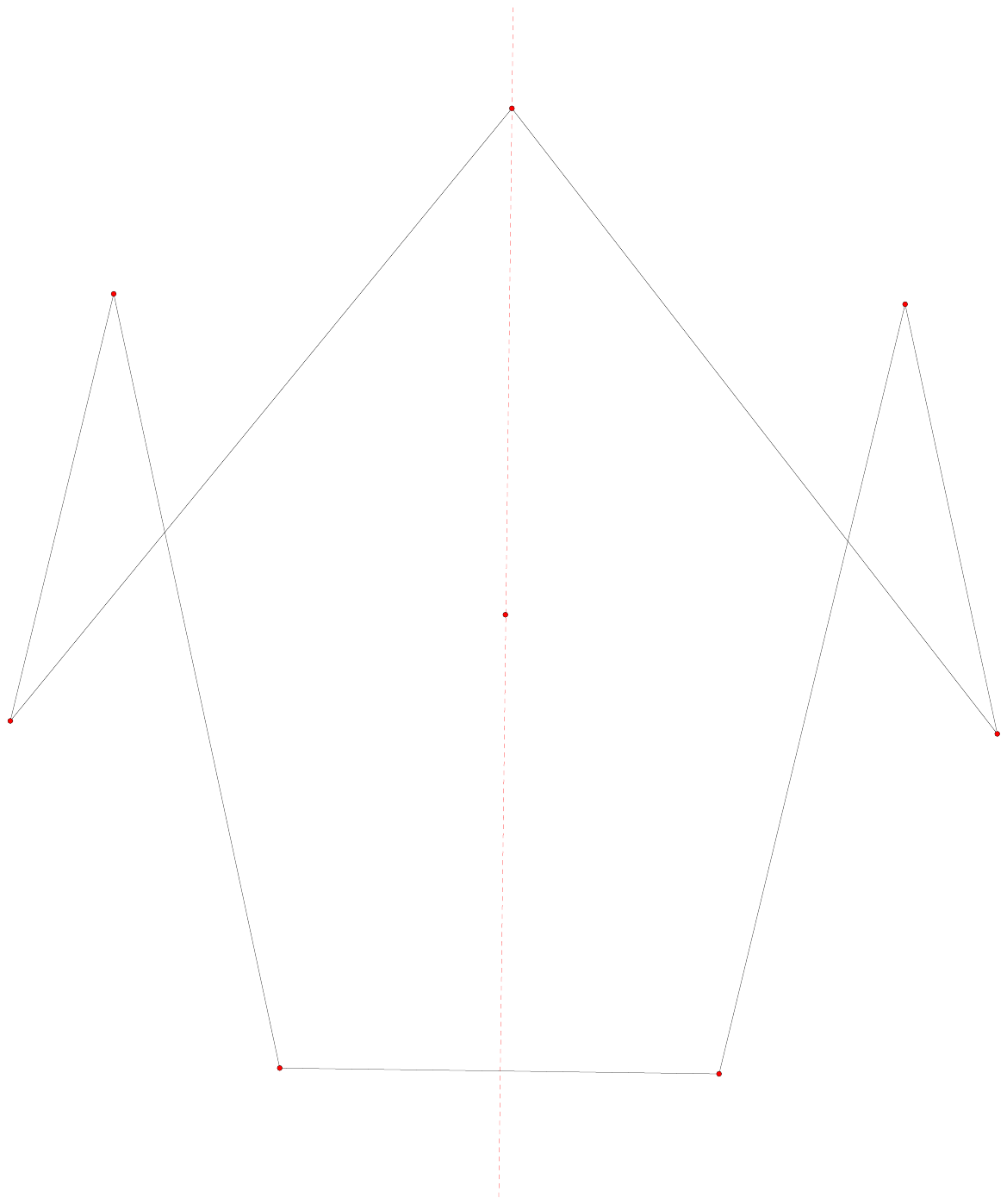} & \includegraphics[width=0.2\textwidth]{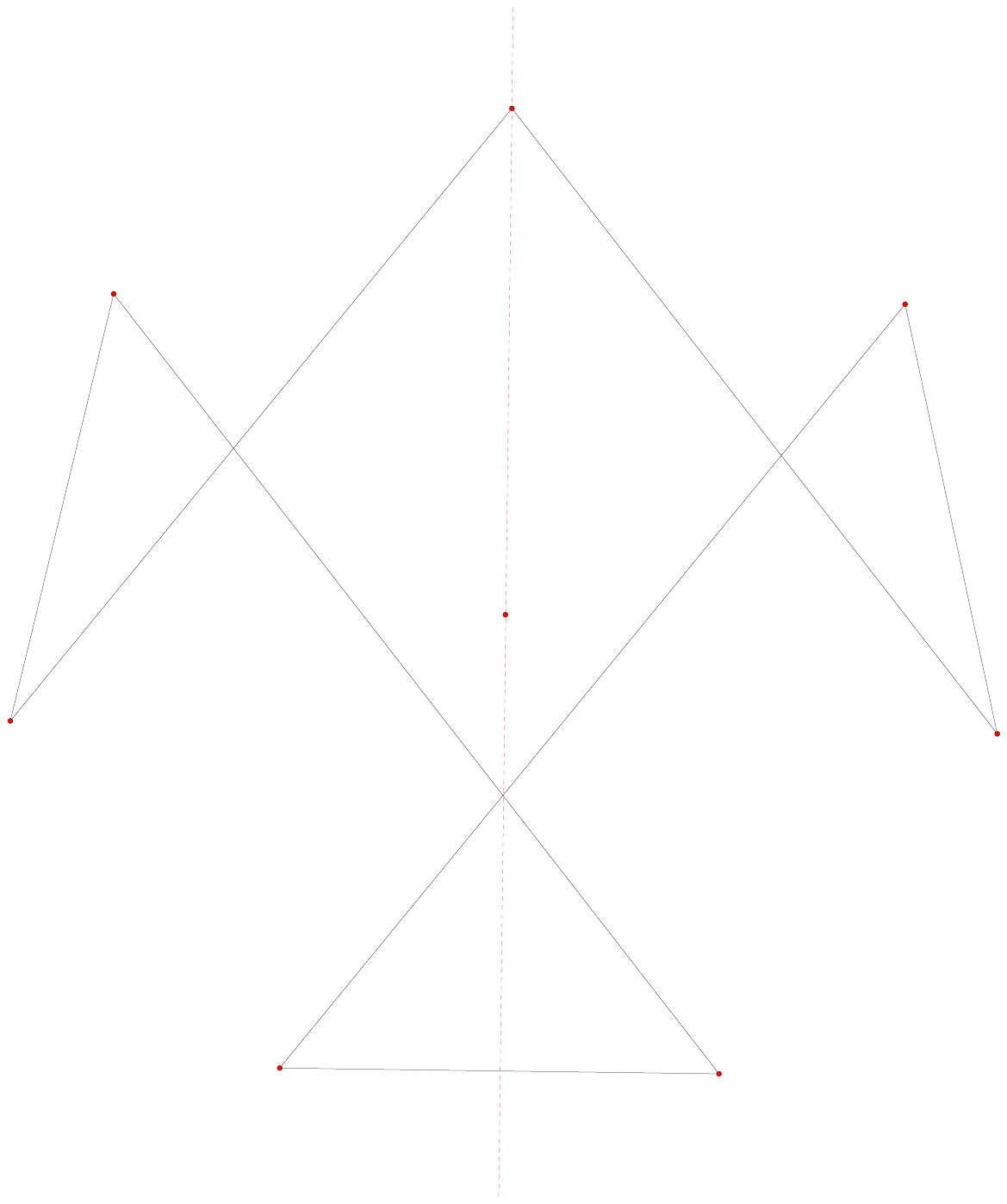} & \includegraphics[width=0.2\textwidth]{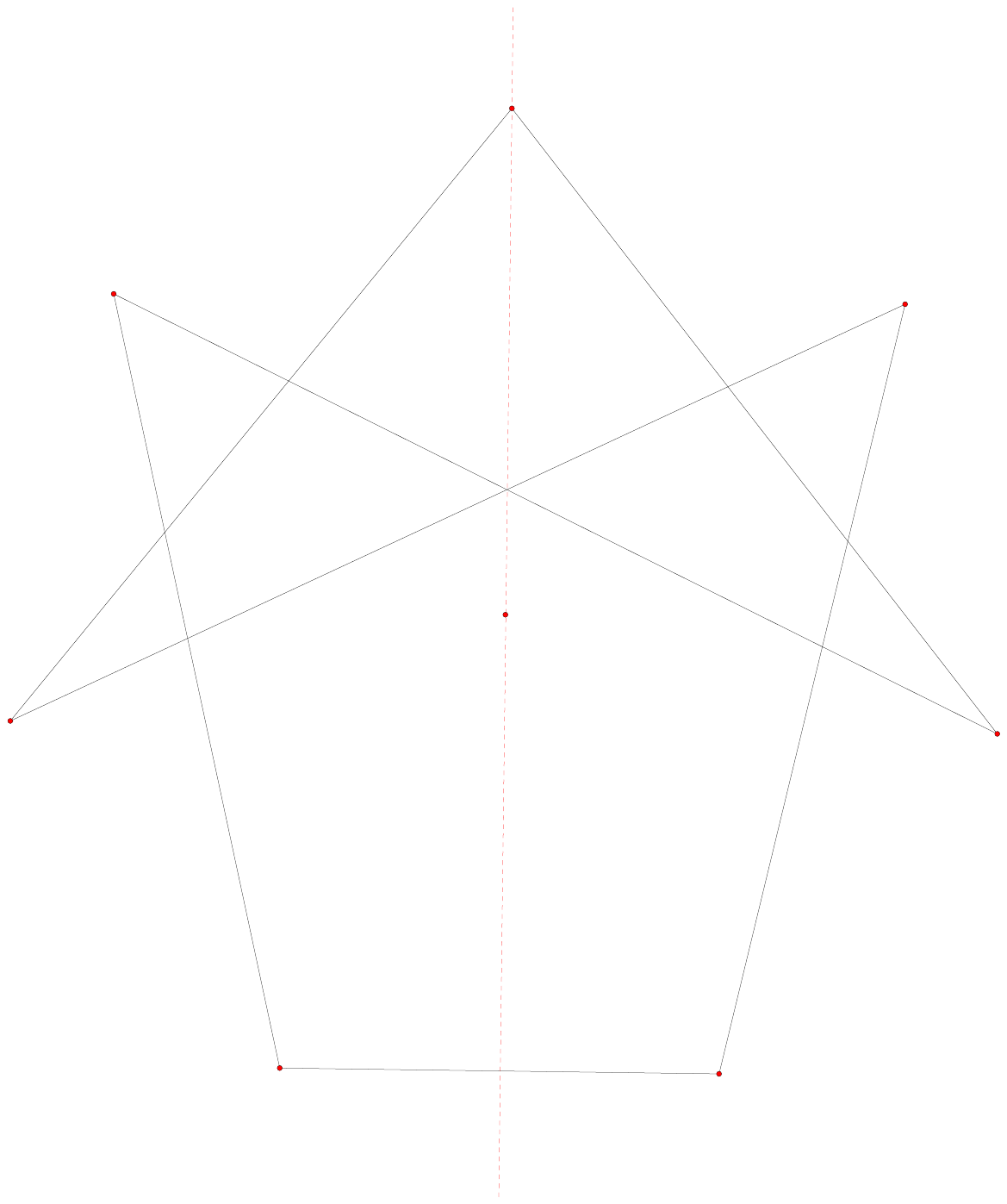} & 
\includegraphics[width=0.2\textwidth]{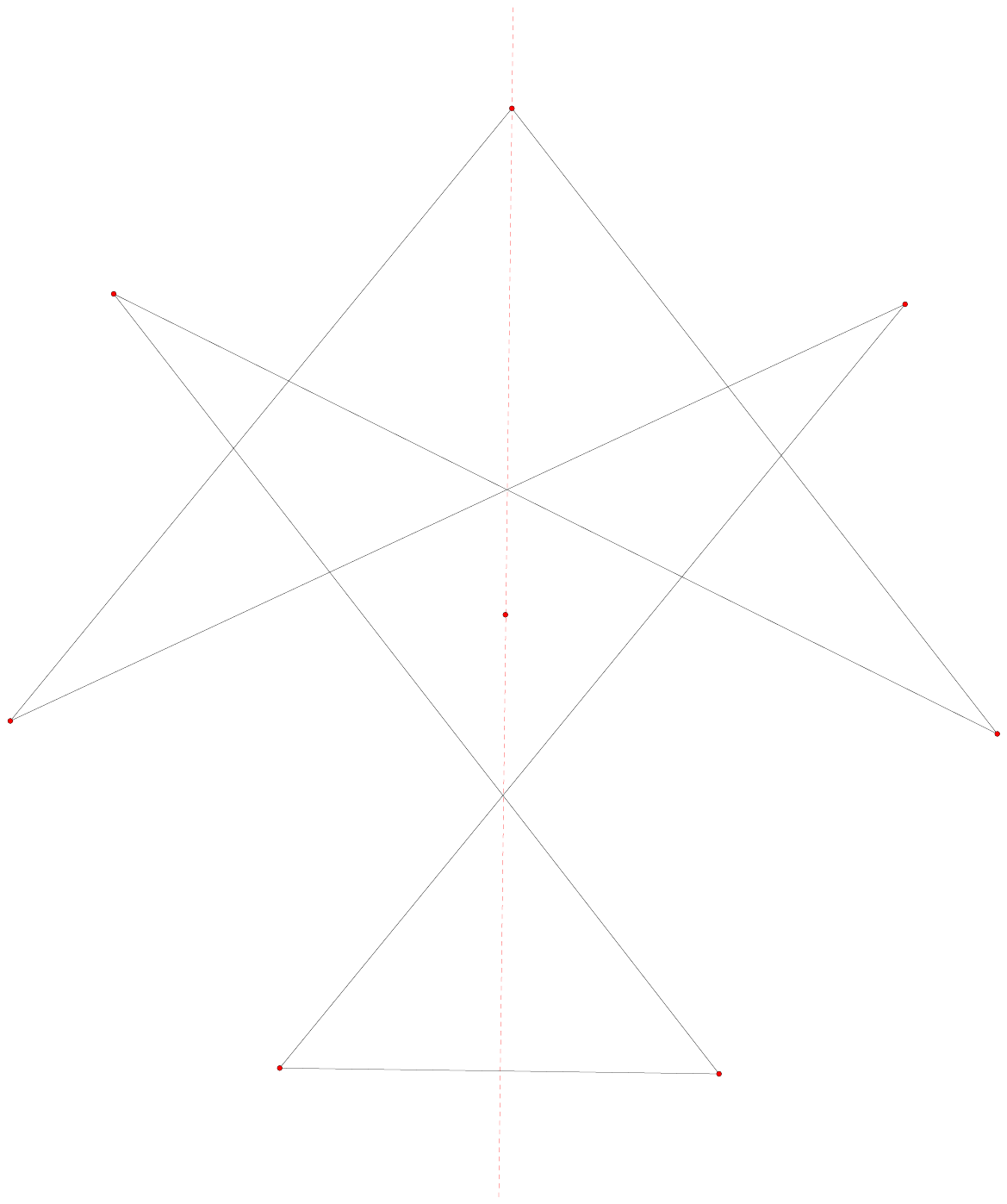}\\ \hline
\includegraphics[width=0.2\textwidth]{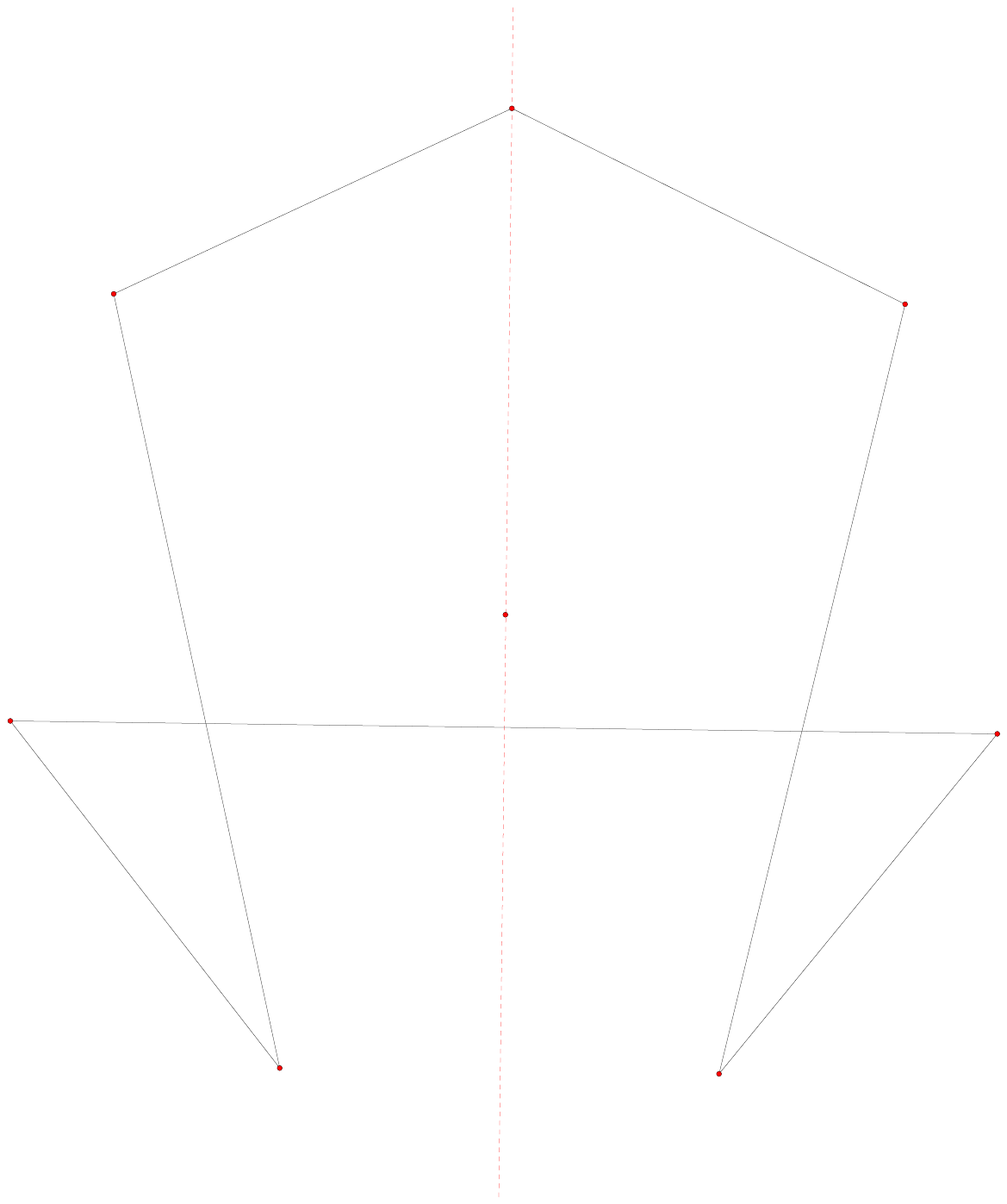} & \includegraphics[width=0.2\textwidth]{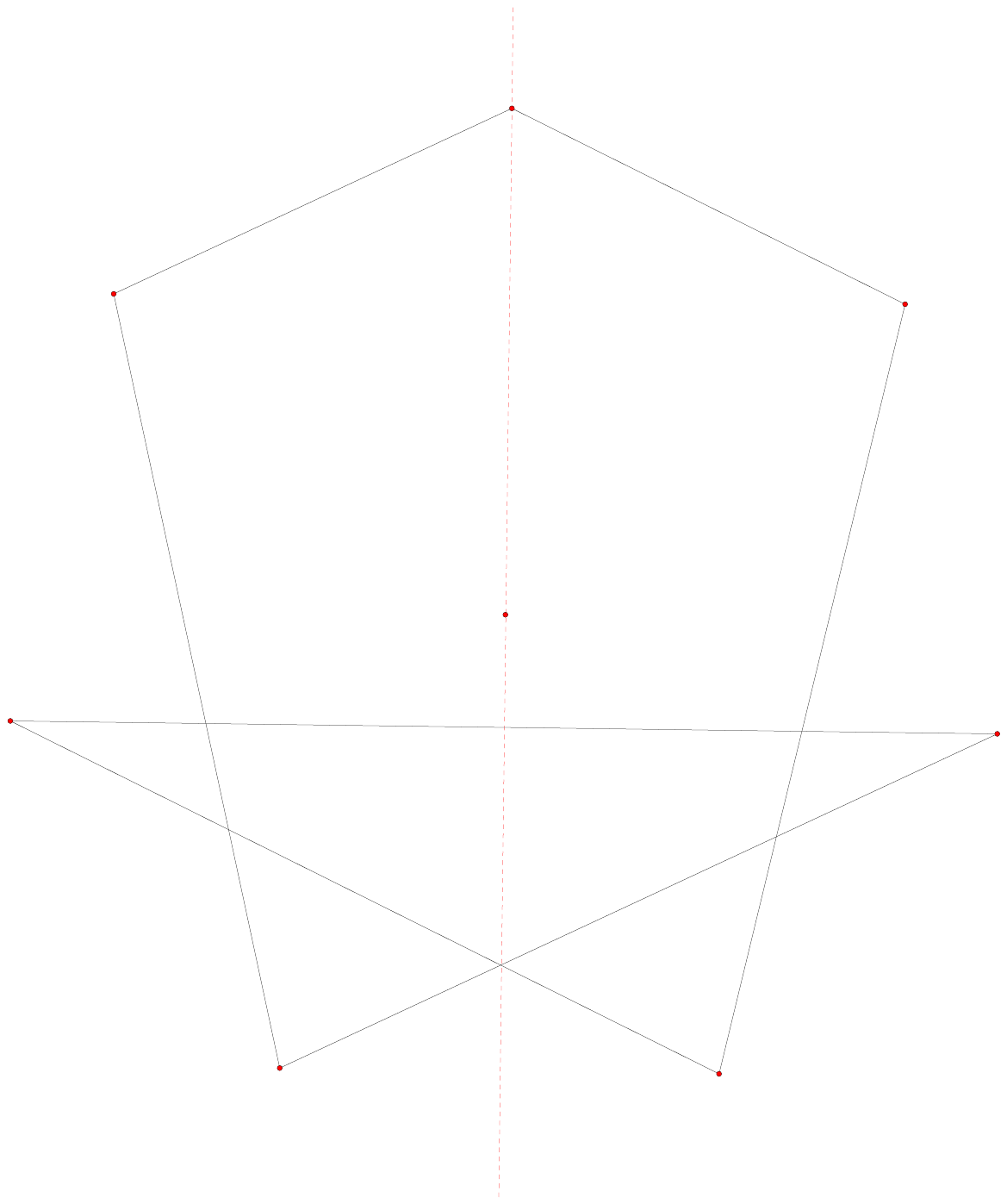} & \includegraphics[width=0.2\textwidth]{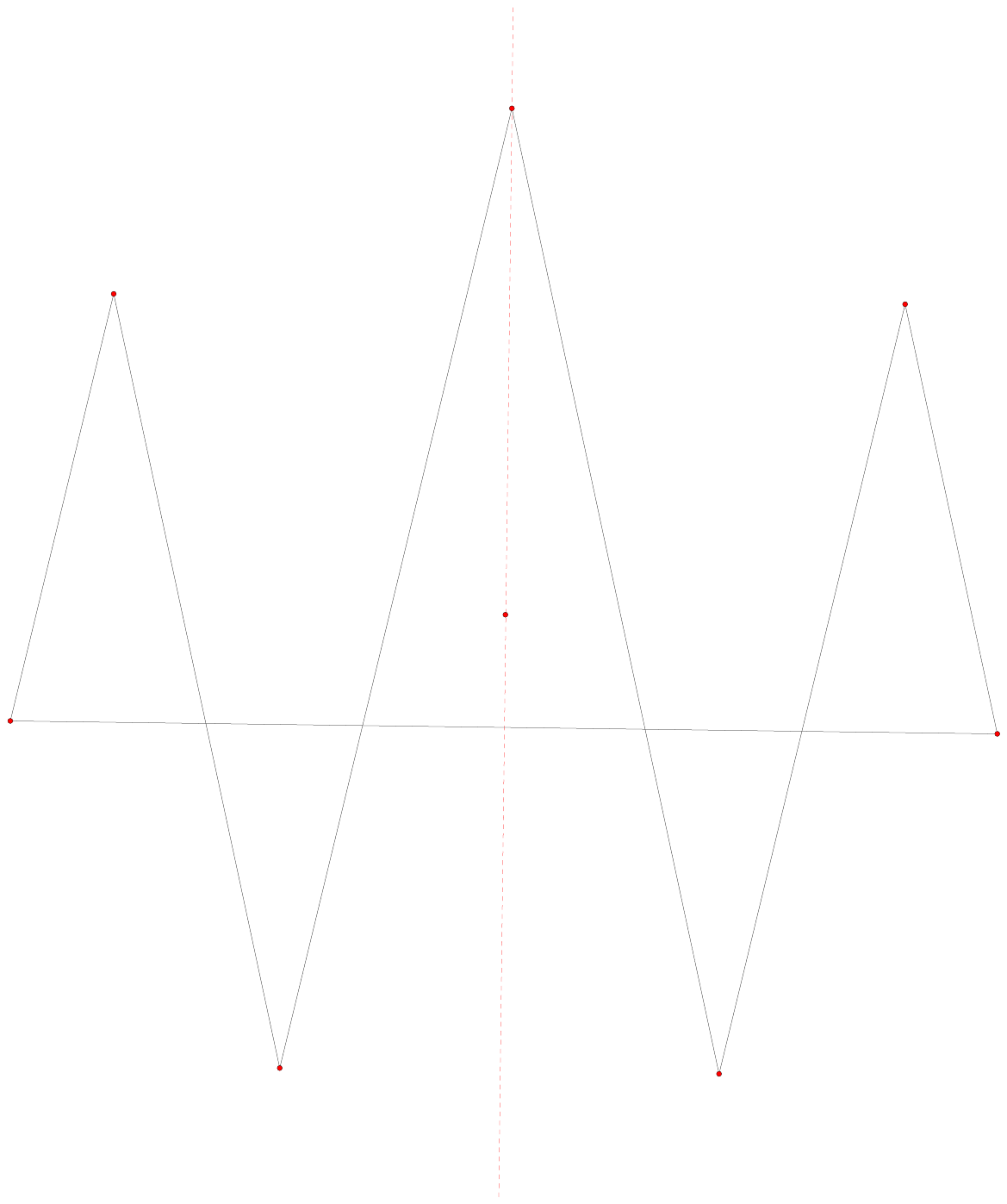} & \includegraphics[width=0.2\textwidth]{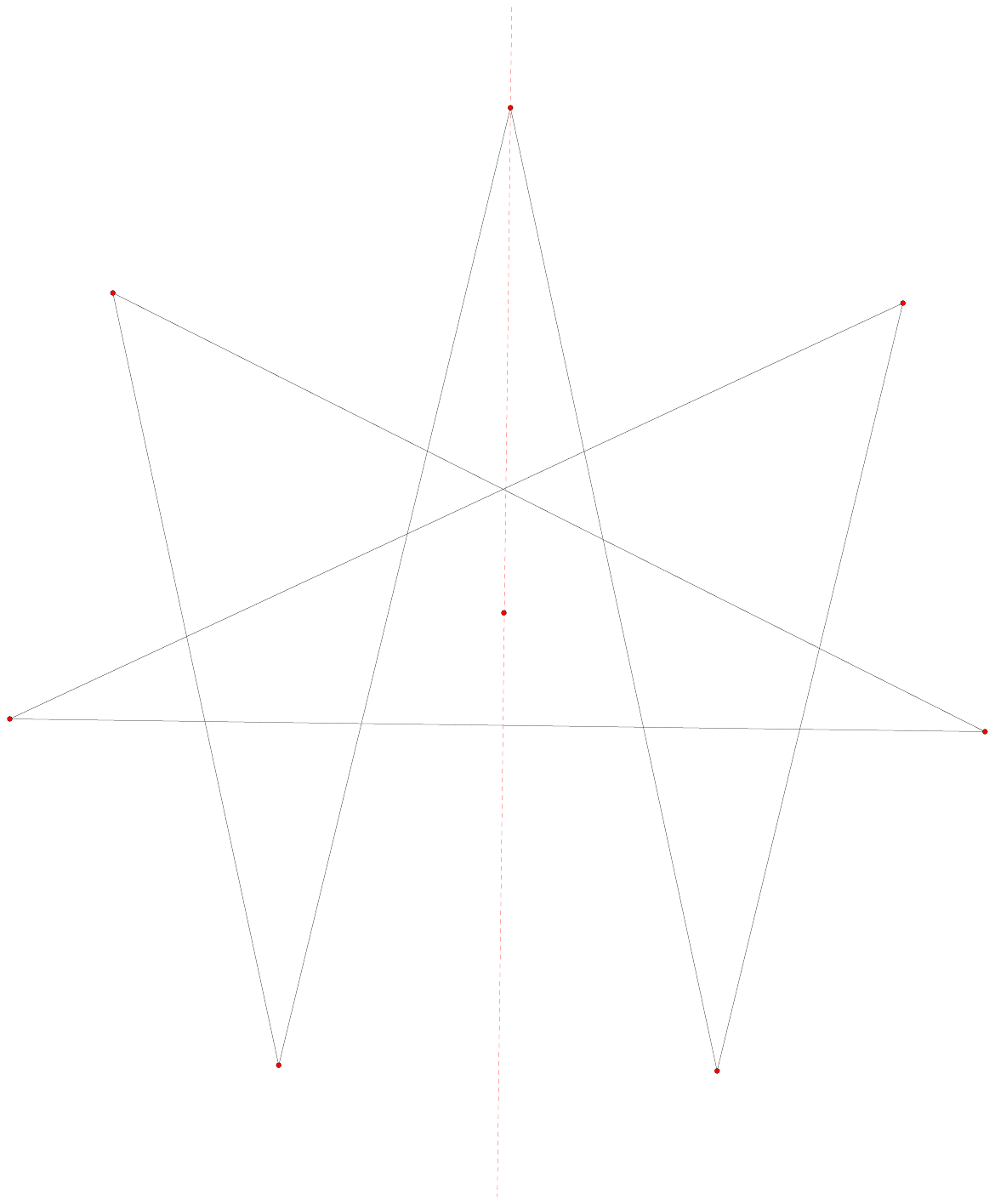}\\ \hline
\end{tabular}
\caption{representatifs of 12 equivalence-classes of the $7$-polygons with at least one axe}
\end{figure}
\end{center}

These symmetrical figures have one vertical drawn axe and represent 12 equivalence-classes.
\newpage
\subsection{Representation of the 24 symmetrical $7$-polygons (follower)}
\label{subsec:representation_of_the_24_symmetrical_7-polygons_follower}
\begin{center}
\begin{figure}[H]
\centering
\begin{tabular}{| c | c | c | c |}
\hline
\includegraphics[width=0.2\textwidth]{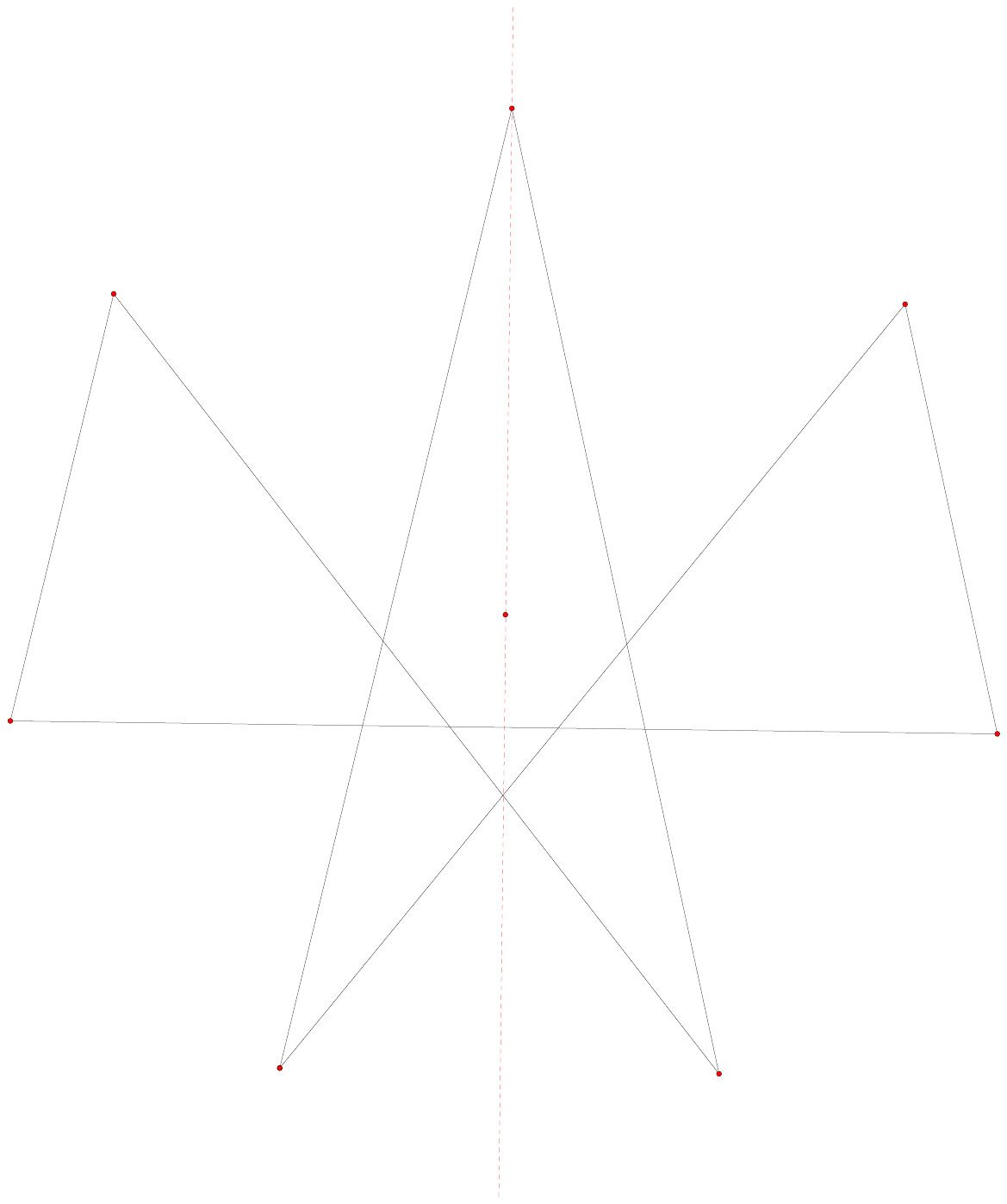} & \includegraphics[width=0.2\textwidth]{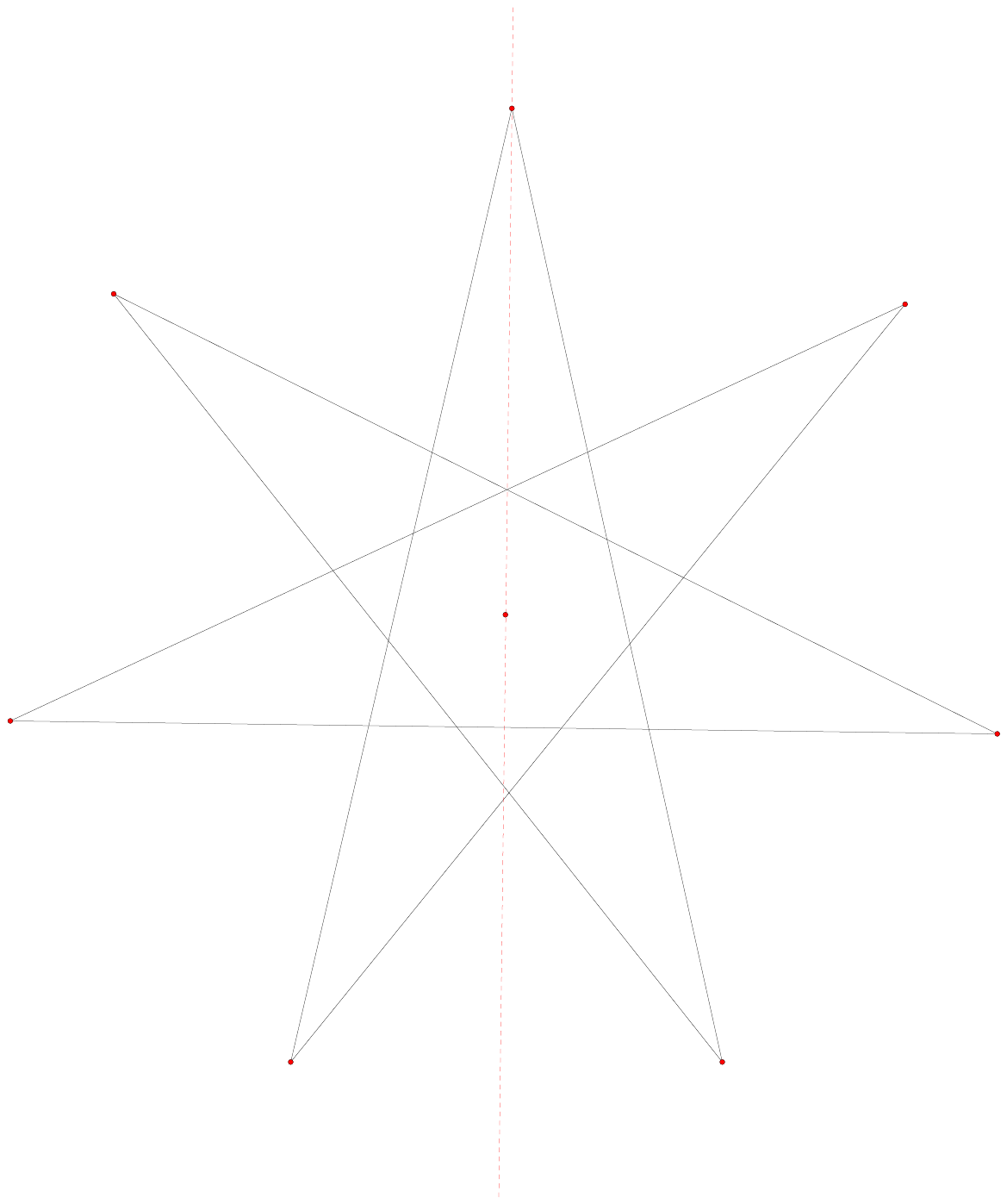} & \includegraphics[width=0.2\textwidth]{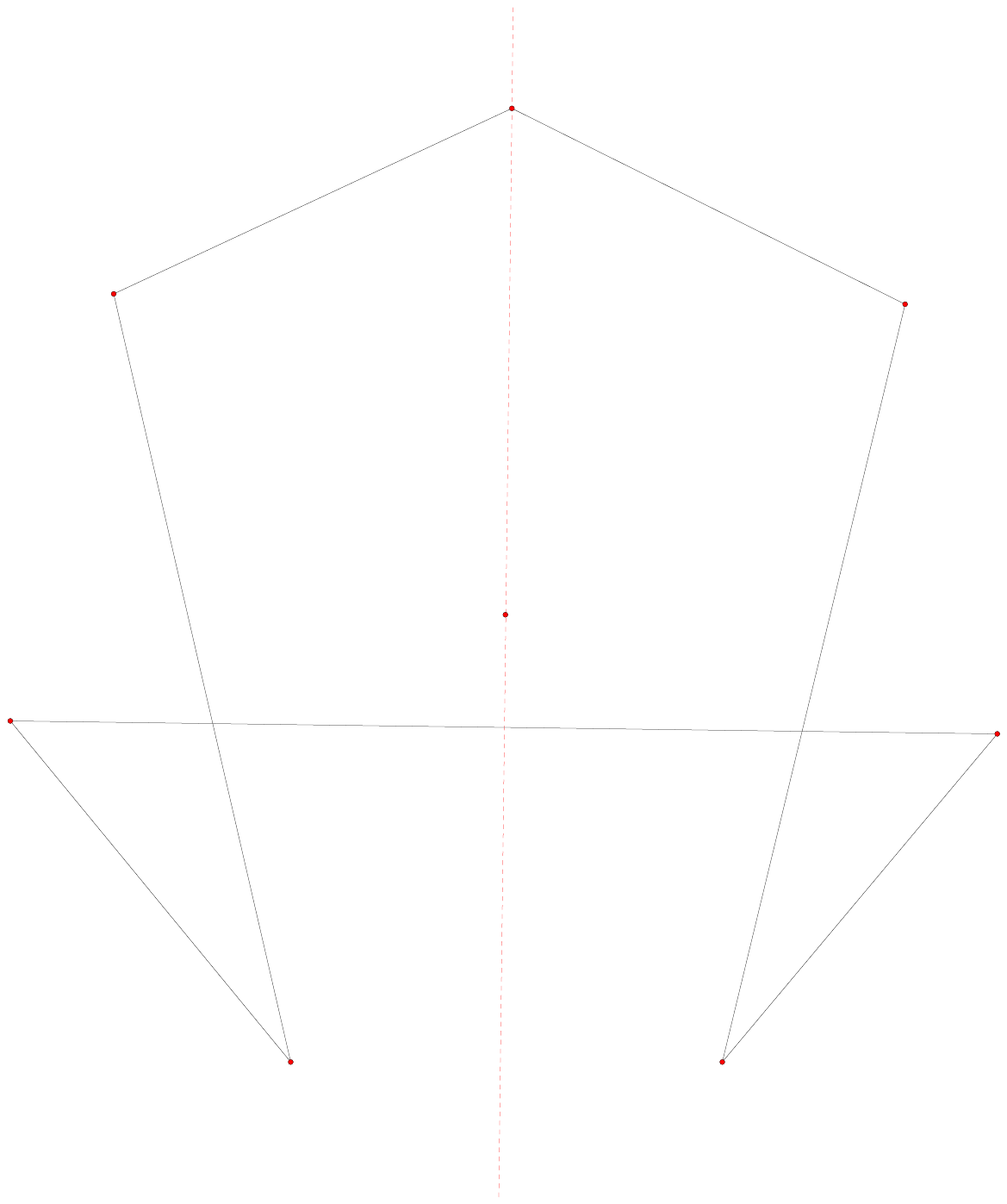} & \includegraphics[width=0.2\textwidth]{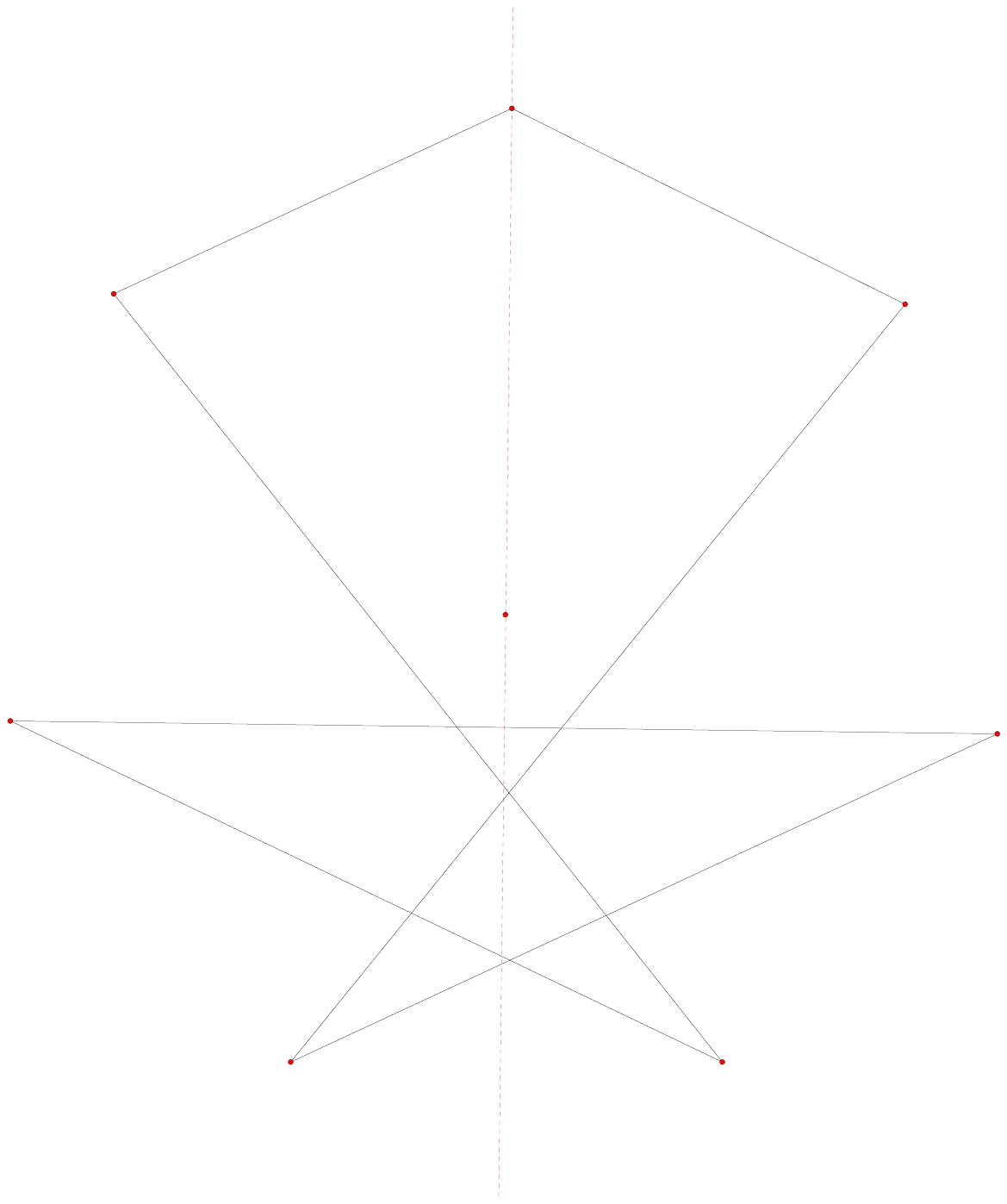}\\ \hline
\includegraphics[width=0.2\textwidth]{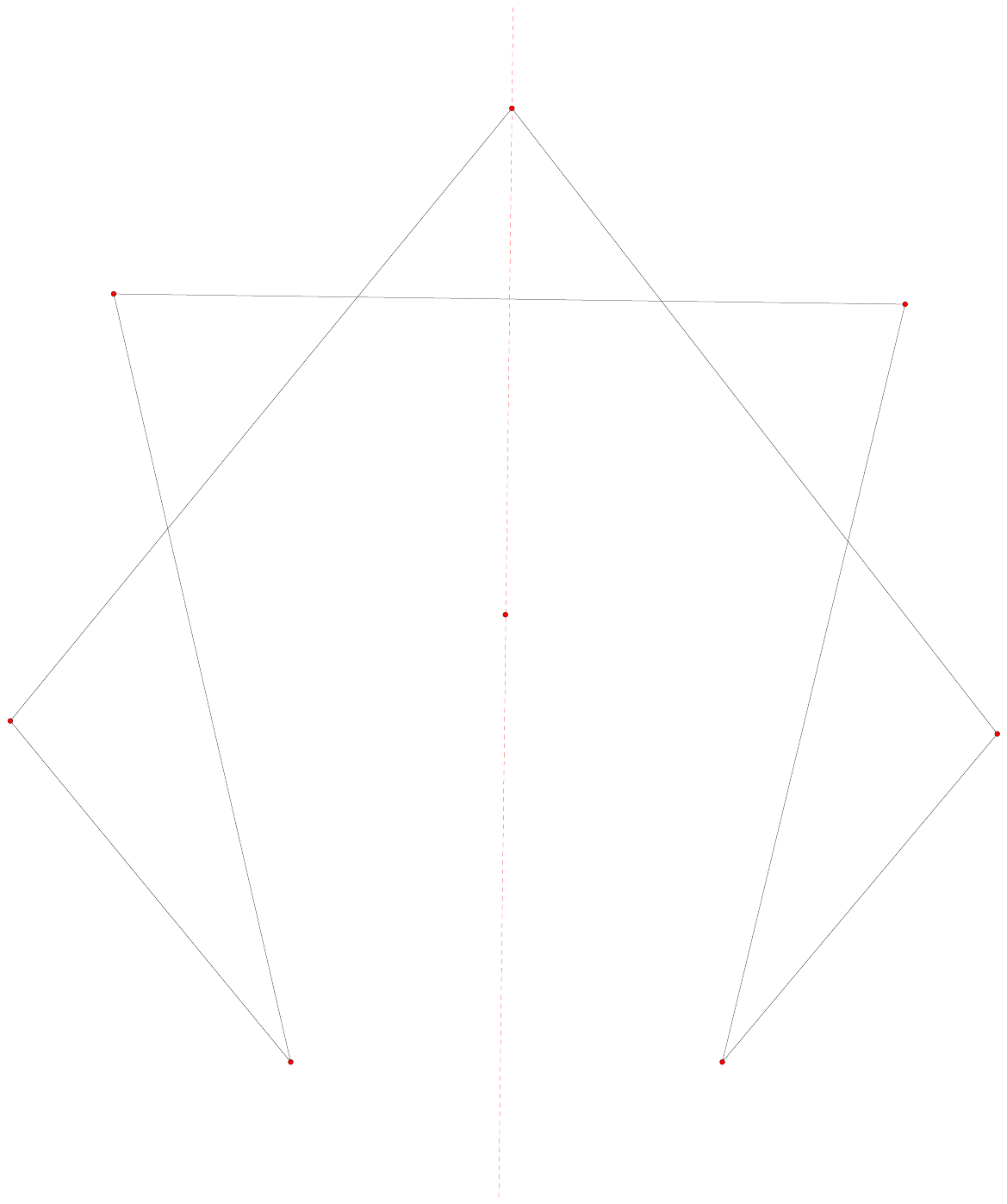} & \includegraphics[width=0.2\textwidth]{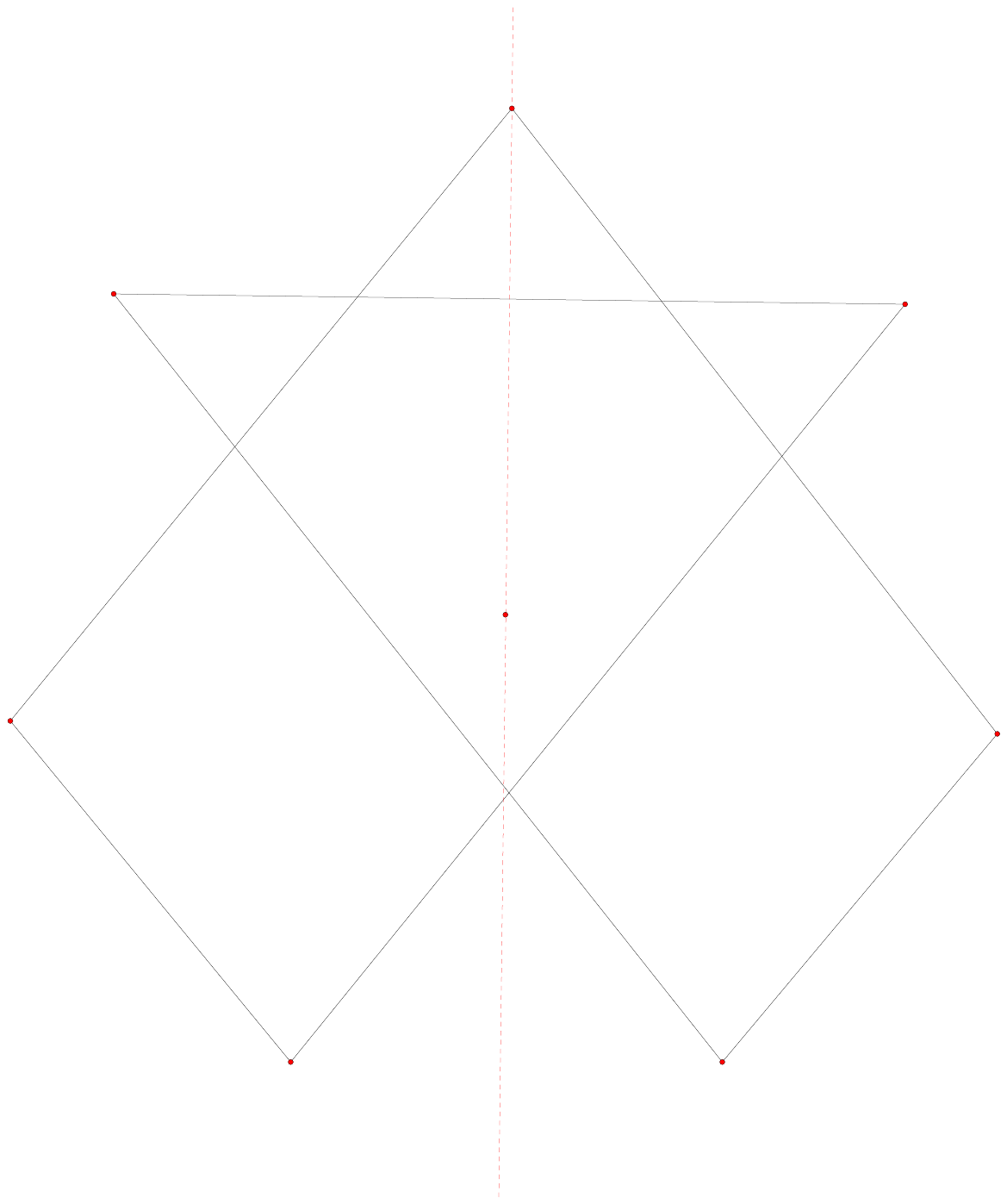} & \includegraphics[width=0.2\textwidth]{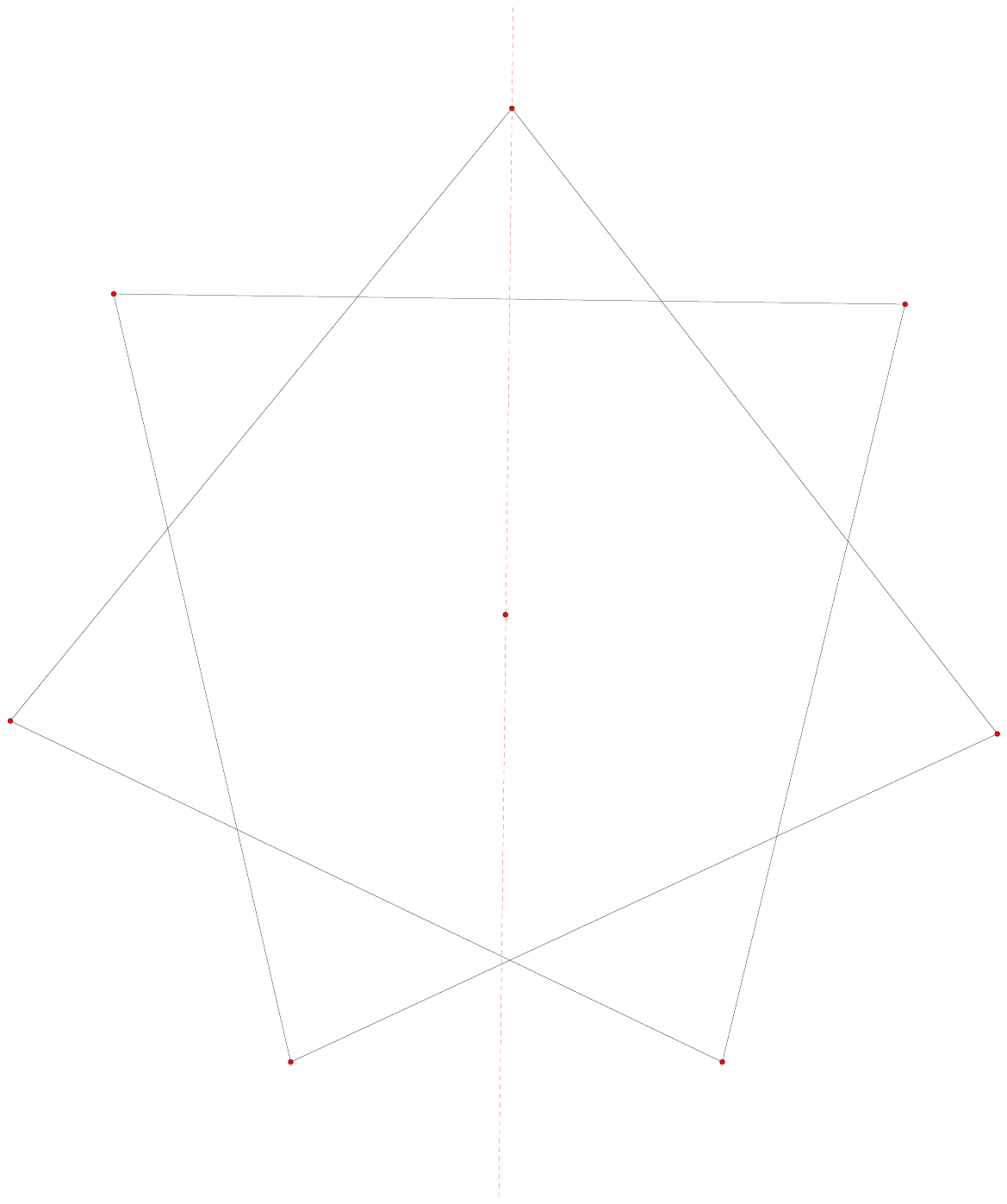} & \includegraphics[width=0.2\textwidth]{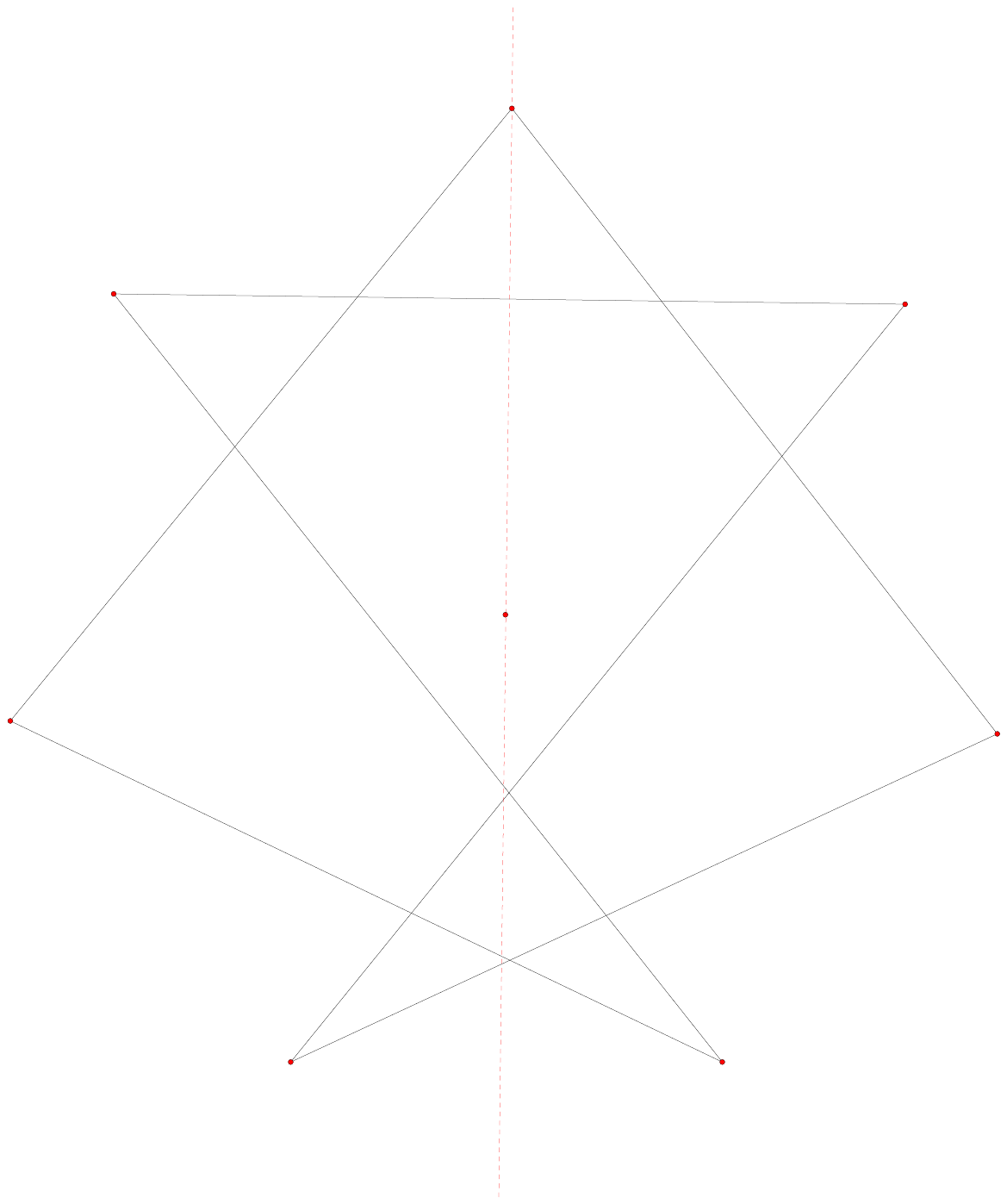}\\ \hline
\includegraphics[width=0.2\textwidth]{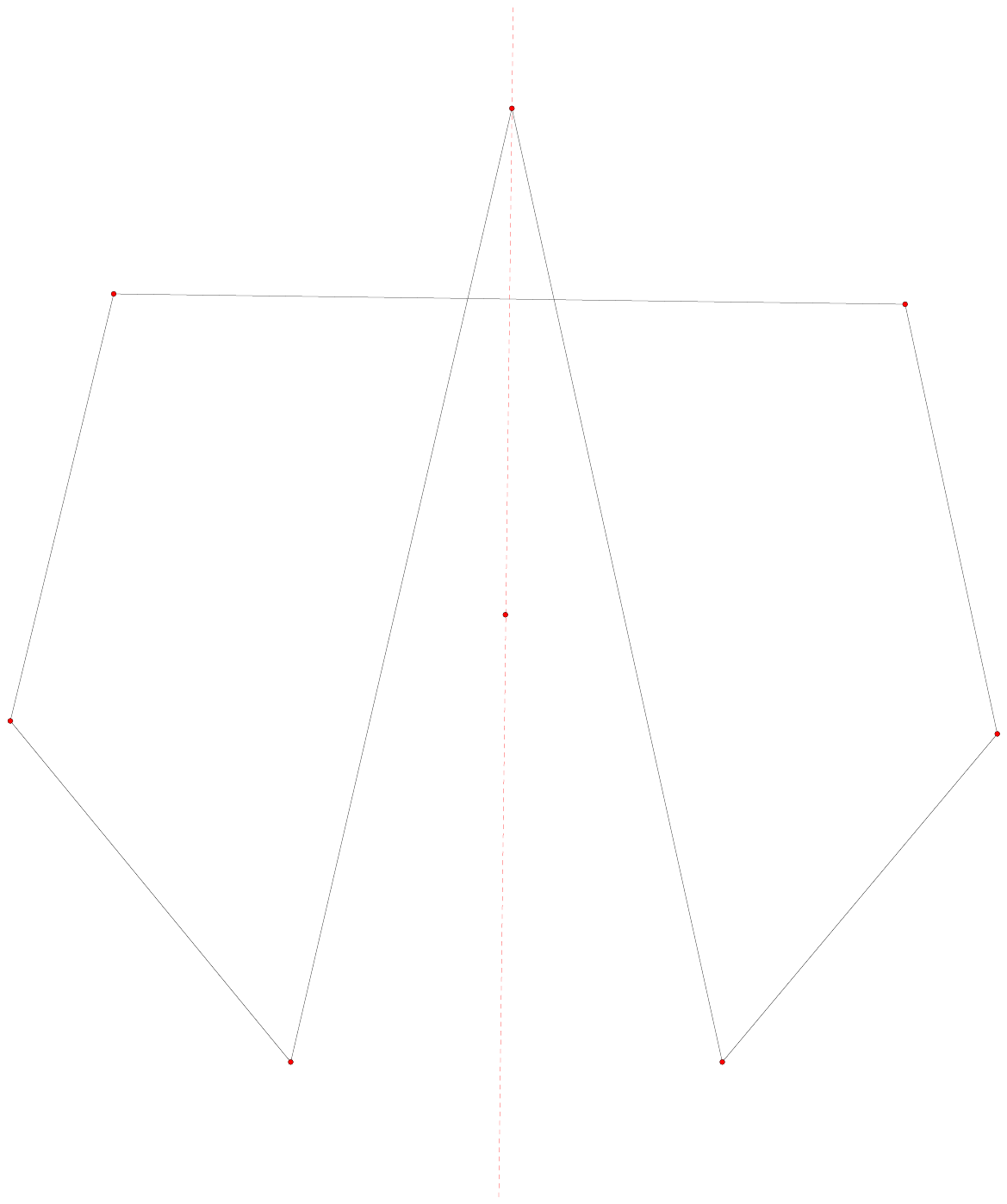} & \includegraphics[width=0.2\textwidth]{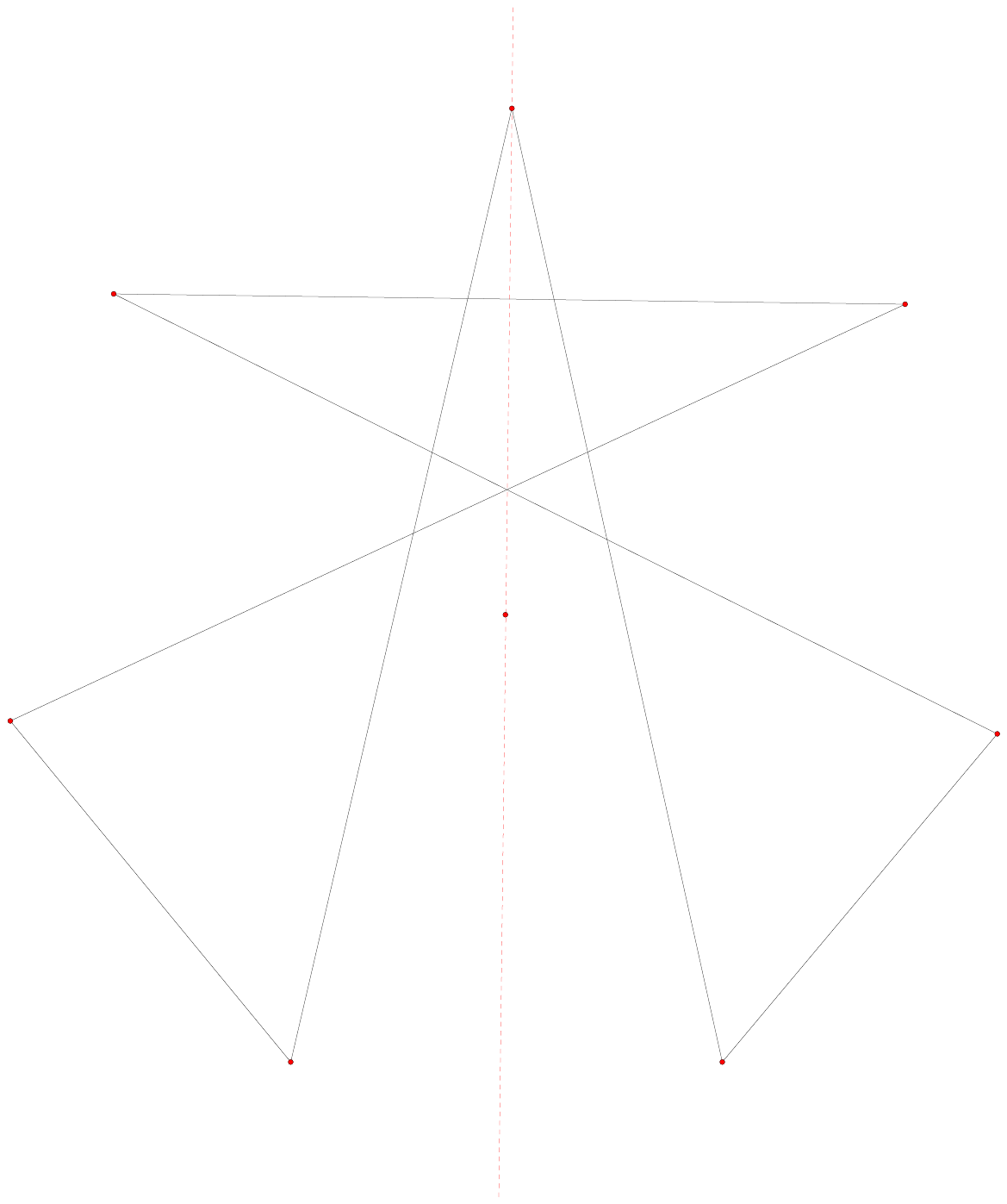} & \includegraphics[width=0.2\textwidth]{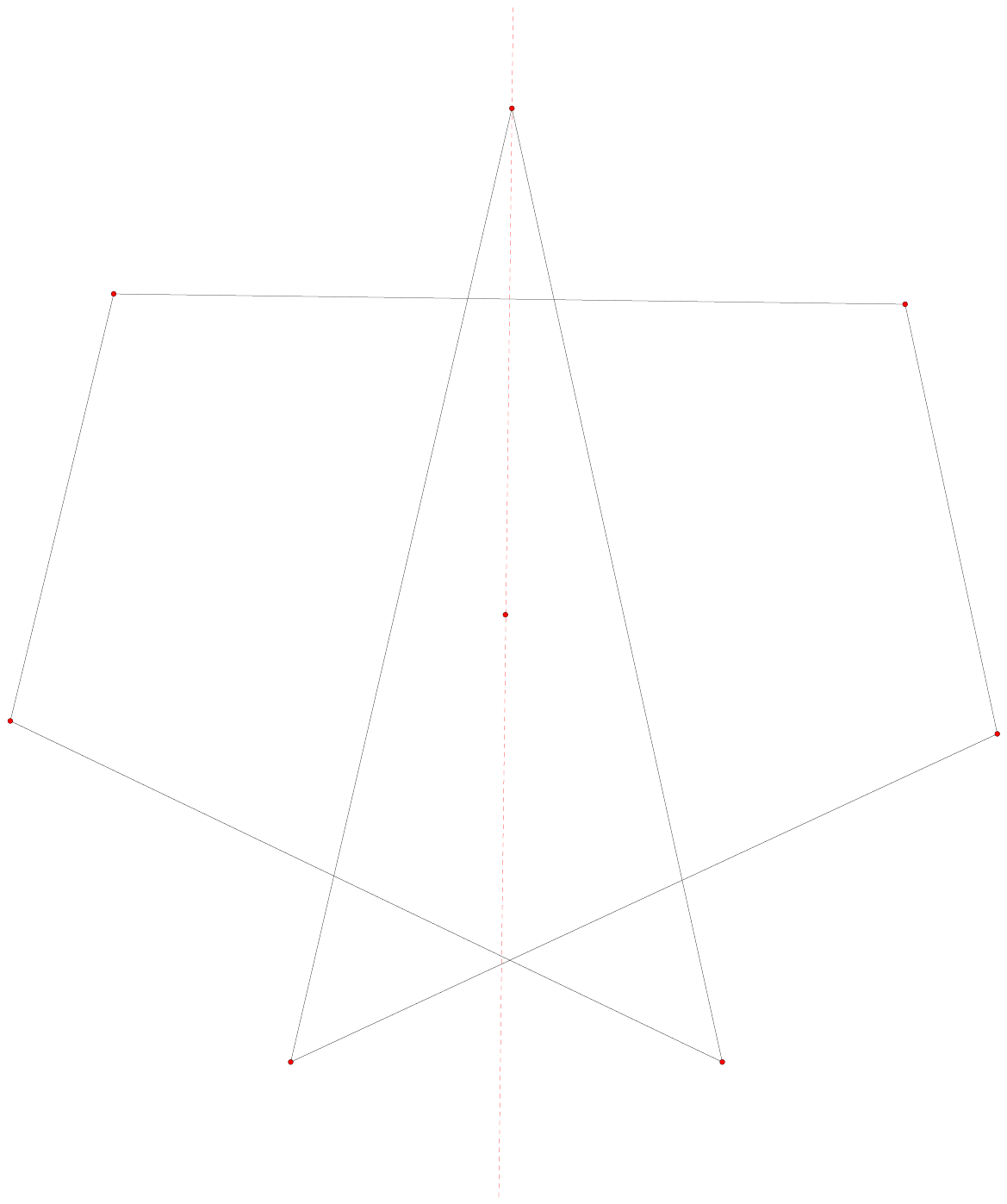} & \includegraphics[width=0.2\textwidth]{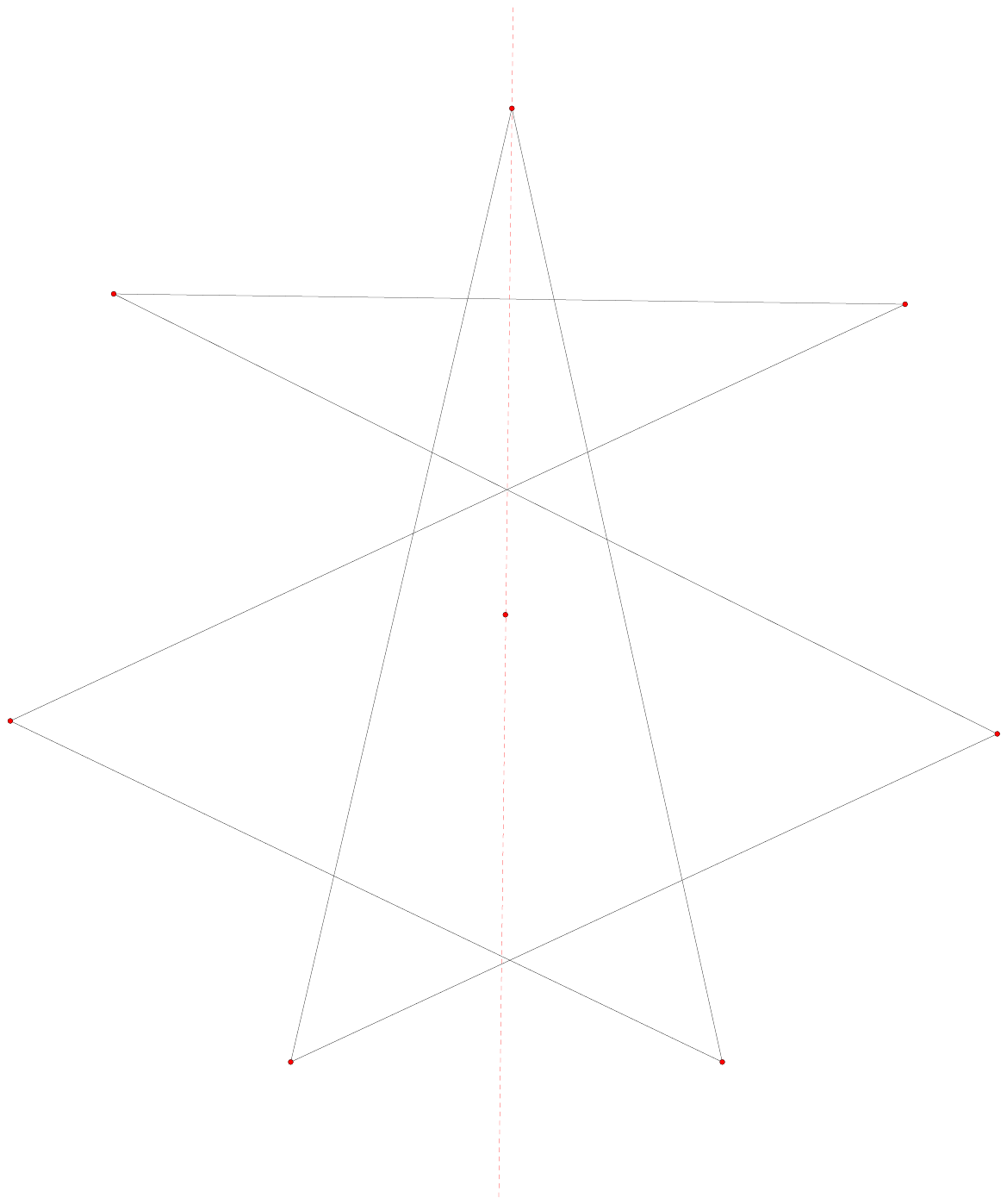}\\ \hline
\end{tabular}
\caption{representatifs of the other 12 equivalence-classes of the $7$-polygons with at least one axe}
\end{figure}
\end{center}

These symmetrical figures have also one vertical drawn axe and represent the other 12 equivalence-classes.
\newpage
\part{$p^{2}$-polygons}
\section{Question to deal with in this part of the article}
\label{sec:questions_to_deal_with_in_this_part_of_the_article}
Let $p \geq 3$ be a prime number. The number $p^{2}$ has the only divisors $p^{2}$ and $p$ and $1$. Therefore there exist only four different kinds of symmetry-categories for $p^{2}$-polygons:
\begin{enumerate}
\item[$\bullet$] equivalence-classes with $p^{2}$ symmetry-axes (socalled regular or star-polygons),
\item[$\bullet$]equivalence classes with $p$ symmetry-axes,
\item[$\bullet$]equivalence-classes with $1$ symmetry-axis,
\item[$\bullet$]equivalence-classes of the $p$-circular $p^{2}$-polygons, which by definition have no symmetry-axis, 
\item[$\bullet$]equivalence-classes of the $p^{2}$-polygons with no symmetry-axis and no circularity, i.e. the completely asymmetric polygons.
\end{enumerate}
We define and will calculate the following numbers: 
\begin{enumerate}
\item[$\bullet$] $\vert X(p^{2})\vert$ of all equivalence-classes,
\item[$\bullet$] $\vert X_{p^{2}}(p^{2}) \vert$ of  the equivalence-classes with $p^{2}$ symmetry-axes,
\item[$\bullet$] $ \vert X_p(p^{2})\vert$ of the equivalence-classes with $p$ symmetry-axes,
\item[$\bullet$] $ \vert X_1(p^{2})\vert$ of  the equivalence-classes with $1$ symmetry-axis,
\item[$\bullet$] $ \vert X_c(p^{2})\vert$ of  the equivalence-classes of the $p$-circular $p^{2}$-polygons which have no symmetry-axis,
\item[$\bullet$] $ \vert X_a(p^{2})\vert$ of  the equivalence-classes of the non $p$-circular $p^{2}$-polygons with no symmetry-axis, the completely asymmetrical $p^{2}$-polygons.
\end{enumerate}
In addition, we give a table of results and illustrate them by some sets of representatifs for the prime number $p=3$.
\section{Results}
\label{sec:results}
\subsection{Preparations}
\label{subsec:preparations}
Besides the numbers, which were defined above and which are the goals of these calculations, we need  the following numbers during the proofs:
\begin{enumerate}
\item[$\bullet$]$\vert X_{p+}(p^{2}) \vert$ is the number of equivalence-classes of $p^{2}$-polygons with at least $p$ symmetry-axes,
\item[$\bullet$]$\vert X_{1+}(p^{2}) \vert$ is the number of equivalence-classes of $p^{2}$-polygons with at least $1$ symmetry-axis.
\end{enumerate}
\subsection{Main sentence of part II}
\label{subsec:main_sentence_of_part_II}
Let be $p\geq3$ a prim number. The number of the equivalence-classes of $p^{2}$-polygons are:
\begin{itemize}
\item[a.)]{$ \vert X(p^{2}) \vert = \dfrac{1}{2 \cdot p^{4}} \cdot \left[p^{4} \cdot (p-1)^{2}+(p-1)^{2} \cdot p! \cdot p^{p}+(p^{2})!    \right]$ as the number of all equivalence-classes,}
\item[b.)]{$\vert X_{p^{2}}(p^{2}) \vert=\dfrac{p\cdot(p-1)}{2}$ as the number of equivalence-classes of the $p^{2}$-polygons with $p^{2}$ axes,}
\item[c.)]{$\vert X_p(p^{2} \vert = \dfrac{p-1}{2} \cdot \left[(2p)^{\dfrac{p-1}{2}} \cdot \left(\dfrac{p-1}{2}\right)!-p \right]$ as the number of equivalence-classes of $p^{2}$-polygons with exact $p$ symmetry-axes,}
\item[d.)]{$\vert X_1(p^{2} \vert = \dfrac{1}{2}\cdot \left[\left(\dfrac{p^{2}-1}{2}\right)! \cdot 2^{\dfrac{p^{2}-1}{2}}-(2p)^{\dfrac{p-1}{2}} \cdot (p-1)\cdot \left( \dfrac{p-1}{2}\right)!\right]$ as the number of equivalence-classes with exactly $1$ symmetry-axis,}
\item[e.)]{$\vert X_c(p^{2} \vert  = \dfrac{p-1}{2}\cdot \left[(p-1)!\cdot p^{p-2}-(2p)^{\dfrac{p-1}{2}}\cdot\left(\dfrac{p-1}{2}\right)!+p-1\right]$ as the number of equivalence-classes of $p$-circular $p^{2}$-polygons,}
\item[f.)]{$\vert X_a(p^{2}) \vert =\dfrac{1}{2 \cdot p^{4}} \cdot \left[p^{4} \cdot (p-1)^{2}+(p-1)^{2} \cdot p! \cdot p^{p}+(p^{2})!    \right]-\dfrac{1}{2}\cdot\left( \dfrac{p^{2}-1}{2}\right)!\cdot 2^{\dfrac{p^{2}-1}{2}}-\dfrac{p-1}{2}\cdot \left[((p-1)!)\cdot p^{p-2}-\left(\dfrac{p-1}{2}\right)!\cdot\left(2p\right)^{\dfrac{p-1}{2}}+p-1\right] $ as the number of equivalence-classes of asymmetrical and not circular $p^{2}$-polygons.}
\end{itemize}
\section{Proofs of part II}
\label{sec:proofs_of_part_II}

\subsection{Main sentence a.)}
\label{subsec:main_sentence_a.)}
\begin{proof}
The number of equivalence-classes of all $p^{2}$-polygons follows from the formula, given by S. W. Golomb and L.R.Welch. \cite{Golomb1960} They proved there the formula for the number $\vert X(n) \vert$ of equivalence-classes of $n$-polygons, if $n$ is an odd number:
$\vert X(n) \vert=\dfrac{1}{2n^2}\left(   \sum \limits_{d \mid n}  \varphi^2 \left( \dfrac{n}{d} \right)\cdot d! \cdot \left( \frac{n}{d} \right)^d \right)$, where $d$ is a divisor of $n$ and $\varphi(n)$ denotes the totient function of $n$.
As $p \geq 3$ is a prime number and $n=p^{2}$ is an odd number with the three divisors $p^{2},p$ and $1$ and we get immediately:
\begin{center}
\fbox{$\vert X(p^{2}) \vert = \dfrac{1}{2 \cdot p^{4}} \cdot \left[p^{4} \cdot (p-1)^{2}+(p-1)^{2} \cdot p! \cdot p^{p}+(p^{2})!    \right]$}
\end{center}
\end{proof}
\subsection{Main sentence b.)}
\label{subsec:main_sentence_b.)}
\begin{proof}
The number of equivalence-classes of the $n$-polygons with $n$ symmetry-axes is well known by Leonard Euler and others and mentioned by H.S.M. Coxeter \cite{Coxeter}: $\vert X_n(n)\vert=\dfrac{\varphi(n)}{2}$. As $\varphi(p^{2})=p^{2} \cdot \left(1-\dfrac{1}{p} \right)=p\cdot (p-1)$, we conclude:
\begin{center}
\fbox{$\vert X_{p^{2}}(p^{2}) \vert=\dfrac{p \cdot(p-1)}{2}$}
\end{center}
\end{proof}
\newpage
\subsection{Main sentence c.) $p$ axes}
\label{subsec:main_sentence_c.)_p_axes}
\begin{proof}
To develop this proof we need construction figures for $p=5$ with $p^{2}= 25$ regularly distributed vertices and regularly distributed axes to make it more easy to follow the construction. One of the $5$ axes is set vertically, the others are drawn after rotations by the angle $\dfrac{360}{p}$.

\begin{figure}[H]                                
\begin{tabular}{| c | c |}
\hline\
sketch 1 & sketch 2 \\
\includegraphics[width=0.5\textwidth]{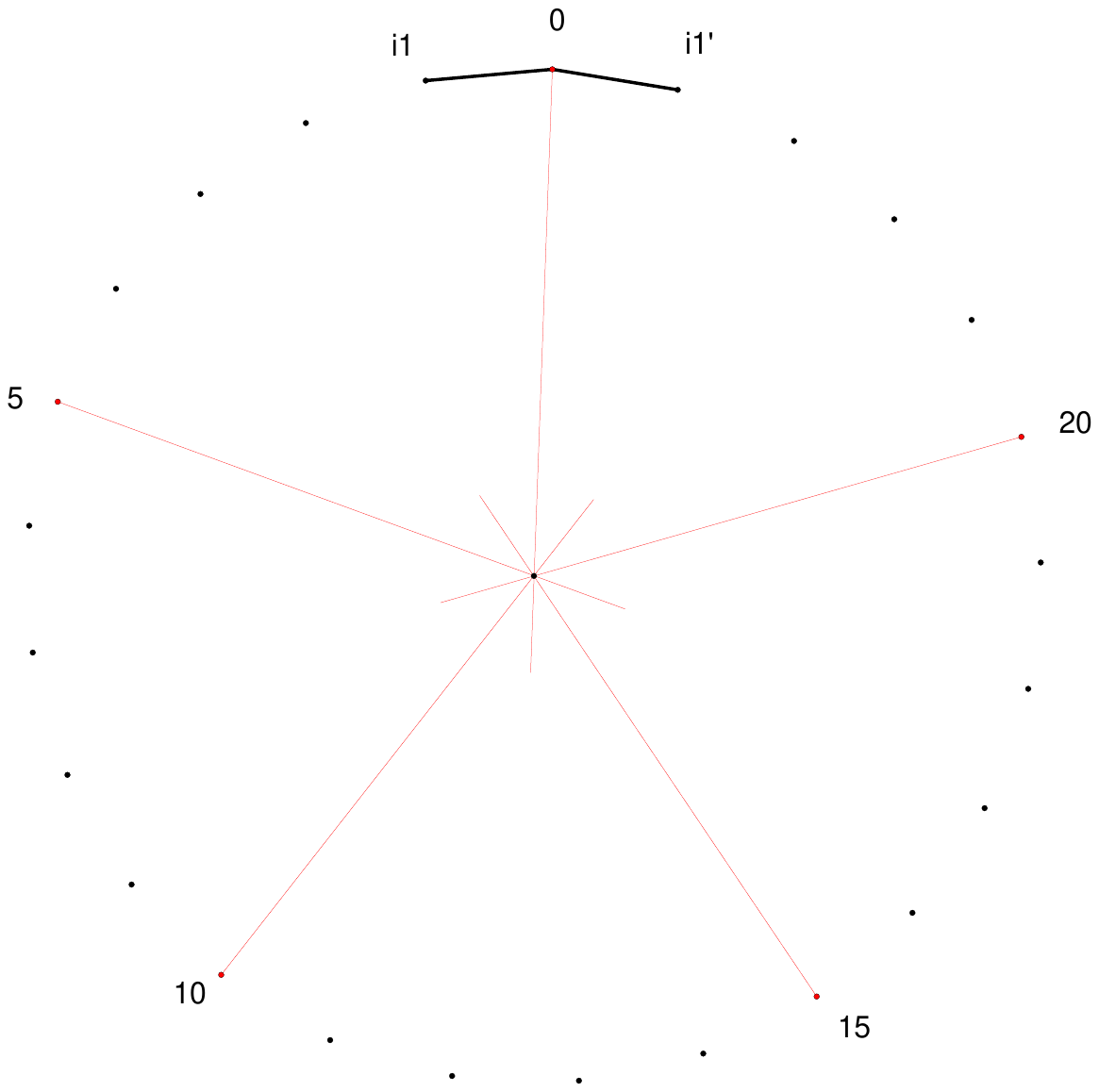} &
\includegraphics[width=0.5\textwidth]{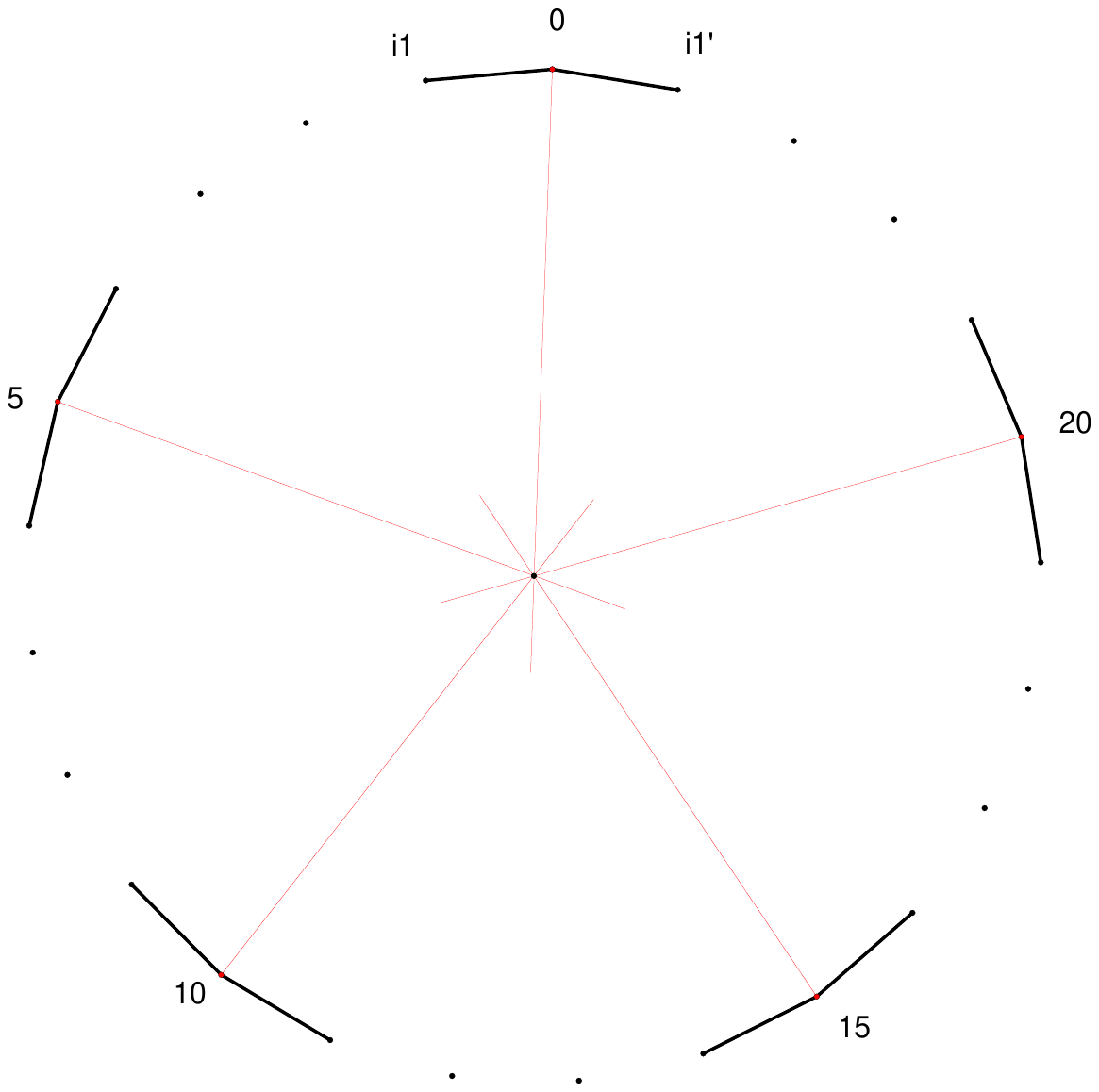} \\ \hline
\end{tabular}
\caption{construction-figures 1 and 2 for the equivalence-classes of the $p^{2}$-polygons with at least $p$ axes}
\end{figure}
Sketch 1 and 2: We choose a pair of vertices $v_{i1}$ and $v_{i1'}=v_{p^{2}-i1}$, which is symmetric to the vertical axis. And we connect them with vertex $v_0$ and get a symmetric open chain $O(3)=\overline{v_{i1'} v_0 v_{i1}}$. We have $\dfrac{p^{2}-p}{2}=\dfrac{p \cdot (p-1)}{2}$ possibilities to choose such a pair of vertices $v_{i1}$ and $v_{i1'}$ symmetric to the vertical axis. The $p$ vertices $v_0, v_p,v_{2p}...,v_{(p-1)p}$ cannot be chosen,because they would be used twice after rotation by a multiple of $\dfrac{360}{p}$. Therefore we have $ \vert O_3(p^{2})\vert=\dfrac{p\cdot(p-1)}{2}$ possible chains $\overline{v_{i1}, v_0, v_{i1'}}$ with middlevertex $v_0$.
Now we rotate this chain $O(3)$ $p$ times by the angle of $\dfrac{360}{p}$ and get $p$ chains $O(3)$ like the first one.
\begin{figure}
\begin{tabular}{| c | c |}
\hline\
sketch 3 & sketch 4 \\
\includegraphics[width=0.5\textwidth]{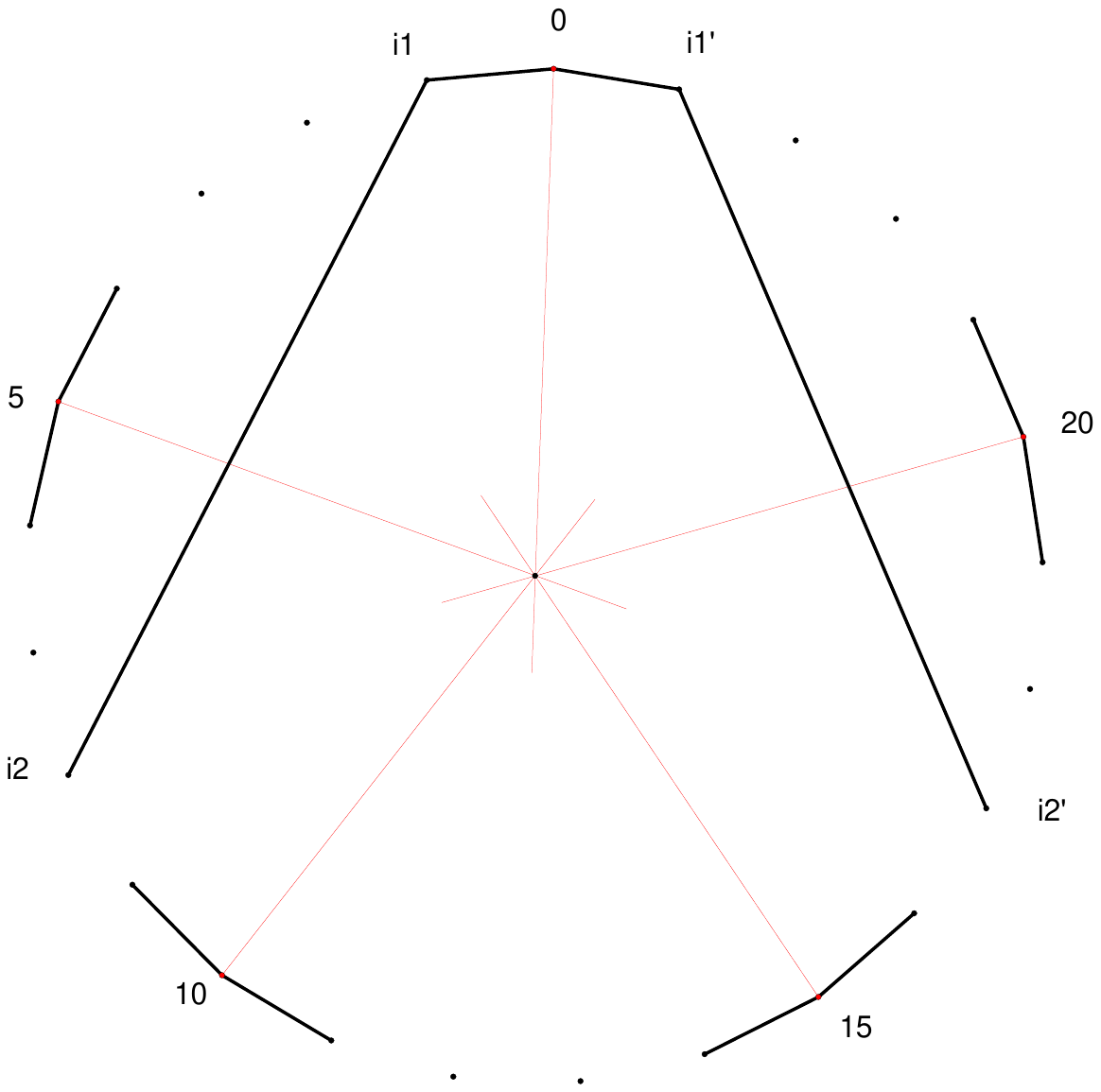} &
\includegraphics[width=0.5\textwidth]{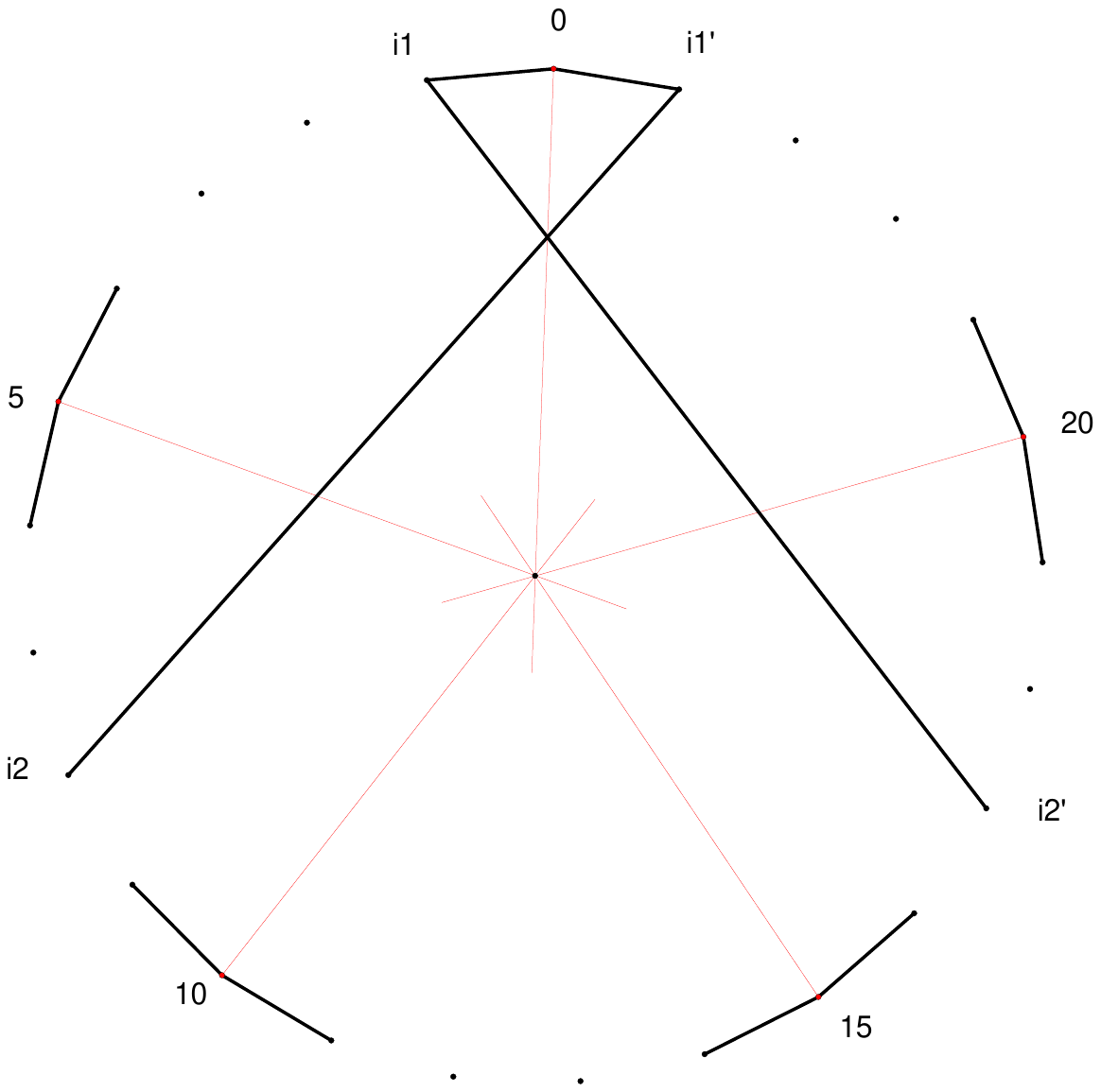} \\ \hline
\end{tabular}
\caption{construction-figures 3 and 4 for the equivalence-classes of the $p^{2}$-polygons with at least $p$ axes}
\end{figure}
\newpage
Sketch 3 and 4: We choose a 2nd pair of vertices ($v_{i2}$ and $v_{i2'}=v_{p^{2}-i2}$), which is symmetric to the vertical axis and not yet involved as a vertex of the former chains: For the choice of this second pair of vertices ($(v_{i2}$ and $v_{i2'}$) we have $\dfrac{p^{2}-3p}{2}= \dfrac{p \cdot(p-3)}{2}$ possibilities.
Now we connect $v_{i2}$  and $v_{i2'}$ with the ends $v_{i1}$ and $v_{i1'}$ to get two chains $O_1(5)$ and $O_2(5)$, which are symmetric to the vertical axis.There are $2$ possibilities of connection: $O_1(5)=\overline{v_{i2}v_{i1}v_0v_{i1'}v_{i2'}}$ and $O_2(5)=\overline{v_{i2}v_{i1'}v_0v_{i1}v_{i2'}}$.
\begin{figure}[H]
\begin{tabular}{| c | c |}
\hline\
sketch 5 & sketch 6 \\
\includegraphics[width=0.5\textwidth]{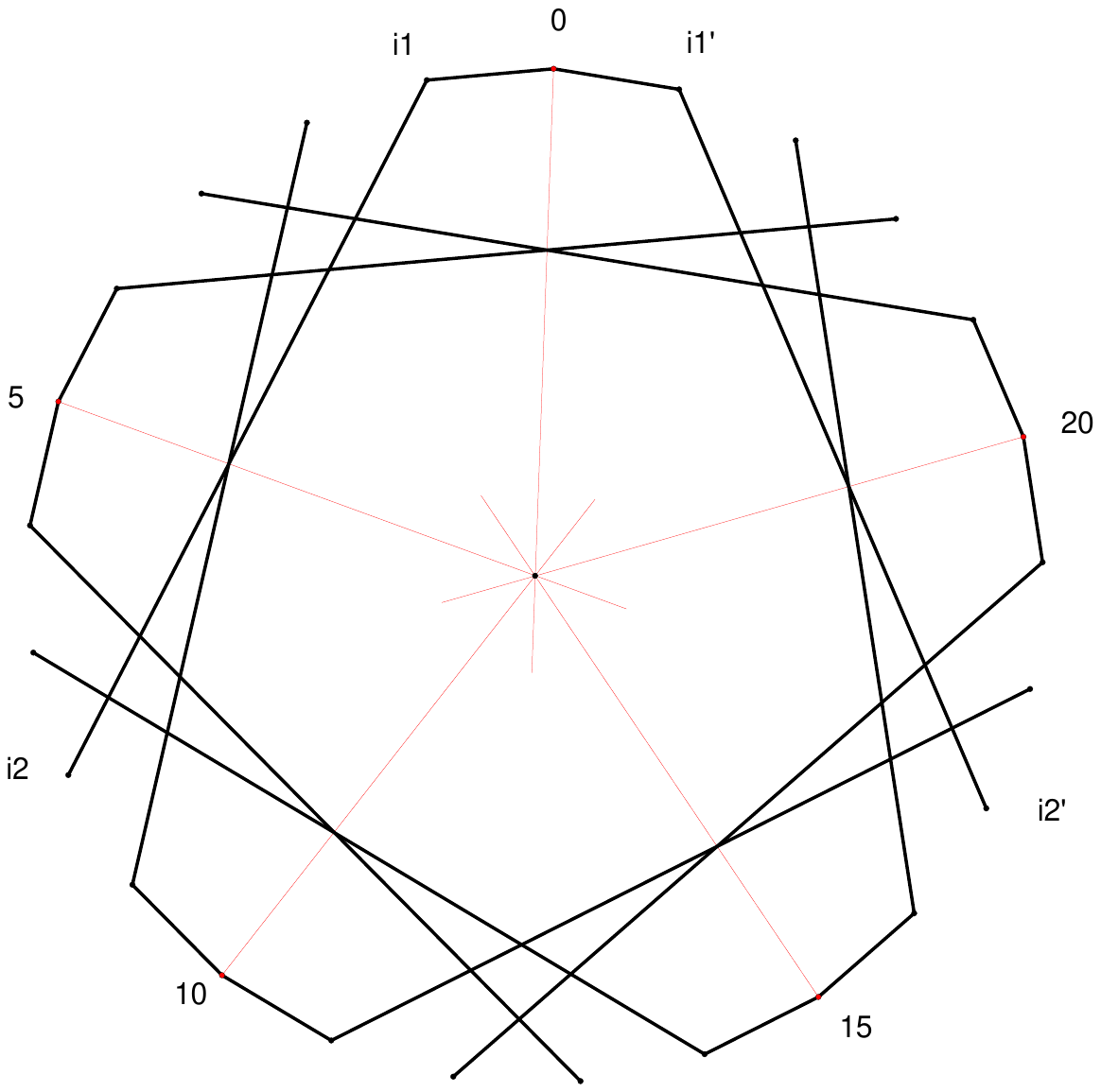} &
\includegraphics[width=0.5\textwidth]{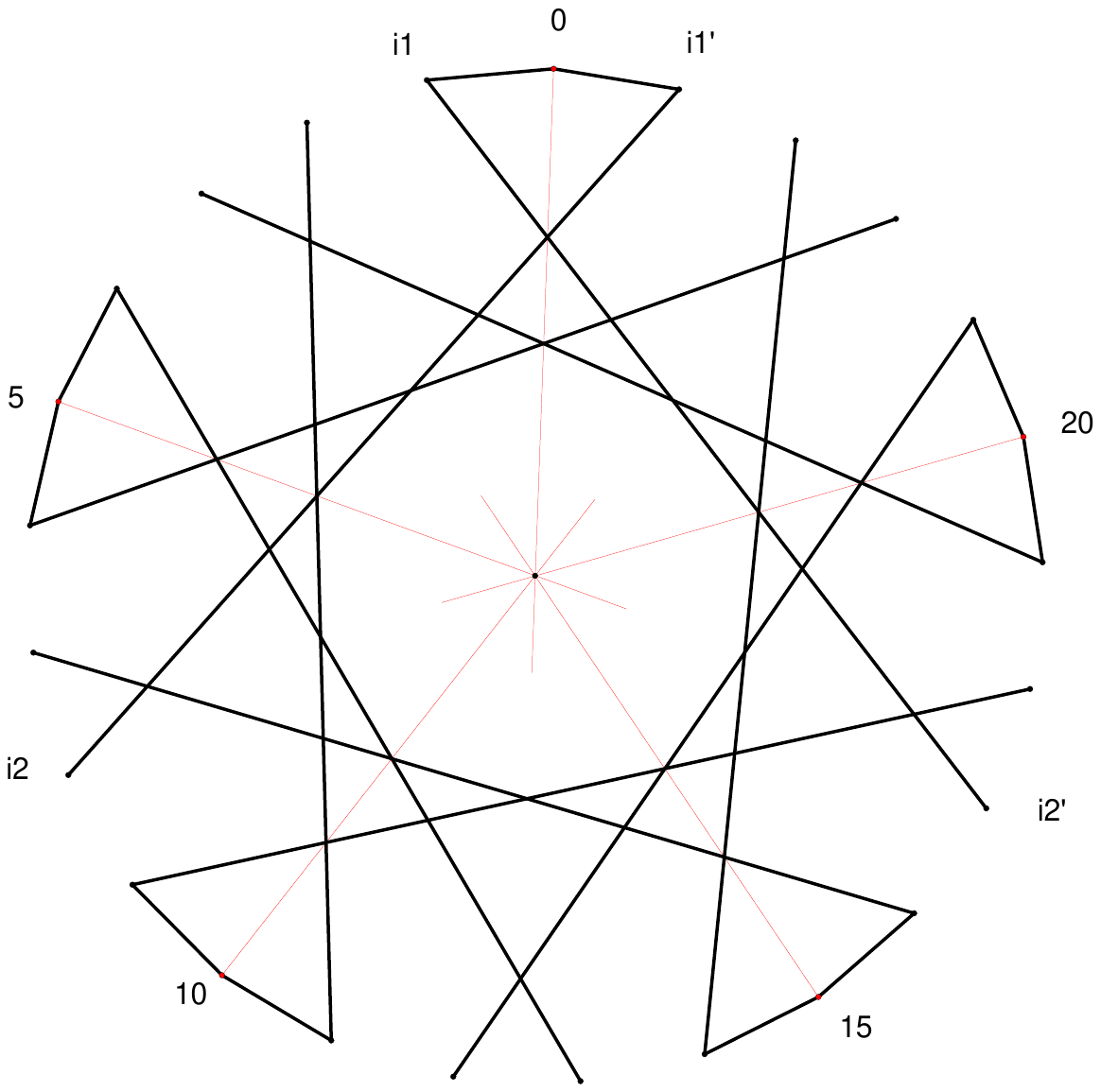} \\ \hline
\end{tabular}
\caption{construction-figures 5 and 6 for the equivalence-classes of the $p^{2}$-polygons with at least $p$ axes}
\end{figure}
Sketch 5 and 6: In both cases we rotate the chains $O_1(5)$ and $O_2(5)$ $(p-1)$- times by the angle of $\dfrac{360}{p}$.

Now we repeat this procedure until there exist no more vertices, which are not involved in the construction. So we constructed in this way $p$ symmetric open chains each with $p$ vertices. The $p$ symmetry-axes go through the vertices $v_0, v_p, v_{2p} \cdots v_{(p-1)\cdot p}$.
\begin{figure} [H]
\begin{tabular}{| c | c |}
\hline\
sketch 7& sketch 8 \\
\includegraphics[width=0.5\textwidth]{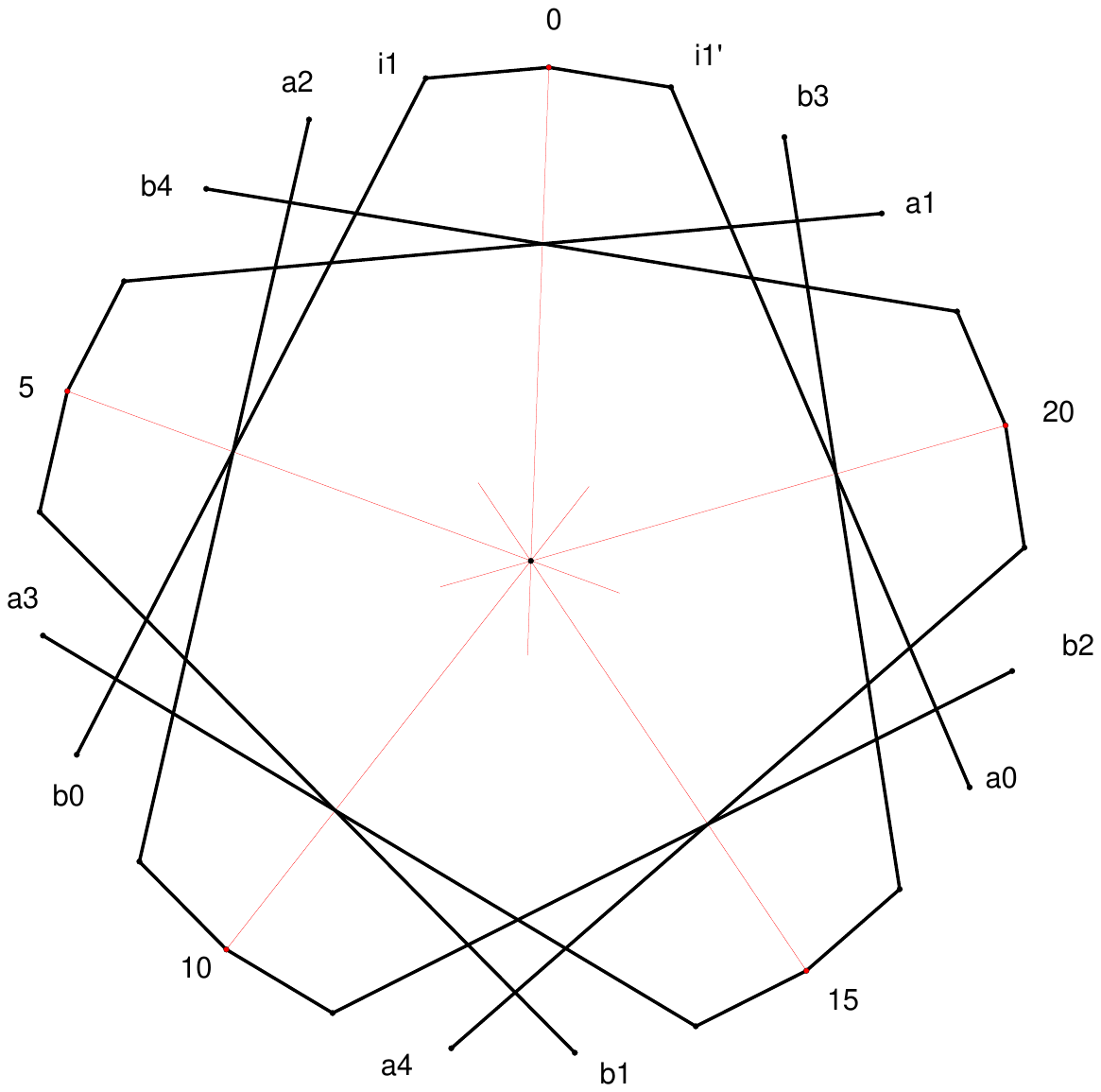} &
\includegraphics[width=0.5\textwidth]{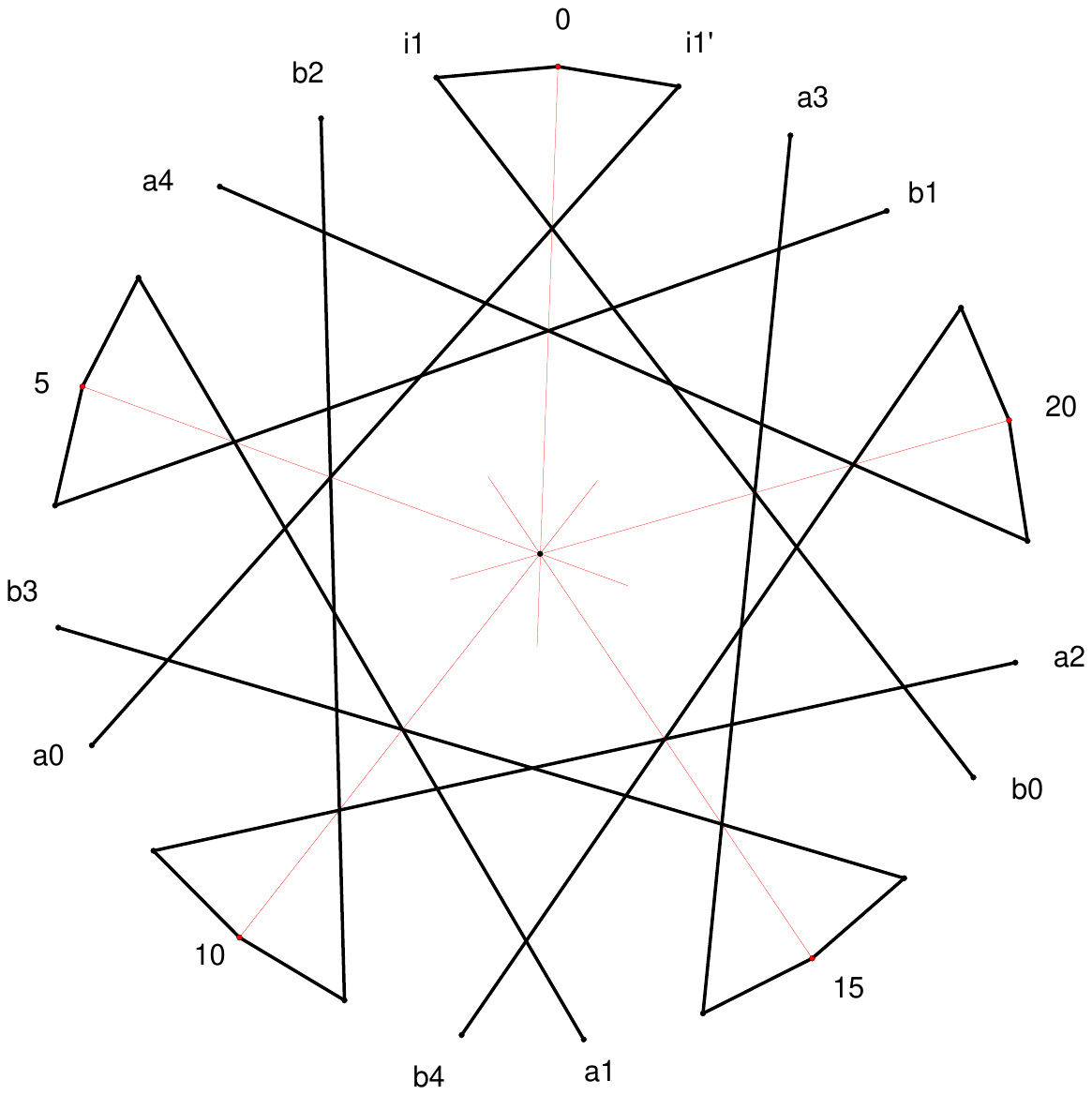} \\ \hline
\end{tabular}
\caption{construction-figures 7 and 8 for the equivalence-classes of the $p^{2}$-polygons with at least $p$ axes}
\end{figure}

The final step is now to connect the open chains to a $p^{2}$-polygon with at least $p$ axes. There are exactly $p-1$ possibilities to do so. We look at an open chain $O(p)$, i.g. let $O(p)=\overline{v_a\cdots v_0 \cdots v_b}$. If we would connect now $v_b$ with $v_a$, we would get no $p^{2}$-polygon. But we can connect $v_b$ with any other vertex $v_x$,where $x \in a+p, a+2p, \cdots a+(p-1)\cdot p$. 
\begin{figure} [H]
\begin{tabular}{| c | c |}
\hline\
sketch 9& sketch 10 \\
\includegraphics[width=0.5\textwidth]{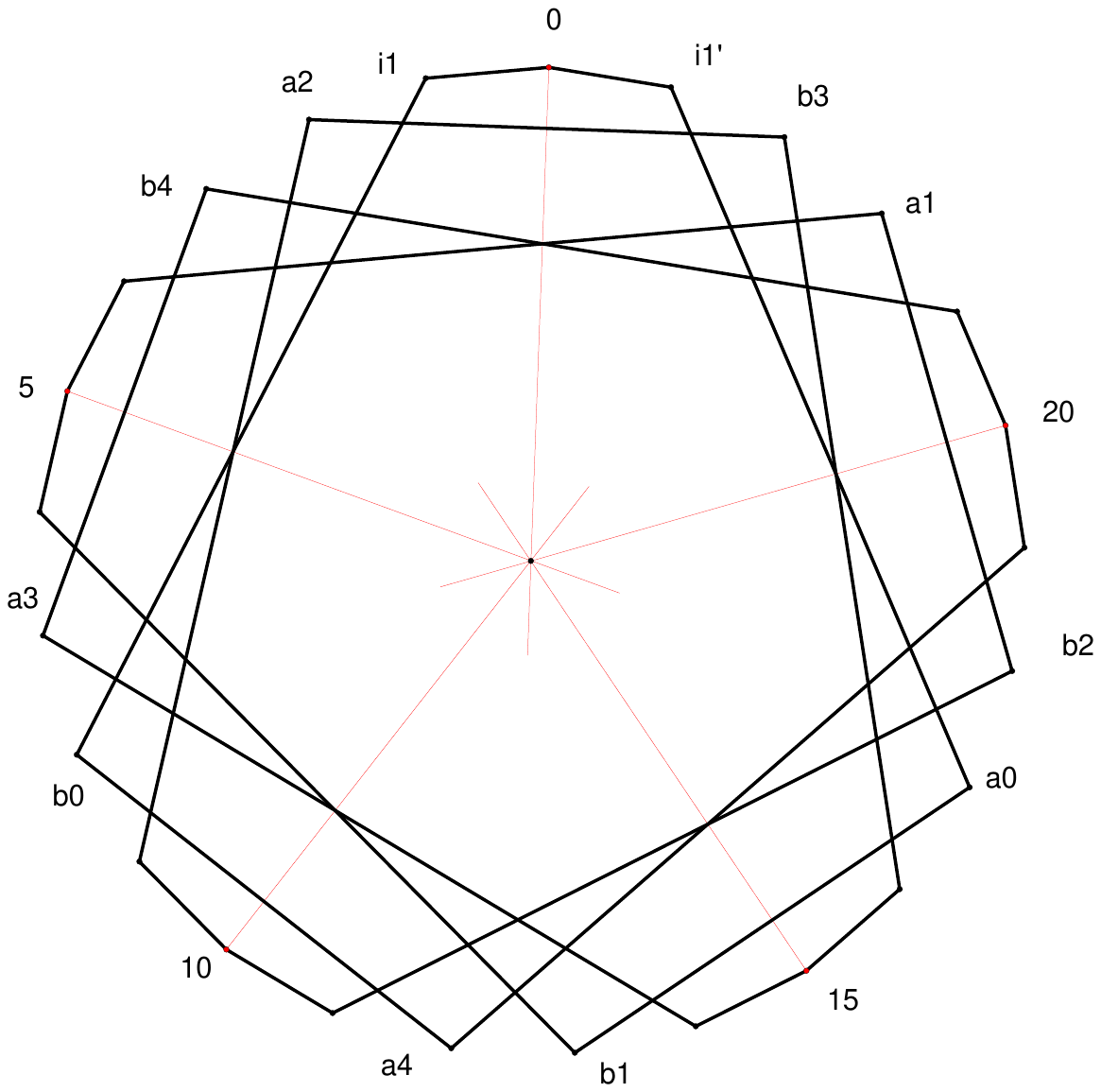} &
\includegraphics[width=0.5\textwidth]{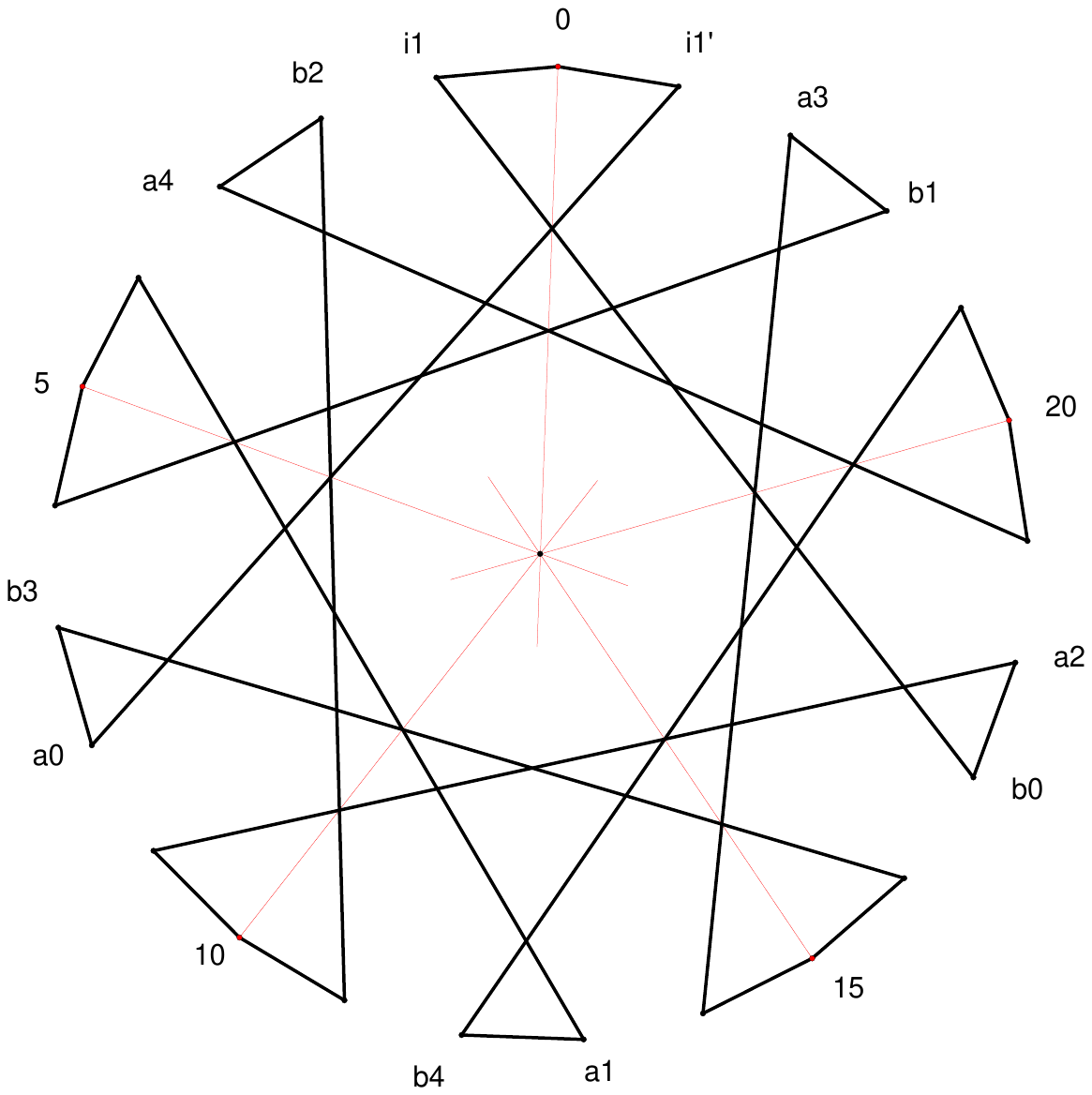} \\ \hline
\end{tabular}
\caption{construction-figures 9 and 10 for the equivalence-classes of the $p^{2}$-polygons with at least $p$ axes}
\end{figure}

Then we can rotate the edge $\overline{v_xv_b}$ $p$ times by the angle $\dfrac{360}{p}$ to get a representatif of an equivalence class of the $p^{2}$-polygons with at least $p$ axes by connecting the $p$ open chains in this way.\\\\
\textbf{Counting of the possibilities and elementary transformations}

$\vert X_{p+}(p^{2}) \vert=(p-1) \cdot \left[ \dfrac{p\cdot(p-1)}{2} \cdot 2 \cdot \frac{p\cdot (p-3)}{2}\cdots 2 \cdot \dfrac{p\cdot 4}{2} \cdot 2 \cdot \dfrac{p \cdot2}{2}\right]$\\\
By extracting the first denominator $2$ out of the brackets and shorten the other $2$-denominators:
$\vert X_{p+}(p^{2}) \vert=\dfrac{p-1}{2} \cdot \left[2p\cdot 4p \cdot \cdots \cdot (p-3)p \cdot (p-1)p \right]$ \\\\
By combining: $\vert X_{p+}(p^{2}) \vert=
(\dfrac{p-1}{2})\cdot (2p)^{\dfrac{p-1}{2}}\cdot\left(\dfrac{p-1}{2} \right)!$

After subracting $X_{p^{2}}((p^{2}) \vert$ we get the final result:
\begin{center}
\fbox{$\vert X_p(p^{2}) \vert =(\dfrac{p-1}{2})\cdot \left[(2p)^{\dfrac{p-1}{2}}\cdot\left(\dfrac{p-1}{2} \right)!-p\right]$}
\end{center}
\end{proof}
\newpage
\subsection{Main sentence d.) $1$ axis}
\label{subsec:main_sentence_d.)_1_axis} 
\begin{proof}
To develop this proof we need two construction figure with $p^{2}= 49$ vertices. The one axis of the $7$ axes is set vertically and red.
\begin{figure}[H]
\includegraphics[width=\textwidth]{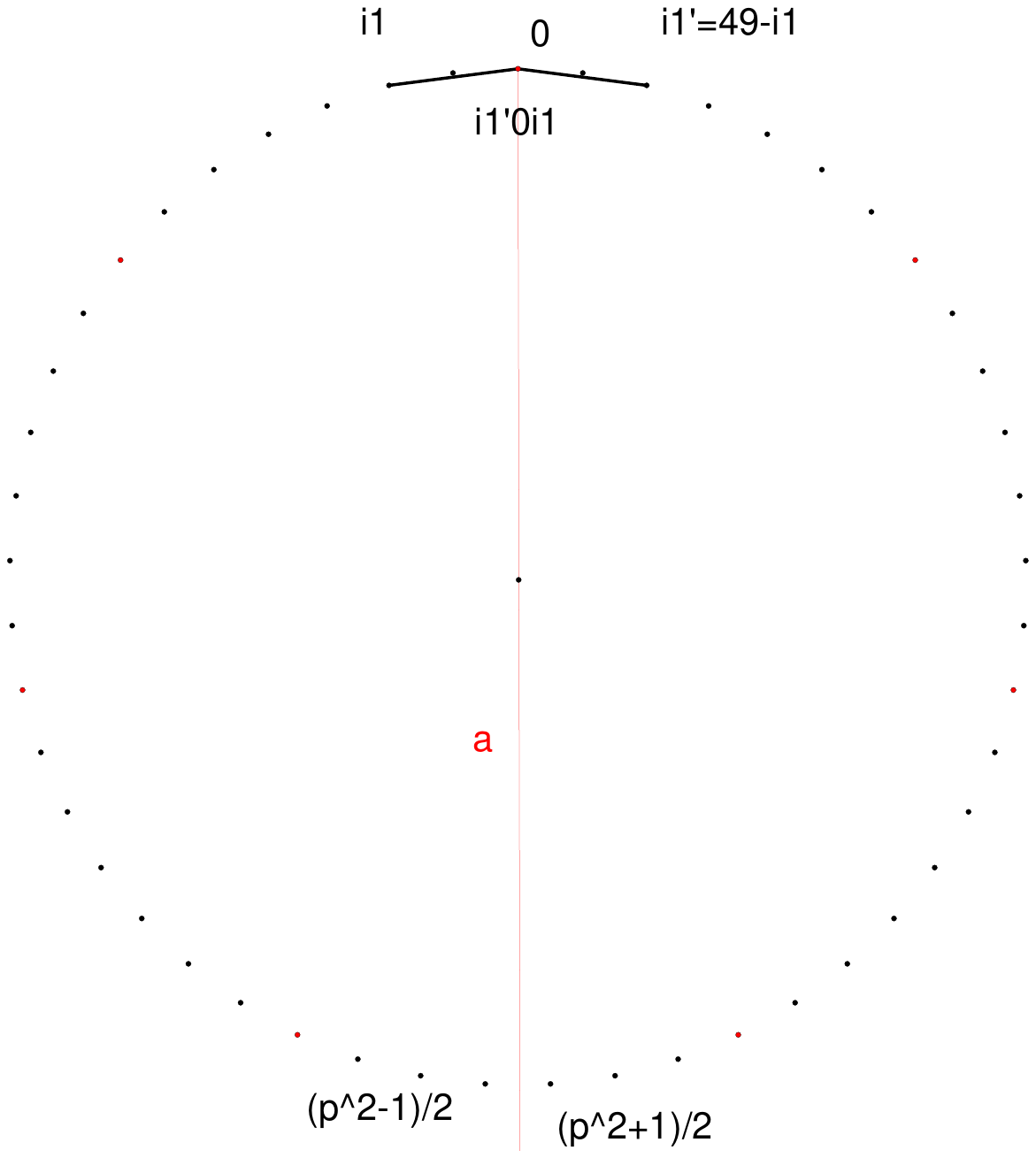} 
\caption{construction-figure 1 for the equivalence-classes of the $p^{2}$-polygons with at least $1$ axis}
\end{figure}
There are $\dfrac{p^{2}-1}{2}$ possibilities to choose the first pair of vertices $v_{i1}$ and $v_{i1'}$, which are situated symetrically to the axis a.We connect them and get an open chain $C(3)=\overline{ v_{i1'}v_0v_{i1}}$.
\begin{figure}[H]
\includegraphics[width=\textwidth]{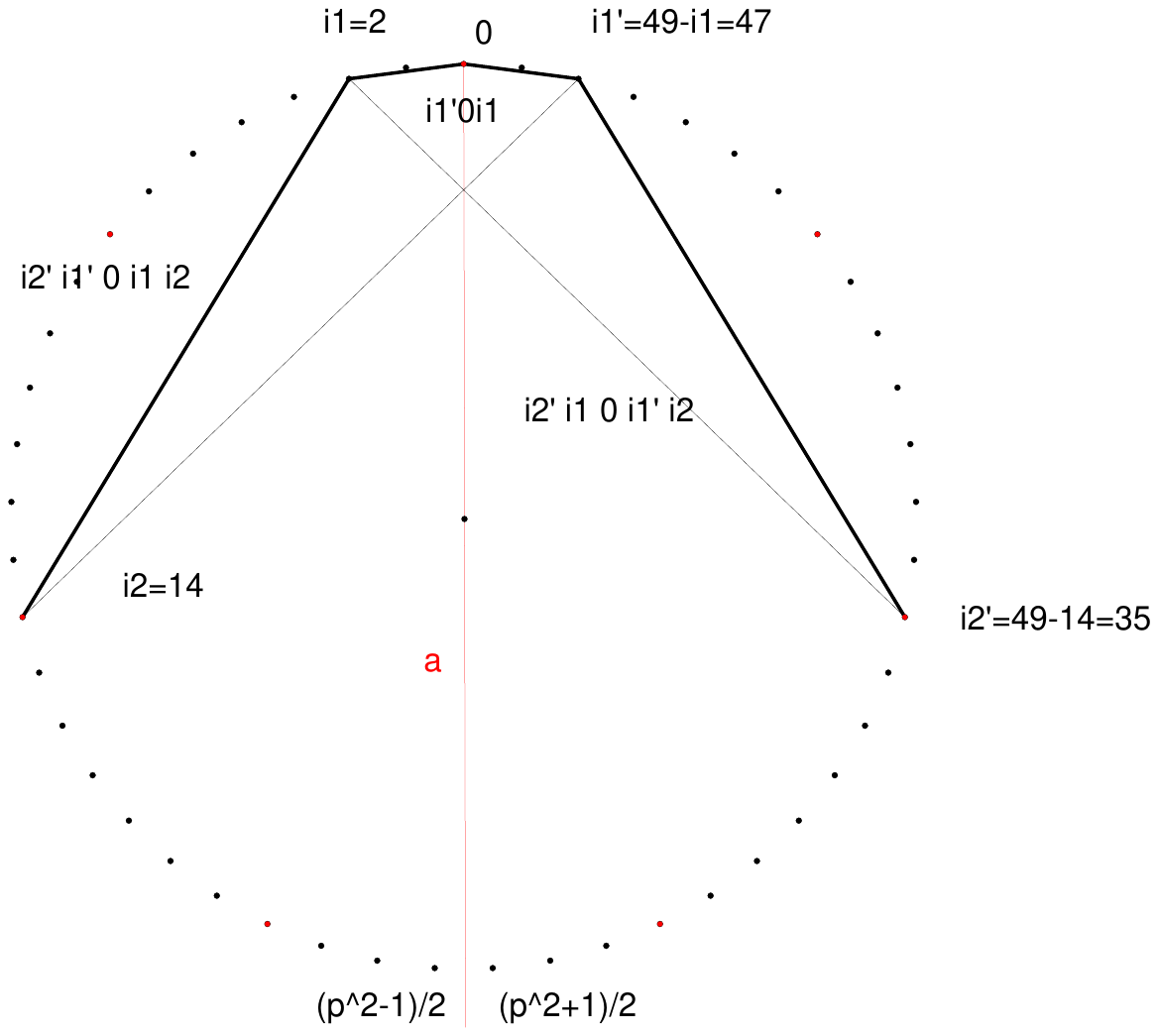} 
\caption{construction-figure 2 for the equivalence-classes of the $p^{2}$-polygons with at least $1$ axis}
\end{figure}

There are $\dfrac{p^{2}-1}{2}-1= \dfrac{p^{2}-3}{2}$ possibilities to choose the second pair of vertices $v_{i2}$ and $v_{i2'}$, which are also situated symetrically to the axis a.
There are now $2$ possibilities to connect $v_{i2}$ and $v_{i2'}$ with the chain $C_3$. We get two open chains $C_1(5)$ and $C_2(5)$: $C_1(5)=\overline{v_{i2'}v_{i1'}v_0v_{i1}v_{i2}}$, drawn thick, and $C_2(5)=\overline{v_{i2'}v_{i1}v_0v_{i1'}v_{i2}}$, drawn thin.
Now we iterate this procedure until we get open chains with all $p^{2}$ vertices.The last pair of symmetrically situated vertices  will be chosen with  only $1$ choice and has also $2$ possibilities of connection with each of the yet constructed $C(p^{2}-2)$ chains. Finally we link the last chosen $2$ vertices by a horizontal edge.
By this construction we get a representatif of every equivalence-class of $p^{2}$-polygons with at least $1$ symmetry-axis:\\\\
$\vert X_{1+}(p^{2} \vert  = \dfrac{p^{2}-1}{2} \cdot 2 \cdot\dfrac{p^{2}-3}{2} \cdot \cdots \cdot 2 \cdot \dfrac{4}{2}\cdot 2 \cdot\dfrac{2}{2}$

$\vert X_{1+}(p^{2} \vert  =2 \cdot 1 \cdot 2 \cdot 2\cdot \cdots \cdot 2 \cdot \dfrac{p^{2}-3}{2}\cdot \dfrac{p^{2}-1}{2}$

$\vert X_{1+}(p^{2} \vert  =\dfrac{1}{2} \cdot 2 \cdot 1 \cdot 2 \cdot 2\cdot \cdots \cdot 2 \cdot \dfrac{p^{2}-3}{2}\cdot 2 \cdot \dfrac{p^{2}-1}{2}$

$\vert X_{1+}(p^{2} \vert  = \dfrac{1}{2}\cdot\left( \dfrac{p^{2}-1}{2}\right)!\cdot 2^{\dfrac{p^{2}-1}{2}}$

After subtracting $\vert X_{p+}(p^{2} \vert $ we get the final result:
\begin{center}
\fbox{$\vert X_1(p^{2} \vert  =\dfrac{1}{2}\cdot\left( \dfrac{p^{2}-1}{2}\right)!\cdot 2^{\dfrac{p^{2}-1}{2}}-(\dfrac{p-1}{2})\cdot (2p)^{\dfrac{p-1}{2}}\cdot\left(\dfrac{p-1}{2} \right)!$}
\end{center}
\end{proof}
\subsection{main-sentence e.) $p$-circular $p^{2}$-polygons}
\label{subsec:main-sentence_e.)_p-circular_p2-polygons}
\begin{proof}
As a preparation we look at the special case $p^{2}=9$:

\begin{figure} [H]
\begin{tabular}{| c | c | c |}
\hline\
sketch 1 & sketch 2 & sketch 3 \\
\includegraphics[width=0.3\textwidth]{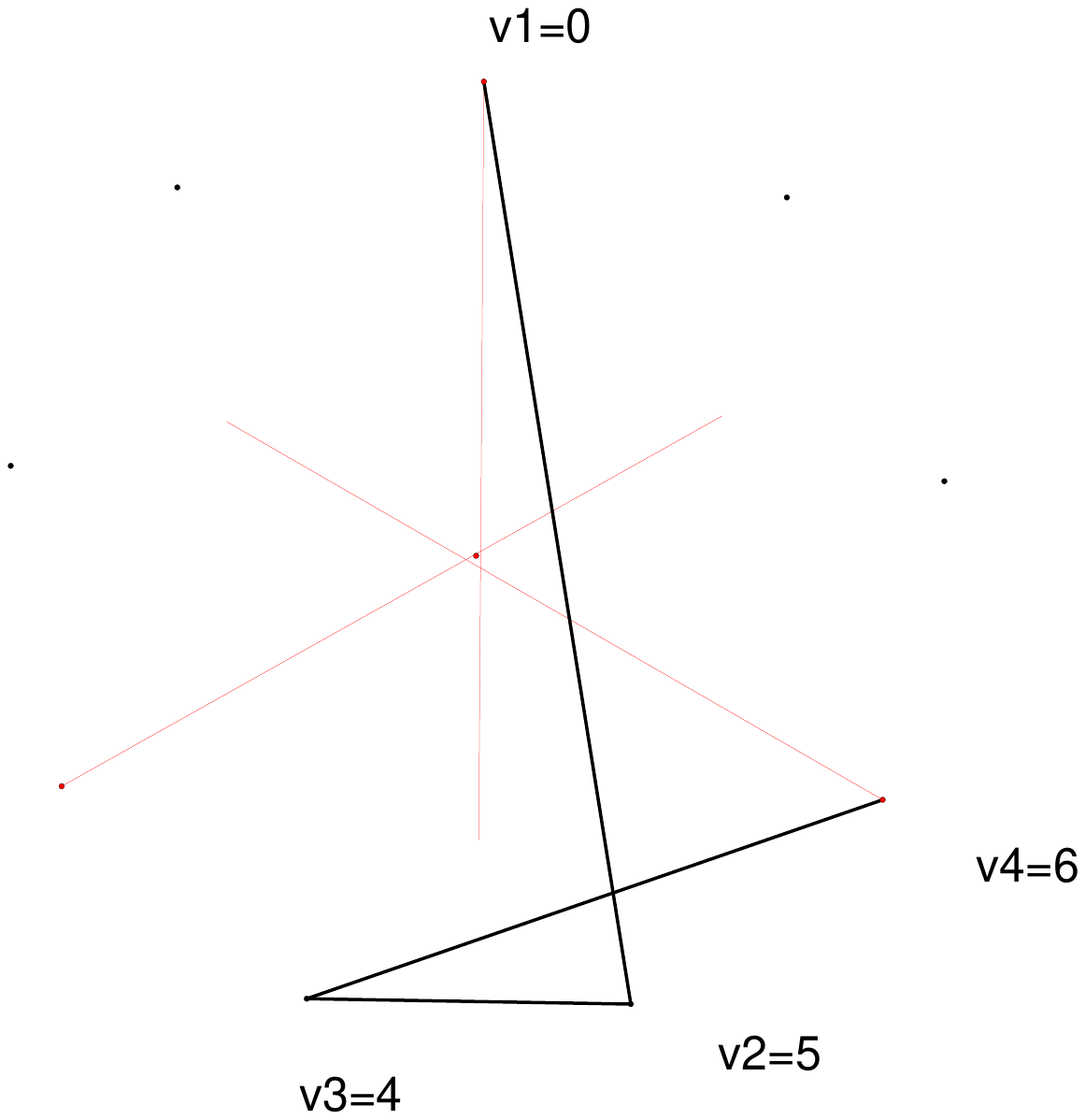} &
\includegraphics[width=0.3\textwidth]{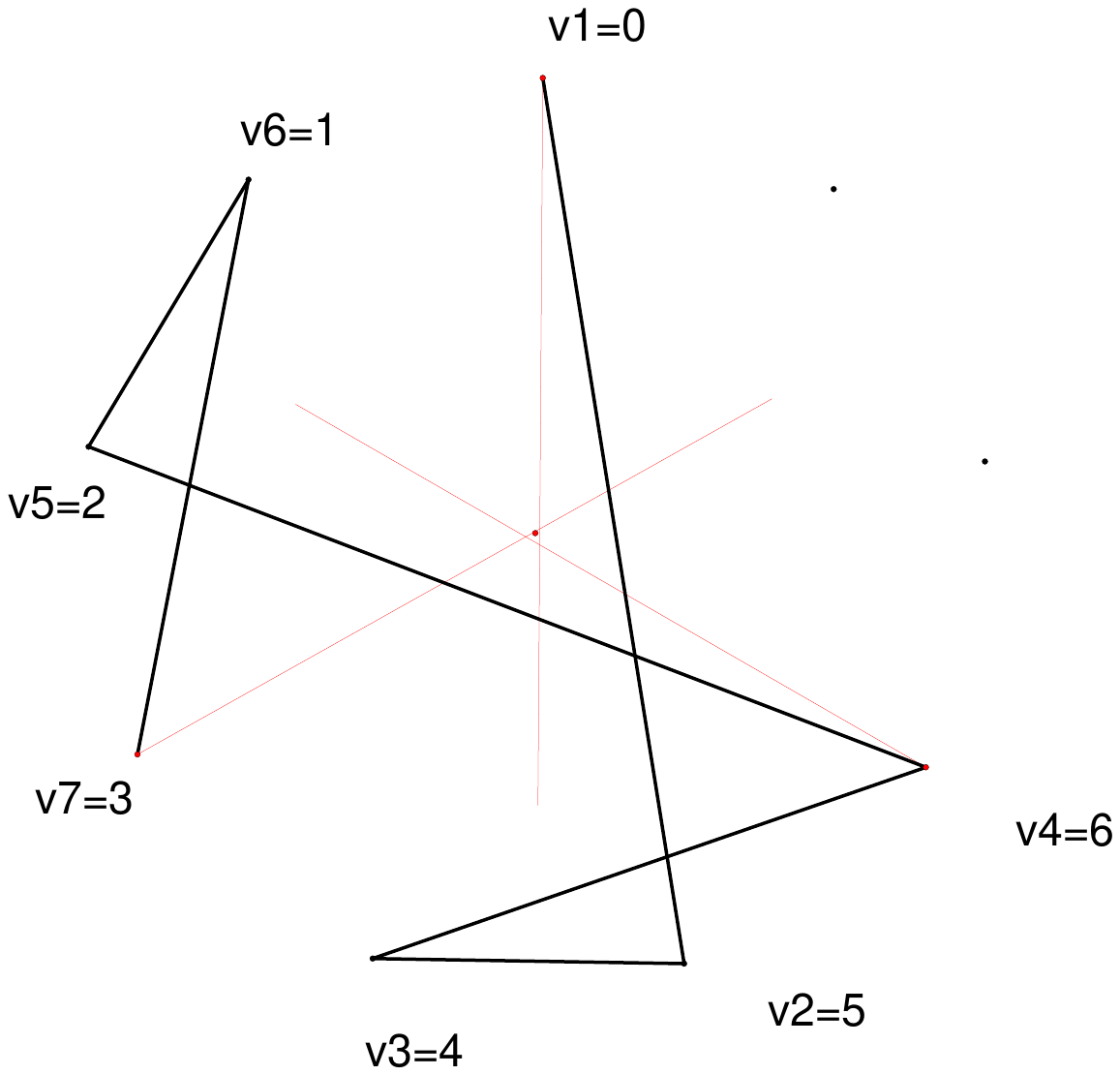} &
\includegraphics[width=0.3\textwidth]{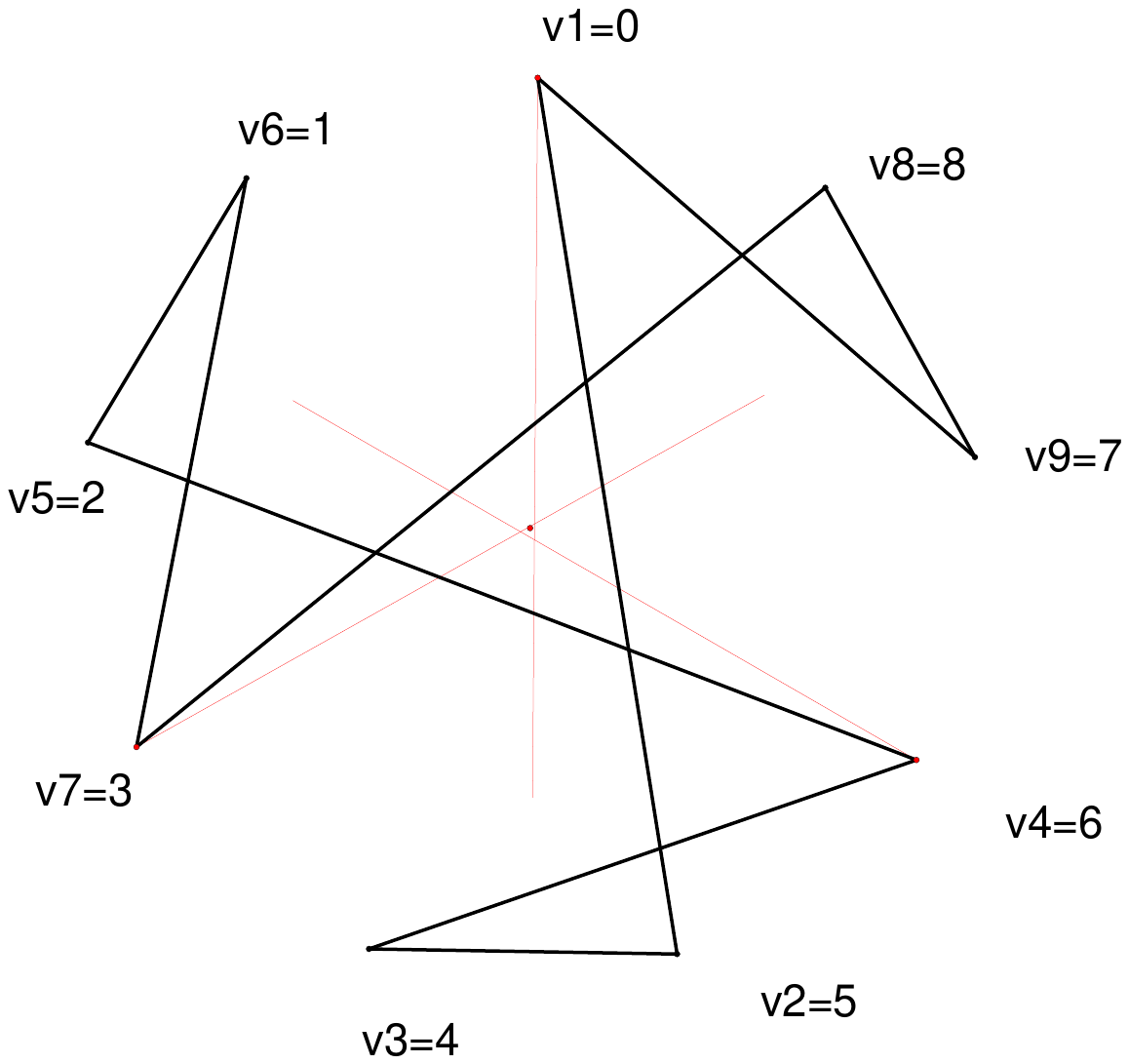} \\ \hline
\end{tabular}
\caption{construction-figures 1-3 for the equivalence-classes of the $p$-circular $p^{2}$-polygons}
\end{figure}
\begin{enumerate}
\item{$v_1=0$ is the starting vertex}
\item{We have $p^{2}-p=(p-1) \cdot p$ possibilities to choose $v_2$ from the points labeled with $v_2=a_2$ with $a_2 \not\equiv 0\pmod{p}$. We choose e.g. $v_2=5$.}
\item{We have $p^{2}-2p=(p-2)\cdot p$ possibilities to choose $v_3=a3$ with $a_3\not\equiv 0\pmod{p}$ and $a_3 \not\equiv a_2 \pmod{p}$. We choose e.g. $v_3=4$.}
\item{We have finally $p-1$ possibilities to choose $v_4$ from the points labeled with ${3,6}$, e.g. $v_4=6$.}
\item{Connecting the chosen vertices by edges we get open chains of $p+1$ vertices and $p$ edges, e.g. $O_4=\overline{v_1 v_2 v_3 v_4}= \overline{0 5 4 6}$}
\end{enumerate}
We constructed so far $p \cdot(p-1) \cdot p \cdot( p-2 )\cdot (p-1) =36 $ such open chains $O_4$.
Now we rotate (several times) the chain $O_4$ by the angle$\dfrac{v_4}{p}\cdot\dfrac{360}{p}=2 \cdot 120=240$ till vertex $v_4$ is coincident with vertex $v_1$ to get a closed chain, a $p^{2}$-polygon: $\overline{0 5 4 6 2 1 3 8 7 0}$, which either has $p^{2}$ axes or $p$ axes or no axes,but $p$-circularity.\\\\
Now we want to examine, how many equivalence-classes exist among the so constructed $p^{2}$-polygons, which have no symmetry-axe, but are circular.
 \begin{enumerate}
 \item{We have always two identic $p^{2}$-polygons: One counterclockwise and one clockwise,e.g. $\overline{0 5 4 6 2 1 3 8 7 0}$ and $\overline{0 7 8 3 1 2 6 4 5 0}$} Therfore the number of equivalence-classes of $9$-polygons with circularity, with no axis or with $3$ axes or $9$ axes is $ \vert X_{c+}(9) \vert =\dfrac{p-1}{2}\cdot p \cdot (p-1) \cdot p\cdot (p-2)=\dfrac{p-1}{2} \cdot p^{2}\cdot (p-1) \cdot (p-2)=1\cdot 9 \cdot 2 \cdot 1=1 \cdot 9 \cdot 2 \cdot 1=18$.
\item{Among the constructed $\dfrac{p-1}{2}\cdot p^{2}\cdot (p-1)\cdot (p-2)=18$ equivalence-classes we constructed also those with $9$ axes, i.e. $\overline{0 2 4 6 8 1 3 5 7 0}$.}
\item{As there exist $\vert X_9(9)\vert=\dfrac{p\cdot (p-1)}{2}=3$ equivalence classes of $9$-polygons with $9$ axes there remain $15$ equivalence-classes of $p^{2}$-polygons beyond the constructed ones, which have $3$ or no symmetry-axes.}
 \item{Finally we have to exclude all sequences of vertices, which lead to $p^{2}$-polygons with $p$ axes. Because every such equivalence-class has $p$ representatifs, we have to subtract from $15$ $3$-times the number of equivalence-classes of $9$-polygons with $3$ axes, i.e. $15-p \cdot \vert X_3(9)\vert=15-3\cdot 3=6$. This is correct, because e.g.$\overline{0 2 7 3 5 1 6 8 4 0}$, $\overline{0 5 1 3 8 4 6 2 7 0}$ and $\overline{0 5 7 3 8 1 6 2 4 0}$ represent the same equivalence-class of $9$-polygons with $3$ axis.}
\item{The $6$ left sequences of vertices belong to the $3$-circular $9$-polygons with no symmetry-axes, and there are $3$ representatifs of each equivalence-class among the remaining $6$ sequences of vertices.}
\end {enumerate}
The proven results for  $p=3$ is:
\begin{center}
$\vert X_c(9)\vert=\dfrac{6}{3}=2$ and
$\vert X_{c+}(9\vert)=18$
\end{center}

\textbf{Proof in general for $ \vert X_{c+}(p^{2}) \vert$ and $\vert X_c(p^{2})\vert$:}
We construct and enumerate first all chains $O(p)$ with $p+1$ edges an $p$ vertices starting in vertex $v_1=0$ and ending in a vertex $v_p= k \cdot p$ with $k \in {1,2, \cdots p-1}$. For this purpose we define ${0,1,\cdots p^{2}-1,p^{2}}$ as the labels of the vertices counterclockwise.

\begin{enumerate}
\item{We choose $v_2=a_2$ with $a_2 \not\equiv 0 \pmod{p}$. We have $p^{2}-p=p \cdot (p-1)$ possible choices.}
\item{We choose $v_3=a_3$ with $a_3 \not\equiv a_2 \pmod{p}$ and $a_3 \not\equiv 0 \pmod{p}$. We have $p^{2}-2p=p \cdot (p-2)$ possible choices.}
\item{We choose $v_4=a_4$ with $a_4 \not\equiv 0 \pmod{p}$ and $a_4 \not\equiv a_2 \pmod{p}$ and $a_4 \not\equiv a_3 \pmod{p}$. We have $p^{2}-3p=p \cdot (p-3)$ possible choices, if $p>3$.}
\item{We continue in this way until we have chosen all sequences of $p-1$ vertices: ${0,a_2,a_3, \cdots a_{p-1}}$}
\item{As vertex $v_p$ we choose $v_p=a_p$ with $a_p\neq 0$ but $a_p \equiv 0 \pmod{p}$. For this step we have $p-1$ choices.}
\item{In total we get $(p-1) \cdot p\cdot(p-1\cdot)p\cdot(p-2) \cdots p\cdot 2 \cdot p \cdot 1$ open chains $O(p)$ starting in vertex $v_1=0$ and ending in $v_p=a_p$ or shortly: $\vert O(p) \vert=(p-1)\cdot p^{p-1}\cdot (p-1)!$.}
\end{enumerate}
Now we rotate the chain $O(p)$ by the angle$\dfrac{v_p}{p}\cdot\dfrac{360}{p}$ till vertex $v_p$ is coincident with vertex $v_1=0$ to get a closed chain, a $p^{2}$-polygon, which has the property of $p$-circularity, but also, perhaps, $p$ or $p^{2}$ axes.\\\\
We see, that  there exist $(p-1)\cdot p^{(p-1)}\cdot ((p-1)!)$ $p^{2}$-polygons, which are constructed in that manner.
\begin{enumerate}
\item{But we have always two identic $p^{2}$-polygons: One counterclockwise and one clockwise.So there are only $\dfrac{p-1}{2}\cdot p^{p-1}\cdot (p-1)!:=\vert C_p(p^{2})\vert$ so constructed $p^{2}$-polygons.}
\item{Among them are also the regular polygons. Because $\vert X_{p^{2}}(p^{2})\vert=\dfrac{p\cdot (p-1)}{2}$ there remain $\dfrac{p-1}{2}\cdot\left[ p^{(p-1)}\cdot ((p-1)!)-p\right]$ so constructed $p^{2}$-polygons, which have less than $p^{2}$ or no symmetry-axes,but have the property of circularity.}
\item{Next we exclude the $p^{2}$-polygons with $p$ axes. Because every equivalence-class of $p^{2}$-polygons with $p$ axes has $p$ representatifs, after the exclusion of all $p^{2}$-polygons with $p$ axes, there remain $\dfrac{p-1}{2}\cdot \left[p^{(p-1)}\cdot ((p-1)!)-p\right]-p\cdot \vert X_p(p^{2})\vert$ polygons with no axes,but the property of $p$-circularity.}
\item{Because every equivalence-class of the circular $p^{2}$-polygons has $p$ representatifs, we divide by $p$ to get the number $\vert X_c(p^{2})\vert$:}
\end{enumerate}
$\vert X_c(p^{2})\vert=\dfrac{p-1}{2p}\cdot \left[p^{p-1}\cdot((p-1)!)-p-p\cdot(\dfrac{p-1}{2})!\cdot \left(2p\right)^{\dfrac{p-1}{2}}+p^{2}\right]$\\\\

$\vert X_c(p^{2}) \vert =\dfrac{p-1}{2}\cdot \left[((p-1)!)\cdot p^{p-2}-1-\left(\dfrac{p-1}{2}\right)!\cdot\left(2p\right)^{\dfrac{p-1}{2}}+p\right]$
\begin{center} 
\fbox{$\vert X_c(p^{2})\vert=\dfrac{p-1}{2}\cdot \left[((p-1)!)\cdot p^{p-2}-\left(\dfrac{p-1}{2}\right)!\cdot\left(2p\right)^{\dfrac{p-1}{2}}+p-1\right]$}
\end{center}
\end{proof}

\subsection{Main-sentence f.) Asymmetric and not circular $p^{2}$-polygons}
\label{subsec:main-sentence_f.)_asymmetric_and_not_circular_p2-polygons}
\begin{proof}
To calculate the number $\vert X_a(p^{2})\vert$ we have to subtract from the number of equivalence-classes of all $p^{2}$-polygons the numbers $X_{1+}(p^{2})$ and $X_c(p^{2})$:
\begin{center}
 $\vert X_a(p^{2})\vert=$
 
 $\dfrac{1}{2 \cdot p^{4}} \cdot \left[p^{4} \cdot (p-1)^{2}+(p-1)^{2} \cdot p! \cdot p^{p}+(p^{2})!    \right]-$
 
$\dfrac{1}{2}\cdot\left( \dfrac{p^{2}-1}{2}\right)!\cdot 2^{\dfrac{p^{2}-1}{2}}-$

$\dfrac{p-1}{2}\cdot \left[(p-1)!\cdot p^{p-2}-\left(\dfrac{p-1}{2}\right)!\cdot\left(2p\right)^{\dfrac{p-1}{2}}+p-1\right]$
\end{center}
\end{proof}
\newpage

\section{Table of the results of part III}
\label{sec:table_of_the_results_of_part_III}
\begin{center}
\begin{table}[!h]
\centering
\begin{tabular}{|| c | c | c | c | c | c | c ||}
\hline
\hline
 $p^{2}$ & $\vert X(p^{2}) \vert$ & $\vert X_{p^{2}}(p^{2})\vert$ & $\vert X_p(p^{2})\vert$ & $\vert X_1(p^{2})\vert$ &  $\vert X_c(p^{2}\vert$ & $\vert X_0(p^{2}) \vert$\\ \hline
9 & 2246 & 3 & 3 & 186 & 2 & 2052\\ \hline
25 & $1.24\cdot 10^{22}$ & 10 & 390 & 980995276400 & 5608 & $1.24 \cdot 10^{22}$\\ \hline
49 &  $1.27\cdot 10^{59}$ & 21& 49371 & $5.21 \cdot 10^{30}$ & 36253746 & $1.27 \cdot 10^{59}$\\ \hline
121 &  $2.76\cdot 10^{196}$ & 55 & 3092179145 & $4.80 \cdot 10^{99}$ & $4.28 \cdot 10^{16}$ & $2.76 \cdot 10^{196}$\\ \hline
\hline
\end{tabular}
\caption{Results for $9\leq p^{2} \leq 121$}
\end{table}
\end{center}

\section{Sets of representatifs of equivalence-classes}
\label{sec:sets_ofrepresentatifs_of_equivalence-classes}
\begin{figure}[H]
\begin{tabular}{|c | c | c |}
 \hline
\includegraphics[width=0.3\textwidth]{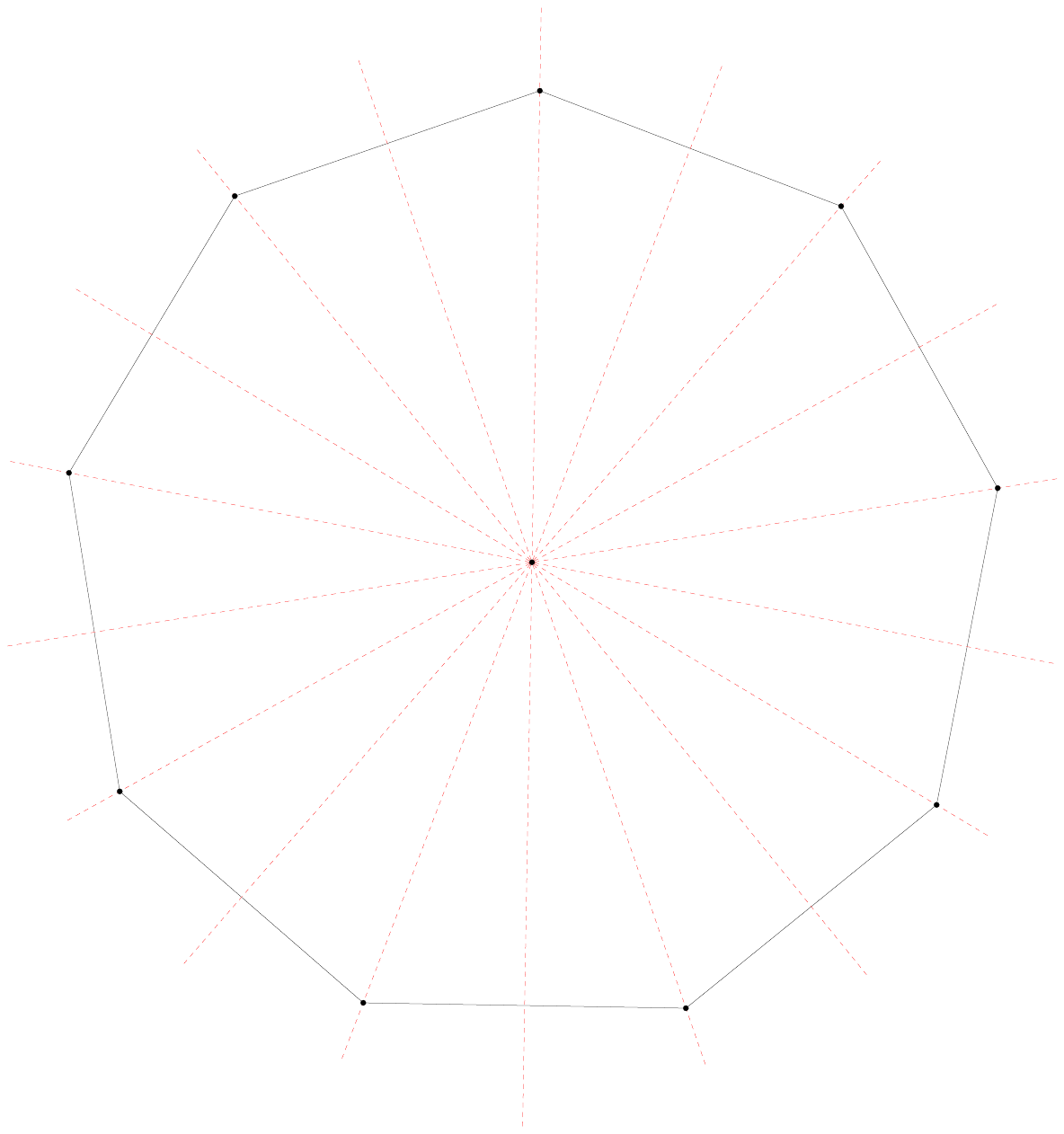} &
\includegraphics[width=0.3\textwidth]{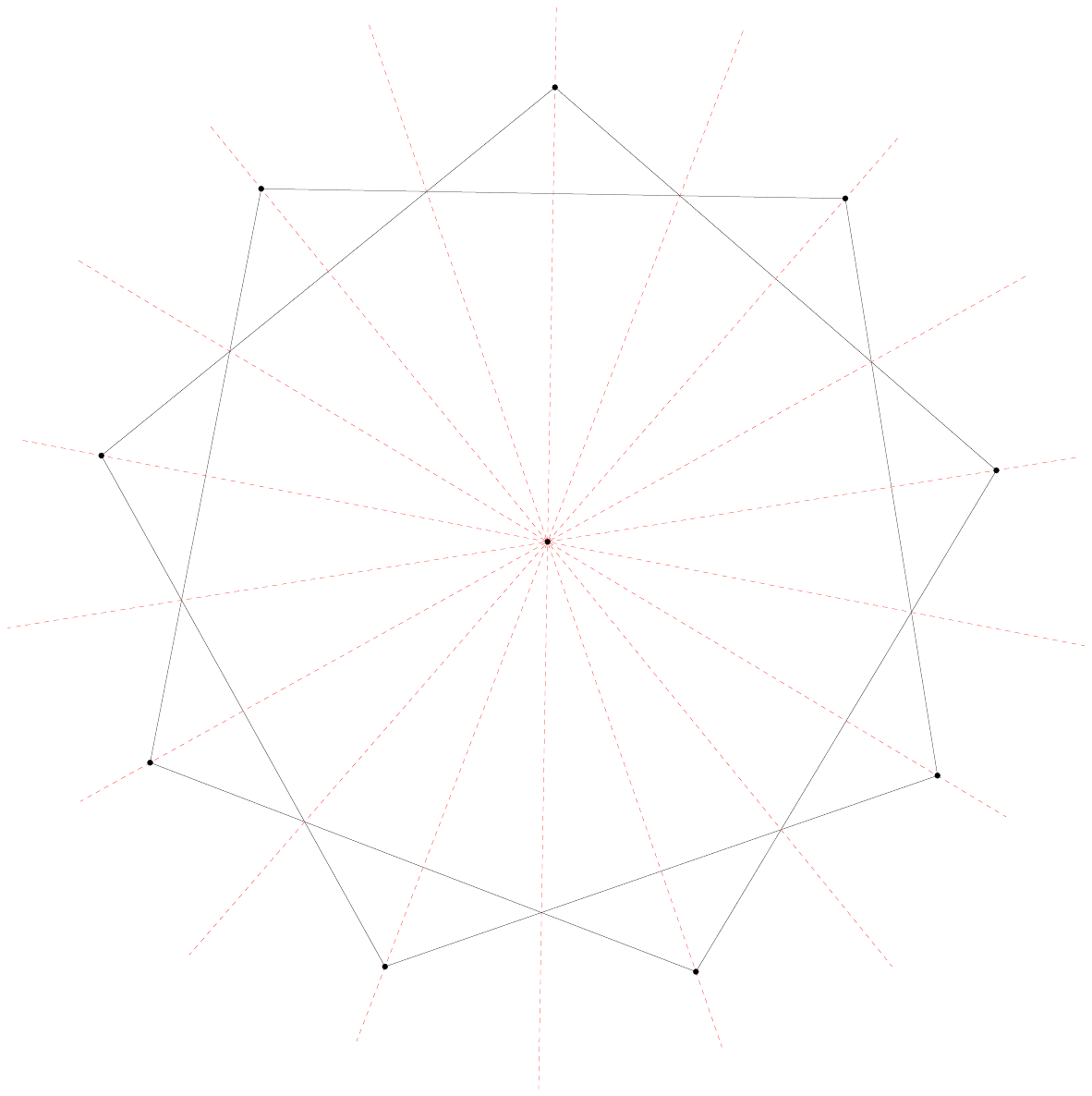} & \includegraphics[width=0.3\textwidth]{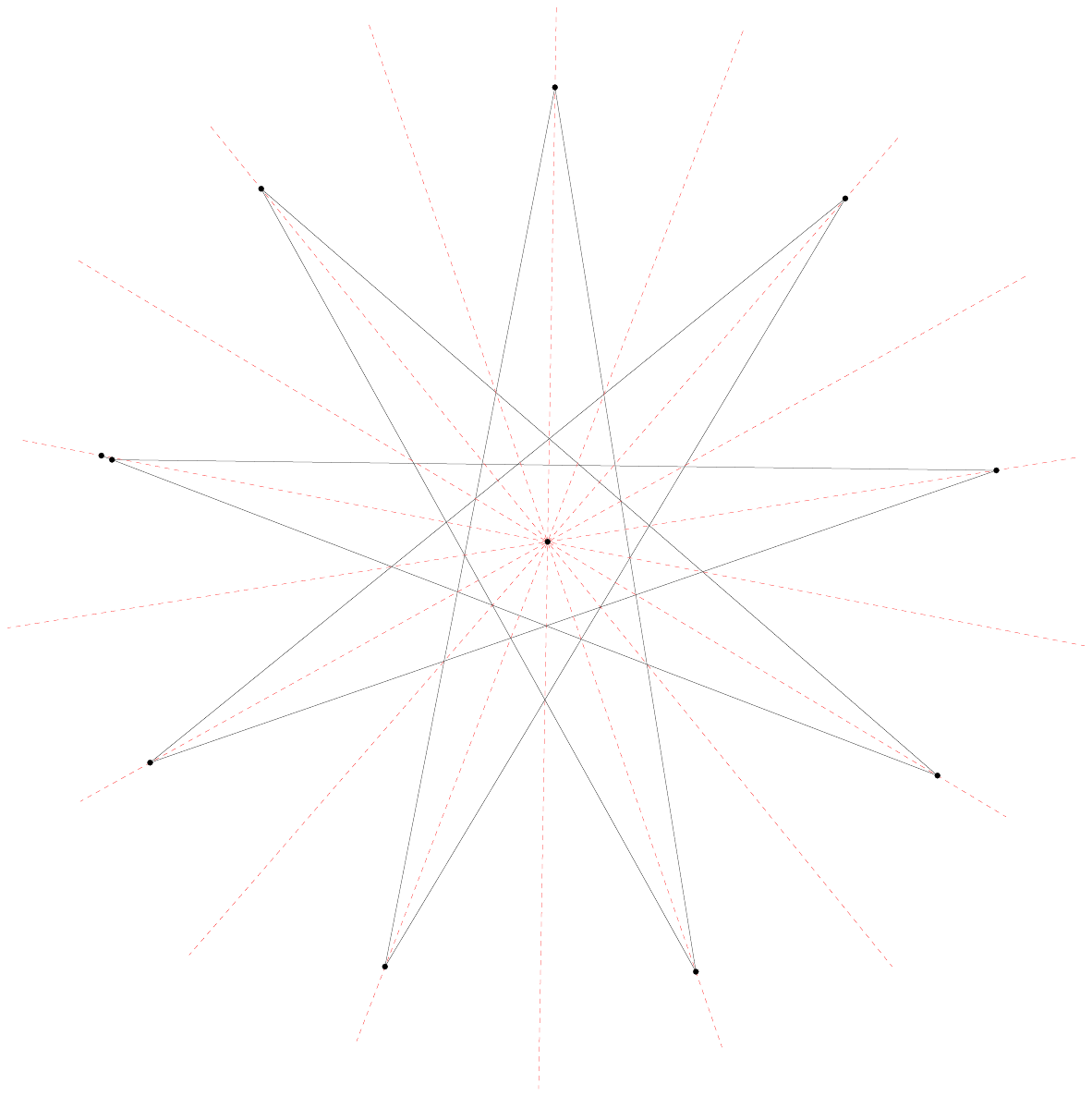}  \\ \hline
\end{tabular}
\caption{$9$-polygons with $9$ axes}
\end{figure}
\begin{figure}[H]
\begin{tabular}{|c | c | c |}
\hline
\includegraphics[width=0.3\textwidth]{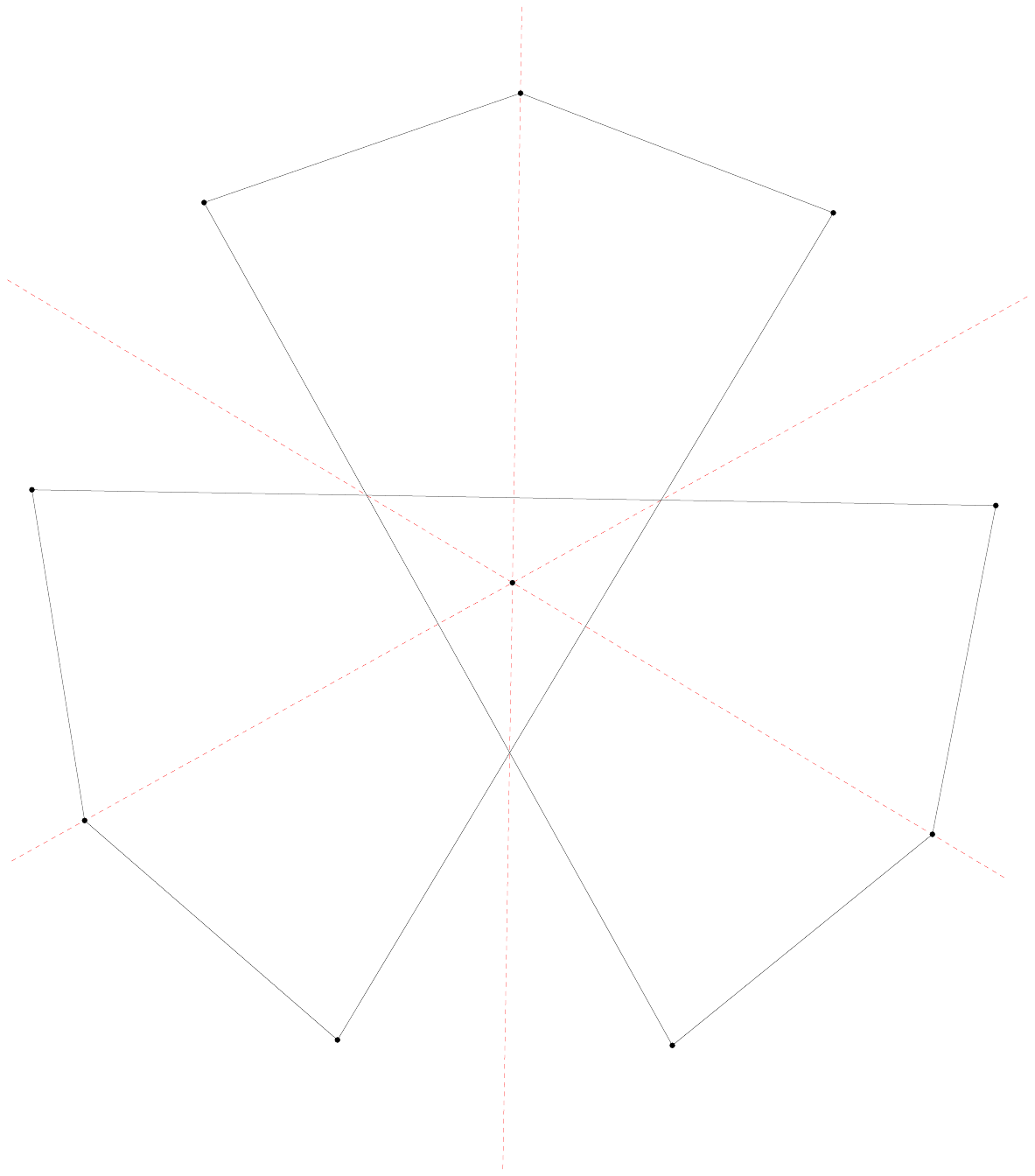} &
\includegraphics[width=0.3\textwidth]{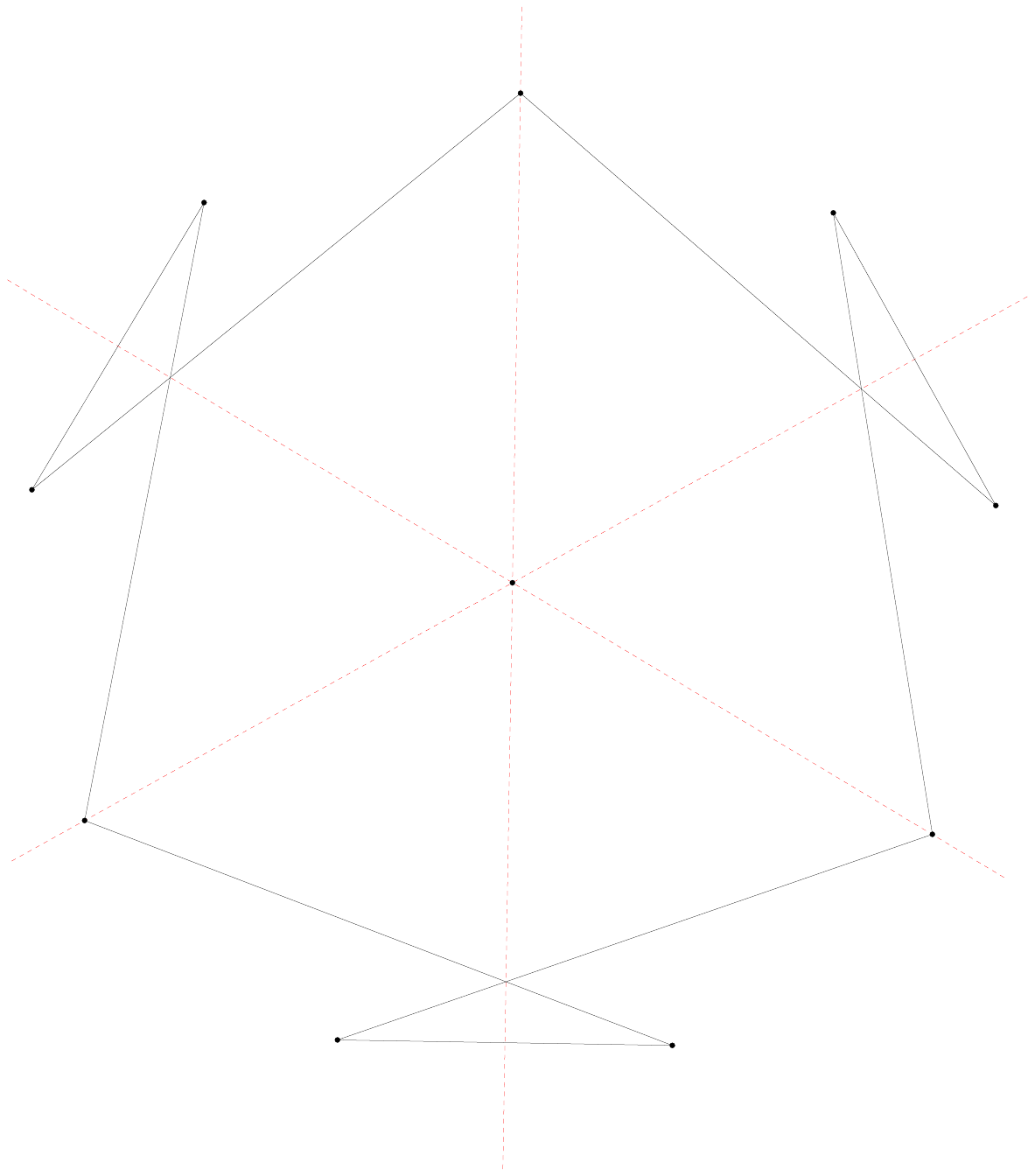} & \includegraphics[width=0.3\textwidth]{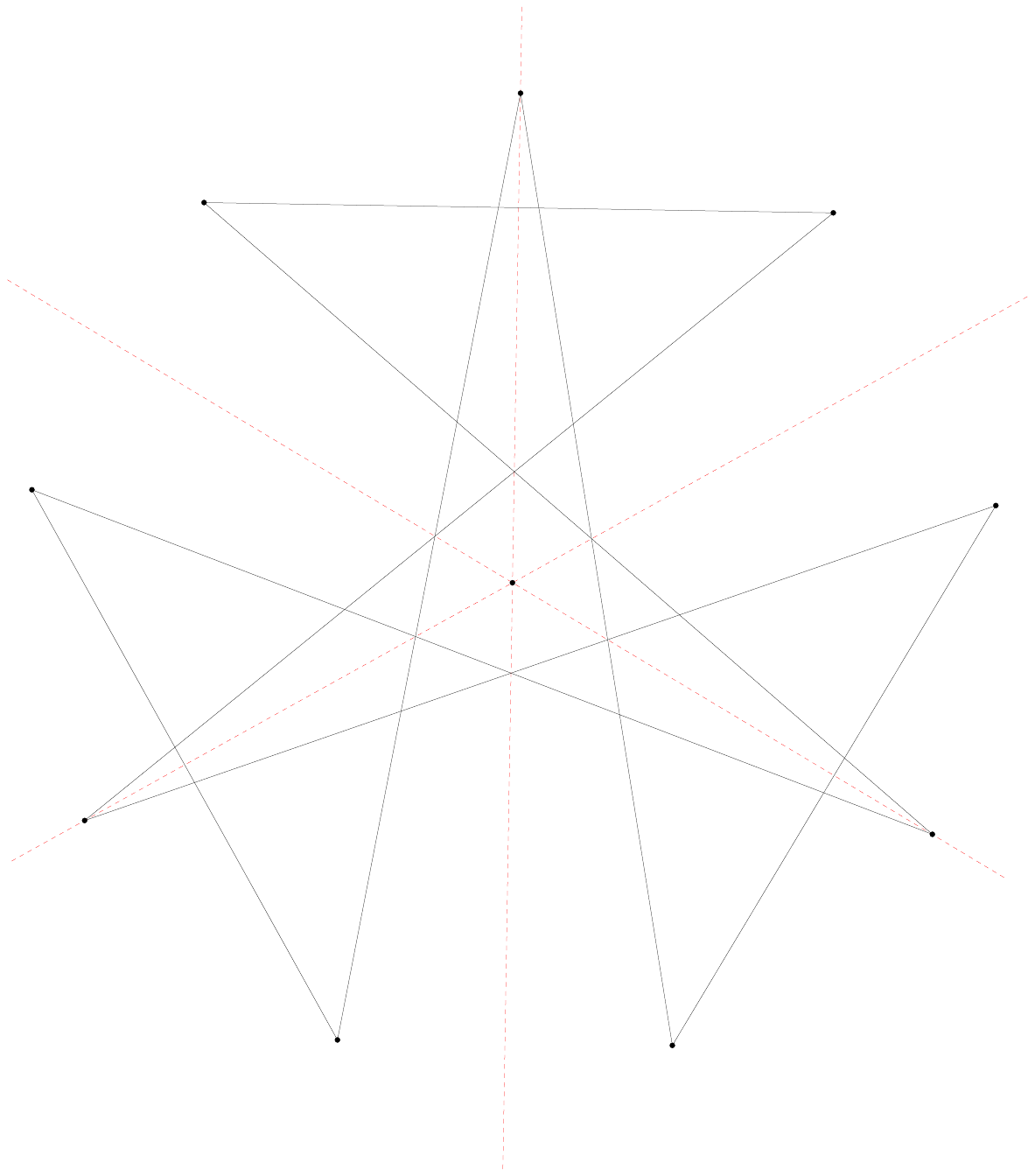}  \\ \hline
\end{tabular}
\caption{$9$-polygons with $3$ axes}
\end{figure}

\begin{figure}[H]
\begin{tabular}{| c | c |}
\hline
\includegraphics[width=0.5\textwidth]{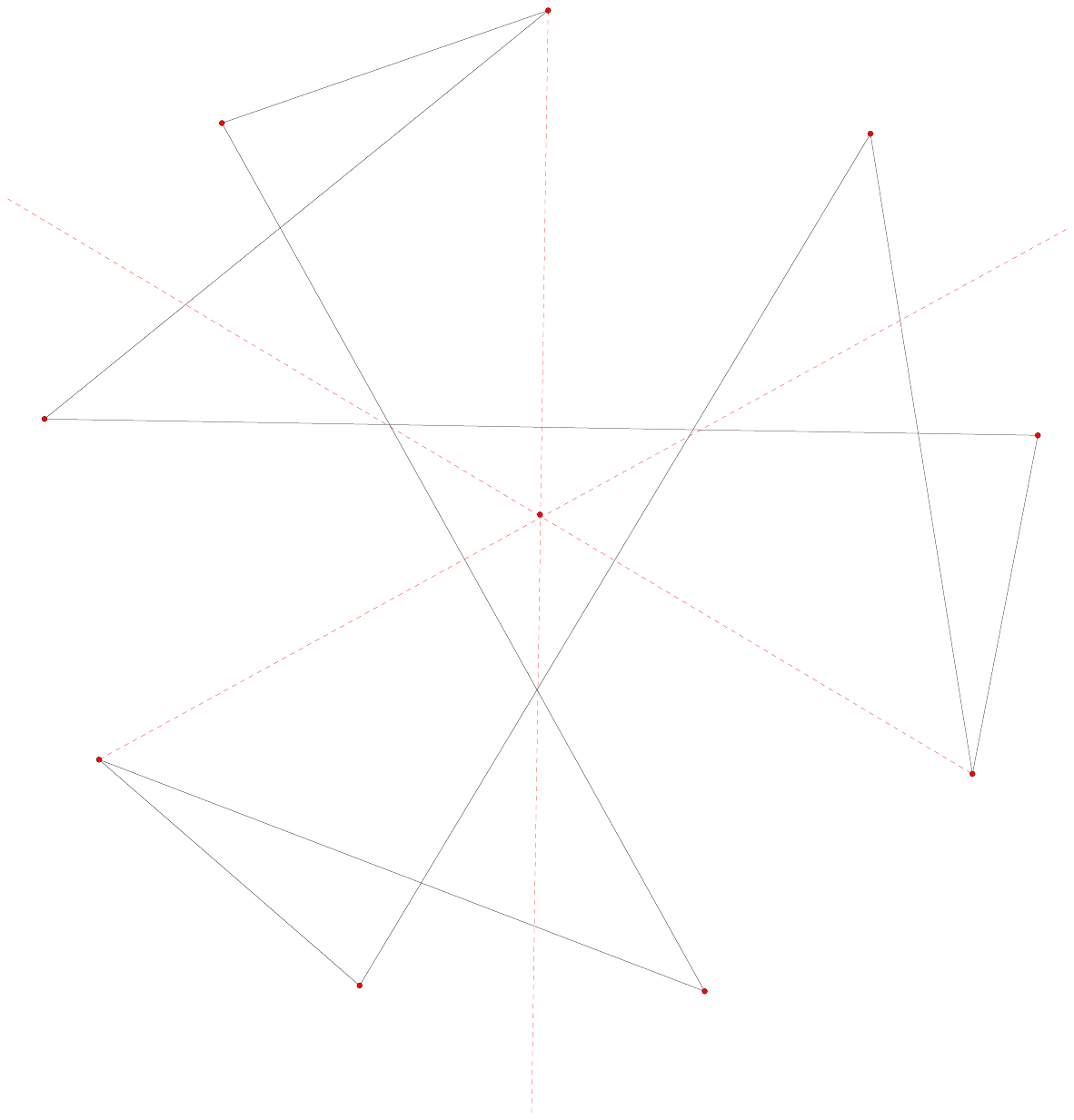} &
\includegraphics[width=0.5\textwidth]{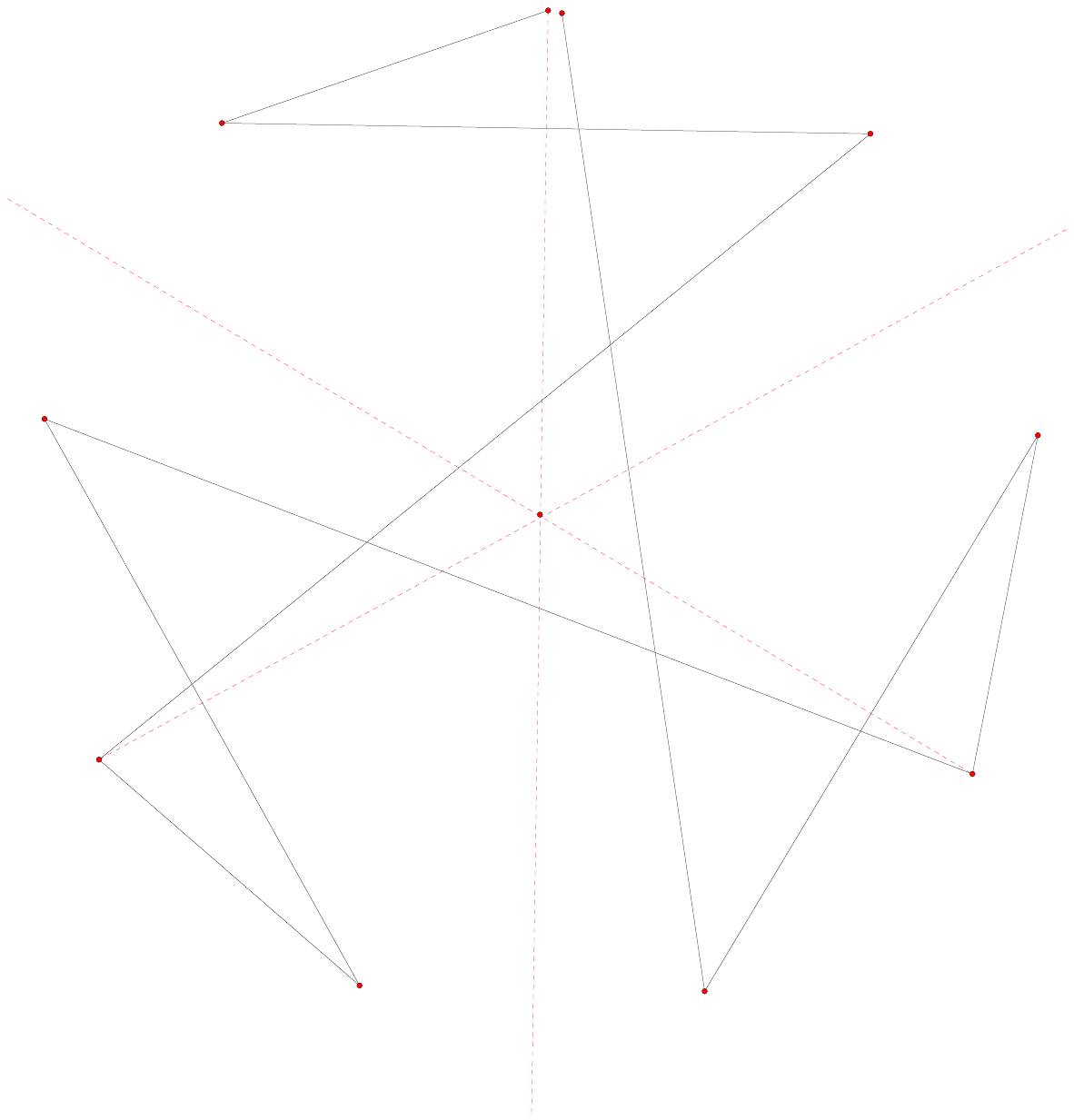} \\ \hline
\end{tabular}
\caption{$3$-circular $9$-polygons}
\end{figure}

Because $\vert X_1(9) \vert =186$ it is not possible to print a set of representatifs of the equivalence-classes of $9$-polygons with exactly $1$ axis. For $p>3$ there are to many equivalence-classes to make a visualisation.

\listoftables
\listoffigures
\bibliographystyle{Alpha}
\bibliography{p_und_pQuadrat_merged_for_Arxiv}
\end{document}